\documentclass[11pt,monochrome,reqno]{amsart}
\usepackage{
	amssymb,
	amsaddr
}
\newcommand{\no}[1]{#1}

\renewcommand{\no}[1]{}  \newcommand{\upDelta}{\Delta} 
\newcommand\norm[1]{\left\lVert#1\right\rVert}
\no{\usepackage{times}\usepackage[
	subscriptcorrection, slantedGreek,
	nofontinfo]{mtpro}
	\renewcommand{\Delta}{\upDelta}
}

\date{\today}
\usepackage{algpseudocode}
\usepackage{algorithm}
\usepackage{color}
\usepackage[pdftex]{graphicx}
\usepackage{subfigure}
\usepackage{epstopdf}
\usepackage{tikz}
\usepackage{xcolor}

\newtheorem{lemma}{Lemma}
\newtheorem{definition}{Definition}

\theoremstyle{remark}
\newtheorem{remark}{Remark}

\newcommand{\be}[1]{\begin{equation}\label{#1}}
\newcommand{\ee}{\end{equation}}

\newcommand{\Scom}[1]{\textcolor{black}{#1}}

\title[Direct QPAT for realistic acoustic media]{Direct quantitative photoacoustic tomography for realistic acoustic media}	
\author[A Javaherian]{Ashkan Javaherian$^{1,2}$}
\author[S Holman]{Sean Holman$^1$} 

\address[1]{$^1$School of Mathematics, University of Manchester, M13 9PL, UK\\ e-mail: sean.holman@manchester.ac.uk}

\address[2]{$^2$Department of Medical Physics and Biomedical Engineering, University College London, WC1E 6BT,UK\\ e-mail: a.javaherian@ucl.ac.uk}

\date{}
\thanks{}

\begin{document}

\maketitle
\begin{abstract}
Quantitative photo-acoustic tomography (QPAT) seeks to reconstruct a distribution of optical attenuation coefficients inside a sample from a set of time series of pressure data that is measured outside the sample. The associated inverse problems involve two steps, namely acoustic and optical, which can be solved separately or as a direct composite problem. We adopt the latter approach for realistic acoustic media that possess heterogeneous and often not accurately known distributions for sound speed and ambient density, as well as an attenuation following a frequency power law that is evident in tissue media. We use a Diffusion Approximation (DA) model for the optical portion of the problem. We solve the corresponding composite inverse problem using three {\textcolor{red}{total variation (TV) regularised}} optimisation approaches. Accordingly, we develop two Krylov-subspace inexact-Newton algorithms that utilise the Jacobian matrix in a matrix-free manner in order to handle the computational cost. Additionally, we use a gradient-based algorithm that computes a search direction using the L-BFGS method, and applies a TV regularisation based on the Alternating Direction Method of Multipliers (ADMM) as a benchmark, because this method is popular for QPAT and direct QPAT. The results indicate the superiority of the developed inexact Newton algorithms over gradient-based Quasi-Newton approaches for a comparable computational complexity.
\end{abstract}

\section{Introduction}\label{Intro}

Quantitative photo-acoustic Tomography (QPAT) is a steadily growing hybrid imaging paradigm that simultaneously takes advantage of the high spatial resolution provided by ultrasound imaging and the rich contrast attributed to optical imaging \cite{Wang-a}. Typically, nanosecond-duration pulses of electromagnetic energy, in the visible or near-infrared ranges, are used to irridiate a sample. Depending on the optical properties of the sample, a fraction of the optical energy is absorbed, and converted into heat \cite{Wang}. The generated heat induces a local increase in pressure via thermo-elastic expansion effects \cite{Wang}. Because of the elasticity of soft tissues, the locally induced pressures propagate outwards as acoustic waves, and carry information about the optical properties of the sample to the surface. These acoustic waves are measured in time by ultrasound sensors located outside the surface of the sample \cite{Wang-a}. Given a set of time series of data at the boundary, the objective in QPAT is to calculate a quantitative image of a distribution of optical absorption coefficients of the sample \cite{Wang-a}.

QPAT involves two distinct inverse problems, namely acoustic and optical. The acoustic inverse problem, often referred to as PAT, seeks to reconstruct a distribution of the spatially varying initial pressure from the boundary data. This is a linear inverse problem, for which a vast number of reconstruction methods are available \cite{Rosenthal}. For acoustically homogeneous media, there are exact inversion methods, e.g., back-projection \cite{Xu-a,Xu-d,Burgholzer,Kunyansky}, or frequency-domain techniques \cite{Xu-b,Xu-c}. Time reversal (TR) is a less restrictive approach since it can be adapted to heterogeneous acoustic media and arbitrary detection surfaces \cite{Finch,Xu-gh,Burgholzer,Hristova-b,Hristova-a}. Model-based iterative inversion approaches, in which the discrepancy between modeled data and the measured data is iteratively reduced, are also widely used in ill-posed cases that arise from data incompleteness, modelling errors, finite sampling or noise \cite{Huang}. The iterative methods can be categorised into convergent Neumann-series methods \cite{Stefanov-b,Qian} or variational (optimisation) approaches \cite{Wang-c,Huang,Ding-g,Belhachmi,Haltmeier-d}.

The  optical inverse problem is devoted to further reconstruction of an image of distributions of optical absorption and scattering coefficients from a recovered initial pressure distribution. This is a highly nonlinear and ill-posed inverse problem, and is commonly solved by iterative model-based approaches \cite{Tarvainen,Tarvainen-b,Gao}. Using model-based approaches for solving these two inverse problems, an accurate, yet efficient, modelling of the associated physical processes is required \cite{tarvainen-n}.

A very accurate model for propagation of light is the Radiative Transfer Equation (RTE), which has been used for the optical portion of 2D QPAT, e.g. \cite{Tarvainen,Saratoon,Haltmeier}. An analysis of the optical inverse problem of QPAT using RTE has been given in \cite{Bal-b}. Since RTE is computationally very expensive, the applicability of this method for 3D QPAT is very limited. A Diffusion Approximation (DA) of the RTE is more efficient than RTE \cite{tarvainen-n}, and is thus more practical for 3D QPAT. In this study, we use the DA for modelling the optical portion of the inverse problem. The DA model is sufficienly accurate when the medium is highly scattering and the scattering is near-isotropic \cite{tarvainen-n}. The DA has been widely used as a light propagation model for QPAT \cite{Tarvainen,Tarvainen-b,Pulkkinen-a,Gao,Hannukainen}.

The DA model is defined as a function of optical absorption and scattering coefficients. It turns out that the recovery of both these coefficients from one optical source and wavelength is non-unique \cite{Bal-a}. To have uniqueness for the inversion, three approaches are used, i.e., the assumption of the scattering coefficient as known \cite{Cox-b}, using more than one optical wavelength \cite{Bal-d,Cox-c,Pulkkinen-a}, or using more than one optical source \cite{Bal-c,Gao,Gao-b}. In this study the latter approach is used.

For modelling the acoustic portion of the problem, the dependence of shape, spectrum and amplitude of propagating acoustic waves on properties of the medium \cite{Cox}, together with the highly nonlinear and ill-conditioned nature of the optical portion of the forward operator, motivates enriching the QPAT problem by simulation of tissue-realistic acoustic properties. To incorporate these effects, we use an acoustic model based on a linear system of three-coupled first-order wave equations which can be adapted to spatially varying sound speed and density \cite{Tabei,Cox}, and include two fractional Laplacian operators in order to account separately for acoustic absorption and dispersion following a frequency power law, which is evident in tissue media \cite{Treeby,Treeby-a}. For a numerical implementation of this acoustic model, we use a k-space pseudo-spectral method, which is popular for iterative PAT because of the high efficiency arising from a requirement of only two grid points per wave-length for defining a field, and a fast computation of the spatial gradient in the frequency domain \cite{Tabei,Cox}.

For ill-posed problems of QPAT, for example, when data are available merely on a part of the boundary, or when acoustic properties of the medium are not known exactly, the reconstruction of initial pressure distribution using the acoustic portion of the inverse problem is not sufficiently accurate to further serve as data for the optical inverse problem \cite{Gao-b}. It turns out that a direct estimation of the optical coefficients from time series of measured data, referred to here as direct QPAT, is more stable than classical variants of QPAT, in which the two inverse problems are solved distinctly. The direct QPAT has thus recently received much attention, where the forward operator is defined using a composite opto-acoustic model \cite{Haltmeier,Gao-b,Pulkkinen-b}. Using direct QPAT,  it will be possible to incorporate prior information about the optical parameters into the acoustic inverse problem \cite{Gao-b}, and mutually the optical inverse problem can utilise information about noise included in the boundary acoustic data \cite{Pulkkinen-b}. {\textcolor{red}{Additionally, the acoustic portion of the inverse problem can benefit from multi-source (resp. multi-wavelength) configurations, since the optical coefficients (resp. chromophore concentrations), as opposed to the initial pressure distribution, are independent from the changes in optical illumination (resp. wavelength) \cite{Ding}.}}

In \cite{Shao,Song}, a linear Born approximation of the DA model based on a Green's function approach is coupled with an acoustic model that uses a free-space Green's function method in order to directly reconstruct perturbations in the optical absorption and diffusion coefficients from time series of measured acoustic data. A simultaneous reconstruction of perturbations in optical absorption coefficient and sound speed given an optical scattering coefficient was also studied using a Born approximation \cite{Ding}. To the best of our knowledge, existing studies for the direct variant of QPAT have been so far limited to homogeneous and non-attenuating acoustic media, and the acoustic portion of the forward operator and its adjoint are computed based on exact formulae using Green's function approaches \cite{Haltmeier,Gao-b,Pulkkinen-b}.

The optimisation approaches for QPAT can be categorised into linearised (Jacobian-based) or nonlinear (gradient-based) approaches. A majority of Jacobian-based methods for the optical portion of QPAT can be fit into a class of Gauss-Newton methods \cite{Schweiger}. For application of these methods to classical (resp. direct) QPAT, see  \cite{Tarvainen,Tarvainen-b,Pulkkinen-a} (resp. \cite{Pulkkinen-b}). These studies utilise an explicit form of the Jacobian matrix. For direct QPAT, because of a very large size of time series of measured data, an explicit computation and storage of the Jacobian matrix is very expensive \cite{Pulkkinen-b}.

To avoid this problem, nonlinear gradient-based approaches have received much interest for the direct problem of QPAT \cite{Gao,Gao-b}. A majority of these approaches use a search direction based on Quasi-Newton methods which utilise only gradient information for an approximation of the Hessian matrix, and are thus memory-efficient \cite{Gao,Gao-b}. The computation of the gradient for direct QPAT is based on an \textit{opto-acoustic} forward operator and an \textit{acousto-optic} adjoint of the Fr\'echet derivative of the forward operator \cite{Haltmeier}. A memory-efficient Quasi-Newton method is Limited-memory BFGS (L-BFGS), for which the inverse of Hessian matrix is approximated using the gradient information stored in a user-adjusted number of iterations. Using a total variation (TV) regularisation approach based on an Alternating Direction Method of Multipliers (ADMM), L-BFGS has been used for computing an associated search direction for the classical QPAT \cite{Gao} and direct QPAT \cite{Gao-b}.

\textbf{Contribution.} Here, we develop two preconditioned inexact Newton (Newton-Krylov) algorithms for solving the direct problem of QPAT. In our first approach, the residual function is iteratively linearised, and each linearised subproblem is solved using a \textit{subspace Krylov} method in a matrix-free manner, for which a total variation matrix is used as a preconditioner for accelerating the convergence of the Krylov method. We use the \textit{Preconditioned Conjugate Gradient} (PCG) as the \textit{Krylov} method, for which a TV-based preconditioning is applied using the \textit{Lagged Diffusivity} (LD) method \cite{Vogel}. Our second approach uses two linearisations, the first of which is enforced to the nonlinear residual function, and the second is applied in order to handle the nonlinearity of a TV functional using a Primal-Dual Interior Point Method (PD-IPM). Using this approach, for each linearisation of the residual function, a sequence of normal equations is derived, and is solved using a \textit{subspace Krylov} method in a matrix-free manner. We solve the arising normal equations using a PCG method, but we emphasise that an extension to other \textit{Krylov} methods, e.g. a preconditioned variant of LSQR \cite{Arridge-c}, is straightforward.

We implement our algorithms by assuming tissue realistic but erroneous properties for the acoustic medium with a limited-view setting for boundary measurements. We model the acoustic portion of the forward operator using a linear system of three-coupled first-order wave equations that can be adapted to heterogeneous media and account for acoustic absorption and dispersion following a frequency power law using two fractional Laplacians \cite{Treeby-a}. This acoustic model is very popular for PAT since it simulates an acoustic attenuation that is evident in many materials of interest including tissues \cite{Treeby-a}. To the best of our knowledge, this is the first study for direct QPAT that uses an acoustic model that accounts for tissue-realistic acoustic properties of the medium. For the acoustic portion of the problem, we include the action of perfectly matched layers (PMLs) in calculation of the adjoint operator. To the best of our knowledge, these effects have not been accounted for in existing studies (See \cite{Huang}).

It is worth mentioning that a singular value decomposition (SVD) analysis on the acoustic forward operator we use has shown that as time steps increase some of the singular values of the forward operator become very small \cite{Treeby}, and make the acoustic forward operator ill-conditioned. However, the use of an acoustically realistic forward operator may be necessary for direct QPAT, since the optical portion of the forward operator is highly nonlinear and ill-conditioned, and thus errors in acoustic modelling may quickly grow in the inversion process, and dominate signal data.

Our numerical results show that the developed preconditioned Newton-Krylov optimisation algorithms perform much better than nonlinear Quasi-Newton methods that have been used in existing studies for direct QPAT. The algorithm we use as a benchmark utilises a TV regularisation based on an ADMM method, together with a search direction using an L-BFGS method. (See \cite{Gao, Gao-b} for application on classical QPAT and direct QPAT, respectively).

\section{Direct QPAT on a continuous domain}
\noindent
In this section, we define our forward operator as a composite map on a continuous domain.
\subsection{Modelling the optical portion of the problem}
The time scale for propagation of acoustic waves is on the order of a micro-second, which is three orders of magnitude larger than the time-scale for illumination, propagation and absorption of light. Therefore, the generated pressure distribution is regarded as instantaneous for the acoustic problem, and is referred to as initial pressure distribution $p_0$. One way to define a forward operator for our QPAT problem is to combine the physics of the optical and acoustic portions of the problem using a simple composition of two maps, one modelling propagation and absorption of optical photons and the other modelling propagation of acoustic waves \cite{Haltmeier}.

Accordingly, let $\Omega\subset\mathbb{R}^d$ be a convex bounded domain with Lipschitz boundary $\partial \Omega$ and $d \in \left\{2,3\right\}$ the spatial dimension. Additionally, let $\phi \in \mathcal{H}^1(\Omega)$ denote the photon density. For modelling the propagation of light, we use a time-independent variant of DA equations with the well-known Robin boundary condition \cite{Tarvainen}. This is written as
\begin{align} \label{DA}
	\begin{split}
		-\nabla \cdot \kappa(r) \nabla \phi (r)+ \mu(r) \phi(r)= 0,& \quad r\in \Omega     \\
		\phi(r)+\frac{1}{2 \gamma_d}  \kappa(r) \frac{\partial \phi(r)}{\partial \hat{n}}=
		I/\gamma_d,& \quad r\in \partial \Omega,
	\end{split}
\end{align}
where $\mu(r),\kappa(r) \in \mathcal{L}_+^\infty (\Omega)$ denote the positive-valued optical absorption and diffusion coefficients, respectively. Here, $r$ denotes the spatial position. Additionally, $\gamma_d$ is a dimension-dependent factor ($\gamma_2=1/\pi$ and $\gamma_3=1/4$), $\hat{n}$ is an outward unit normal, and $I$ is an inward directed diffuse boundary current \cite{Schweiger-b,tarvainen-n}. Following the absorption of photons, a spatially varying heating field $h(r)$ is generated in the form
\begin{align} \label{h}
	h [\kappa,\mu,I](r)=\mu(r) \phi[\kappa,\mu,I](r).
\end{align}
Because of thermo-elastic expansion effects, the induced spatially varying heating field causes an instant local increase in pressure that follows
\begin{align}
	p_0(r)=
	\begin{split}
		\begin{cases}
			\Gamma(r) h(r), \quad &r\in\Omega\\
			0, \quad &r\in \mathbb{R}^d\backslash\Omega,
		\end{cases}
	\end{split}
\end{align}
where $\Gamma$ is the Gruneisen parameter, and describes the efficiency of conversion of heat into pressure \cite{Haltmeier}. Here, we assume $\Gamma(r)$ constant and rescaled to $1$, and thus {\textcolor{red}{$p_0(r)=h(r)$}} \cite{Haltmeier}.

\subsection{Modelling the acoustic portion of the problem}   \label{sda}
We use a linear system of three-coupled first-order equations for describing the propagation of acoustic wavefields in an acoustically heterogeneous and lossy medium \cite{Treeby,Treeby-a}. To explain this, we define our fields, which are the acoustic pressure field $p(r,t)$, particle velocity vector $v(r,t)$ and acoustic density $\rho(r,t)$, where $r\in \mathbb{R}^d$ and $t \in [0,\infty)$ denote the spatial position and time. Additionally, we define the medium's acoustic parameters as sound speed $c_0(r)$, ambient density $\rho_0(r)$, attenuation coefficient $\alpha_0(r)$, and frequency power law exponent $y$. The acoustic wavefield propagation is now modeled by three equations, i.e., linearised equation of motion (conservation of momentum)
\begin{equation}\label{mom}
	\frac{\partial v}{\partial t}(r,t) = - \frac{1}{\rho_0(r)} \nabla p(r,t),
\end{equation}
linearised equation of continuity (conservation of mass)
\begin{equation} \label{mass}
	\frac{\partial \rho}{\partial t}(r,t) = - \rho_0(r) \nabla \cdot v(r,t),
\end{equation}
and equation of state
\begin{equation}\label{state}
	\begin{split}
		p(r,t) = c_0(r)^2 \Big \{& 1 - \tau(r) \frac{\partial}{\partial t} (-\nabla^2)^{\frac{y}{2}-1}\\
		&- \eta(r) (-\nabla^2)^{\frac{y-1}{2}} \Big \} \rho (r,t)
	\end{split}
\end{equation}
with initial conditions
\begin{equation} \label{initial}
	p(r,0) = p_0(r), \quad v(r,0) = 0.
\end{equation}
Here, $\tau(r)$ and $\eta(r)$ are, respectively the absorption and dispersion proportionality coefficients, and are given by
\begin{align}
	&\tau(r)=-2\alpha_0 c_0(r)^{y-1}, \  \ \eta(r)=2\alpha_0 c_0(r)^y \tan(\pi y/2).     \label{Twelve}
\end{align}

\subsection{Opto-acoustic forward operator} \label{skp}
Having given the models for describing the optical and acoustic portions of the forward operator, we now define our \textit{opto-acoustic} forward operator. To do this, we also require the measurement operator $\mathcal{M}$, which at each time step maps the pressure $p(r,t)$ to the measured data at sensors. The time series of acoustic data are denoted by $\hat{P} \in \mathbb{R}^{N_s N_t}$ with $N_s,N_t \in \mathbb{N}$ the number of detectors and the number of measurement time instants, respectively.

\begin{definition}
	For a fixed illumination $I$, the optical portion of the forward operator is a nonlinear map in the form
	\begin{align}   \label{opt_map}
		\begin{split}
			&\Lambda_o:\mathcal{L}_+^\infty(\Omega) \times \mathcal{L}_+^\infty(\Omega) \rightarrow \mathcal{L}^2(\Omega)  \\
			&{\Lambda_o}[\kappa,\mu](r)= h(r),
		\end{split}
	\end{align}
	where $h(r)$ is given by \eqref{h}.
	Also, the acoustic portion of  the forward operator is a linear map in the form \cite{Arridge}
	\begin{equation}   \label{acoust_map}
		\begin{split}
			&\Lambda_a: \mathcal{L}^2(\Omega)   \rightarrow \mathbb {R}^{N_s N_t},\\
			&\Lambda_a [p_0](r,t) =  \mathcal{M} p(r,t).
		\end{split}
	\end{equation}
	Using \eqref{opt_map} and \eqref{acoust_map}, the coupled opto-acoustic forward operator is defined by the composite map $\Lambda$, i.e.,
	\begin{equation}
		\begin{split}
			&\Lambda:  \mathcal{L}_+^\infty(\Omega) \times \mathcal{L}_+^\infty(\Omega)   \rightarrow \mathbb {R}^{N_s N_t},\\
			&\Lambda [\kappa,\mu]                                                              = \Lambda_a \big[ \Lambda_o[\kappa,\mu] \big] .
		\end{split}
	\end{equation}
\end{definition}

\subsection{Model-based approach for inverse problem}
The inverse problem of QPAT is a simultaneous reconstruction of $\mu, \kappa$ from $\hat{P} = \Lambda[\kappa,\mu]$. Applying a model-based approach for solving this problem, the objective is to minimise an error functional, the sum of squared differences between modeled data and the measured data via an iterative adjustment of the unknown optical coefficients $\mu, \kappa$. Using $N_q$ optical illuminations, the error functional is defined by \cite{Gao-b}
\begin{align} \label{obj_couple}
	\epsilon (\kappa,\mu)= \frac{1}{2} \sum_{q=1}^{N_q} \norm{\Lambda_q[\kappa,\mu]- \hat{P}_q }^2 ,
\end{align}
where, $q$ indexes a set of $N_q$ illuminations $I_q$ with corresponding acoustic data $\hat{P}_q$ \cite{Gao-b}. We combine the unknown parameters as $x=[\kappa,\mu] \in \mathcal{L}_+^\infty(\Omega) \times \mathcal{L}_+^\infty(\Omega)$.

Let $D_x \Lambda$ denote the Fr\'echet derivative of the forward operator at $x$. The Fr\'echet derivative of $\epsilon$ \Scom{at} $x$ is given by
\begin{align} \label{fre}
	D_x \epsilon =\sum_{q=1}^{N_q} D_x^* \Lambda_q \left(\Lambda_q[x]-\hat{P}_q \right).
\end{align}
Here, $D_x^* \Lambda$ denotes the adjoint of the Fr\'echet derivative of the forward operator, and is given by
\begin{align}
	D_x^* \Lambda [\hat{P}]= D_x^* \Lambda_o [\Lambda_a^*[\hat{P}]],
\end{align}
where $\Lambda_a^*$ denotes the adjoint of the linear acoustic forward operator $\Lambda_a$ \eqref{acoust_map} and $D_x^* \Lambda_o$ represents the adjoint of the Fr\'echet derivative of the optical forward operator \eqref{opt_map} both with respect to the $\mathcal{L}^2(\Omega)$ inner product. Formulae for $\Lambda_a^*$ can be found in \cite{Javaherian}, while for $D_x^* \Lambda$ we have the next lemma.

\begin{lemma} \label{ltwo}
	Let us denote the solution of \eqref{DA} for the fixed illumination $I$, diffusion $\kappa_0$ and absorption $\mu_0$ by $\phi_0$. Then the Fr\'echet derivative of the optical portion of the forward operator $D_x \Lambda_o$ at $x_0=[\kappa_0,\mu_0]$ applied to the perturbations $\delta \kappa$ and $\delta \mu$ is given by
	\begin{align}  \label{fr}
		D[\kappa_0,\mu_0]{\Lambda_o}  \left (
		\begin{array}{c}
			\delta \kappa\\
			\delta \mu
		\end{array}
		\right ) = \delta \mu(r) \phi_0(r) + \mu_0(r) \delta \phi(r),
	\end{align}
	where $\delta \phi(r)$ satisfies
	\begin{align} \label{dPhi}
		\begin{split}
			- \nabla \cdot \kappa_0(r) \nabla \delta \phi(r) + \mu_0(r) \delta \phi(r) = \nabla \cdot \delta \kappa(r) \nabla \phi_0(r) - \delta \mu(r) \phi_0(r), \quad &r \in \Omega\\
			\delta \phi(r) + \frac{1}{2 \gamma_d}  \kappa_0(r) \frac{\partial \delta \phi(r)}{\partial \hat{n}}   = - \frac{1}{2 \gamma_d} \delta \kappa(r) \frac{\partial \phi_0(r)}{\partial \hat{n}}, \quad &r \in \partial \Omega.
		\end{split}
	\end{align}
	
	The adjoint map $D_x^* \Lambda_o$ can then be calculated from
	\begin{equation}   \label{opt_adjoint}
		D[\kappa_0,\mu_0]^*{\Lambda_o} h(r) = \left (
		\begin{array}{c}
			\nabla \phi_0(r) \cdot \nabla \tilde{h}(r) \\
			\phi_0(r) \tilde{h}(r)+ \phi_0(r) h(r)
		\end{array}
		\right ),
	\end{equation}
	where the adjoint field $\tilde{h}(r)$ satisfies
	\begin{equation} \label{th}
		\begin{split}
			- \nabla \cdot \kappa_0(r) \nabla \tilde{h}(r) + \mu_0(r) \tilde{h}(r) = -\mu_0(r) h(r), \quad     &r\in  \Omega \\
			\tilde{h}(r) + \frac{1}{2 \gamma_d}  \kappa_0(r) \frac{\partial \tilde{h}(r)}{\partial \hat{n}}=0,\quad     &r\in {\partial \Omega}.
		\end{split}
	\end{equation}
\end{lemma}

\noindent Lemma \ref{ltwo} can be proven using integration by parts.

\section{Numerical computation}
\subsection{Numerical computation of the optical operators ($ \Lambda_o$, $ D_x \Lambda_o$ and $D_x^* \Lambda_o$ )}
\subsubsection{Variational formulae}     \label{vari}
We use a first-order Galerkin finite element method (FEM) for approximation of the optical portion of our QPAT problem. For an approximation of $\Lambda_o$, a variational form of \eqref{DA} is derived, i.e.,
\begin{align}  \label{var_fo}
	\int_\Omega \mu_0 \phi_0 \nu\ dr +\int_\Omega \kappa_0   \nabla\phi_0 \cdot \nabla \nu\ dr + 2 \gamma_d  \int_{\partial\Omega}  \phi_0 \nu\ ds =  \int_{\partial\Omega} 2 I_s \nu\  ds,
\end{align}
where $\nu \in \mathcal{H}^1 (\Omega)$ is a test function. Additionally,
for an approximation of the Fr\'echet derivative operator $D_x \Lambda_o$ using FEM, a variational form of \eqref{dPhi} is derived, i.e.,
\begin{align}  \label{var_fr}
	\begin{split}
		\int_\Omega \mu_0 \delta \phi \nu dr+  \int_\Omega  \kappa_0  \nabla \delta \phi \cdot \nabla \nu dr+  2 \gamma_d \int_{\partial\Omega} \delta \phi \nu ds
		=  - \int_\Omega \delta \kappa \nabla \phi_0  \cdot \nabla\nu dr   - \int_\Omega  \delta \mu \phi_0 \nu dr.
	\end{split}
\end{align}
In the same way, for an approximation of $D_x^* \Lambda_o$, a variational form of \eqref{th} is obtained, i.e.,
\begin{align}  \label{var_f_adj}
	\int_\Omega \mu_0\tilde{h} \nu dr +\int_\Omega \kappa_0 \nabla\tilde{h} \cdot \nabla \nu dr+  2 \gamma_d  \int_{\partial\Omega} \tilde{h} \nu ds
	=-\int \mu_0  h \nu dr.
\end{align}

\subsubsection{Discretisation of optical coefficients and fields}        \label{dis_forward}
Let $T$ denote a triangulation of $\Omega$ with $N_e$ elements, i.e., $T=\big\{t_j\:|\: j=1,...,N_e\big\}$.
Applying an approximation using FEM, we discretise the optical coefficients in a piecewise-constant basis $\{\chi_j=1_{t_j}\:|\: j=1,\,...\,,N_e\}$. Using this, the optical fields are approximated as \cite{Tarvainen}
\begin{align}
	\begin{split}
		\kappa_0(r) \approx\kappa^e(r) =\sum_{j=1}^{N_e} \hat{\kappa}_j \chi_j(r)\\
		\mu_0(r) \approx \mu^e(r)  =\sum_{j=1}^{N_e} \hat{\mu}_j \chi_j(r),
	\end{split}
	\label{elemental}
\end{align}
where $\hat{\kappa}_j$ and $\hat{\mu}_j$ denote the discretised absorption and diffusion coefficients at element $t_j$.
Additionally, $\phi_0(r)$ is approximated in a piecewise-linear basis $\{\varphi_k \:|\: k=1,\,...\,,N_n\}$ in the form
\begin{align} \label{piecewise}
	\phi_0(r) \approx \phi_0^h(r)= \sum_{k=1}^{N_n} \Phi_{0,k} \varphi_k(r),
\end{align}
where $\Phi_{0,k}$ denotes the discretised photon density at node $k$, and $N_n$ is the total number of nodes.
We also approximate the adjoint field $\tilde{h}(r)$ in a piecewise-linear basis \Scom{as}
\begin{align} \label{piecewise2}
	\tilde{h}(r) \approx \tilde{h}^h(r)= \sum_{k=1}^{N_n} \tilde{H}_k \varphi_k(r).
\end{align}
\Scom{In the sequel, a field that is discretised at nodes as in \eqref{piecewise} and \eqref{piecewise2} (resp. elements as in \eqref{elemental}) is called a nodal (resp. elemental) vector.}
In the same way as the continuous formulae, we use $\delta$ for signifying a perturbation in a discretised coefficient or field.

\subsubsection{Matrix form of variational formulae}
For a discretisation of the problem, a matrix form of the variational formulae in section \ref{vari} is derived. \cite{tarvainen-n,Tarvainen}.
To do this, we define a system matrix $\boldsymbol{A}_o$ in the form
\begin{align} \label{system_forward2}
	\boldsymbol{A}_o[\kappa^e,\mu^e] = K[\kappa^e] + C[\mu^e] + R
\end{align}
with
\begin{align}
	\begin{split}
		&K_{kp}[\kappa^e] = \int_\Omega  \kappa^e \nabla\varphi_k\cdot \nabla \varphi_p \ dr  \\
		&C_{kp}[\mu^e]=  \int_\Omega \mu^e \varphi_k \varphi_p  \ dr  \\
		&R_{kp}= 2 \gamma_d \int_{\partial\Omega}  \varphi_k \varphi_p \ ds\\
		&G_p=\int_{\partial\Omega}2 I_s\varphi_p \ ds,
	\end{split}
\end{align}
where $p,k=1,...,N_n$ denote \Scom{nodal indices}.
We also define matrix $\boldsymbol{A}_{\delta,o}$ in the form
\begin{align}
	\boldsymbol{A}_{\delta,o}[\delta \kappa^e,\delta \mu^e]=K[\delta \kappa^e]+C[\delta \mu^e].
\end{align}
Using the above, we now define the discretised optical forward operators. We stress that this definition is setting the notation for the discretised operators, which will be described in detail below.
\begin{definition}
	A discretisation of the optical forward operator $\Lambda_o$  gives a map from a vector space of discretised (elemental) optical coefficients to a vector space of discretised (elemental) heating field coefficients in the form
	\begin{align}   \label{opt_des}
		\begin{split}
			&\mathbb{H}_o:\mathbb{R}^{N_e} \times \mathbb{R}^{N_e}  \rightarrow  \mathbb{R}^{N_e}\\
			&H=\mathbb{H}_o[\hat{\kappa},\hat{\mu}].
		\end{split}
	\end{align}
	Additionally, a discretisation of the Fr\'echet derivative operator $D_x \Lambda_o$ at $X=[\hat{\kappa},\hat{\mu}]$ applied on perturbation  $\delta X=[\delta \hat{\kappa},\delta \hat{\mu}]$ gives
	\begin{align}   \label{optfr_des}
		\begin{split}
			&\mathbb{J}_o:\mathbb{R}^{N_e} \times \mathbb{R}^{N_e}  \rightarrow  \mathbb{R}^{N_e}\\
			&\delta H=\mathbb{J}_o[\hat{\kappa},\hat{\mu}]\left(
			\begin{array}{c}
				\delta \kappa\\
				\delta \mu
			\end{array}\right).
		\end{split}
	\end{align}
\end{definition}

\noindent
Now, we give further details on these operators.

\subsubsection{Discretised forward operator $\mathbb{H}_o[\kappa_0,\mu_0]$}
Plugging \eqref{piecewise} into \eqref{var_fo}, together with taking $\nu(r)$ to be a basis function $\varphi_p(r)$, gives a linear system
\begin{align} \label{system_forward}
	\boldsymbol{A}_o[\kappa^e,\mu^e]\Phi_0 =G.
\end{align}
Using this, the heating field is approximated in \Scom{a} piecewise constant basis \Scom{as} \cite{Tarvainen,Gao,Gao-b}
\begin{align} \label{forward_out}
	H = \hat{\mu} \circ \mathbb{I} \Phi.
\end{align}
\Scom{Here, $\circ$ denotes an element-wise product and $\mathbb{I}$ is a \textit{node-to-element} map. For $\mathbb{I}$ we use the $\mathcal{L}^2$ orthogonal projection from the space of nodal representations to elemental representations.} The action of $\mathbb{I}$ on a nodal vector $\theta$ restricted at element $j$ is in the form
\begin{align} \label{intp1}
	\left( \mathbb{I} \theta \right)_j=
	\frac{1}{d+1} \sum_{p   \in  \ell(j)} \theta_p
\end{align}
where $p$ denotes the nodal index, and $\ell(j)$ is a set of $d+1$ nodes that correspond to element $j$. We also define $\mathbb{I}^+$ as a map from the space of elemental vectors to nodal vectors. The action of $\mathbb{I}^+$ on an elemental vector $\Theta$ restricted at node $p$ is given by
\begin{align} \label{intp}
	\left(\mathbb{I}^+ \Theta \right)_p=\sum_{j=1}^{N_e} \Theta_j \int_{t_j}\varphi_p(r) dr= \frac{1}{d+1}\sum_{j \in l(p)} S_j  \Theta_j,
\end{align}
where $S$ is the vector of volume of elements, and $l(p)$ is the set of elements that are connected to node $p$. Also note that from \eqref{intp1} and \eqref{intp},  $ \mathbb{I}^+ \Theta = \mathbb{I}^T \left( S \circ \Theta \right) $, where $T$ denotes the transpose. In the sequel, we will use the notation $X=[\hat{\kappa},\hat{\mu}]$ and $\delta X=[\delta\hat{ \kappa},\delta \hat{\mu}]$.

\subsubsection{The matrix-free action of $\mathbb{J}_o[X]$ on perturbation $\delta X$}
Here, we explain how the action of $\mathbb{J}_o[X]$  on a perturbation $\delta X$ is approximated in a matrix-free manner. To do this, we approximate the perturbation field $\delta \phi(r)$  in a piecewise-linear basis in the same way as \eqref{piecewise}. Plugging the nodal vector $\delta \Phi$ into \eqref{var_fr} gives a linear system for calculation of $\delta \Phi$ in the form
\begin{align}
	A_o[\kappa^e,\mu^e]\delta \Phi=-\boldsymbol{A}_{\delta,o}[\delta \kappa^e, \delta \mu^e] \Phi_0,
\end{align}
where $\Phi_0$ has been computed using \eqref{system_forward}. Finally, the perturbation \Scom{in} the heating field $\delta H$ is computed as
\begin{align}  \label{dh_des}
	\delta H = \delta \hat{\mu} \circ  \mathbb{I} \Phi_0+\hat{\mu} \circ \mathbb{I} \delta \Phi.
\end{align}

\subsubsection{The matrix-free action of $\mathbb{J}_o^*[X]$  on $H$}
Here we explain how the action of the adjoint of  Fr\'echet derivative $\mathbb{J}_o^*$ on an elemental vector $H$ can be approximated using an \textit{adjoint-then-discretise} method. Plugging $\tilde{H}$ from \eqref{piecewise2} into the variational form of the adjoint formula \eqref{var_f_adj} gives a linear system in the form
\begin{align} \label{system_adjoint}
	\boldsymbol{A}_0[\kappa^e, \mu^e]   \tilde{H} = -\mathcal{G},
\end{align}
where
\begin{align}  \label{system_adjoint1}
	\mathcal{G}&=\mathbb{I}^+ (\hat{\mu}_0 \circ   H).
\end{align}
\Scom{Finally, given the nodal vectors $\tilde{H}$ and $\Phi_0$, we will approximate the action of the adjoint using \eqref{opt_adjoint}. For this we must choose how to calculate the products that appear in \eqref{opt_adjoint}, and then how to project onto the space of elemental vectors. We will use the following Lemma which also makes use of the matrices
	\begin{align}    \label{dak}
		\frac{\partial{\boldsymbol{A}}_o}{\partial{\hat{\kappa}_j}}=\int_{t_j}  \nabla \varphi_k \cdot \nabla  \varphi_p \  dr
	\end{align}
	and
	\begin{align}   \label{dam}
		\frac{\partial{\boldsymbol{A}}_o}{\partial{\hat{\mu}_j}}=\int_{t_j}  \varphi_k \varphi_p \  dr.
	\end{align}
}

If the products in \eqref{opt_adjoint} are calculated by first multiplying the nodal functions, and then using an $\mathcal{L}^2$-orthogonal projection on the space of elemental functions, then for the discretisation of the adjoint we have
\begin{align}\label{prodtproj}
	\begin{split}
		\mathbb{J}_o^* (H)_j= \left(
		\begin{array}{c}
			\frac{1}{S_j}  \tilde{H}^T \frac{\partial \boldsymbol{A}_o}{\partial \hat{\kappa}_j } \Phi_0  \\
			\frac{1}{S_j}  \tilde{H}^T \frac{\partial \boldsymbol{A}_o}{\partial \hat{\mu}_j } \Phi_0  + H^T (\mathbb{I} \Phi_0 )
		\end{array}
		\right).
	\end{split}
\end{align}

\begin{remark}
	Using a \textit{discretise-then-adjoint method}, the adjoint of  Fr\'echet derivative of the optical forward operator will be in the form (cf. \cite{Gao-b}, eq. (27))
	\begin{align} \label{num_adj}
		\begin{split}
			\mathbb{J}_{o,\text{dis}}^* (H)_j= \left(
			\begin{array}{c}
				\left[ (\mathbb{I}\boldsymbol{A}_o^{-1})^T( \hat{\mu}_0 \circ  H)   \right]^T \left(-\frac{\partial \boldsymbol{A}_o}{\partial \hat{\kappa}_j }\right) \Phi_0  \\
				\left[ (\mathbb{I}\boldsymbol{A}_o^{-1})^T( \hat{\mu}_0 \circ H)   \right]^T \left(-\frac{\partial \boldsymbol{A}_o}{\partial \hat{\mu}_j }\right)  \Phi_0  + H^T (\mathbb{I} \Phi_0 )
			\end{array}
			\right).
		\end{split}
	\end{align}
	Using \eqref{intp} and the fact that $ \mathbb{I}^+ \Theta = \mathbb{I}^T \left( S \circ \Theta \right) $, it can be shown that \eqref{num_adj} matches \eqref{prodtproj}.
\end{remark}

\begin{remark}
	From a theoretical point of view, it is perhaps more natural to make the projections from nodal to elemental representations in the $\mathcal{L}^\infty$ norm since $\kappa,\mu \in \mathcal{L}_+^\infty$. The $\kappa$ component of the adjoint, i.e., $\nabla  \phi_0(r) \cdot \nabla \tilde{h}(r)$, is constant at each element, and thus the $\mathcal{L}^\infty$ projection is the same as the $\mathcal{L}^2$ orthogonal projection. However, the two terms included in the $\mu$ component of the adjoint, i.e., $\phi_0(r) h(r)+\phi_0(r) \tilde{h}(r)$, will be different if $\mathcal{L}^\infty$ projection is used. Our simulation experiment showed us that an $\mathcal{L}^\infty$ projection gives almost the same reconstruction as the $\mathcal{L}^2$-orthogonal projection, \Scom{but dramatically increases the computational cost}. Therefore, we used the $\mathcal{L}^2$-orthogonal projection for approximation of the unknown optical coefficients.
\end{remark}

\subsection{Numerical computation of the acoustic operators ($\Lambda_a$ and its adjoint)}
Here, for a numerical computation of the acoustic portion of the problem, we used a k-space pseudo-spectral time-domain (PSTD) method \cite{Tabei,Cox}. Applying this method, the spatial gradients are approximated in a frequency domain using a Fast Fourier Transform (FFT), and the temporal gradients are approximated using finite difference schemes \cite{Tabei,Cox}.

\subsubsection{Discretisation of acoustic fields}
We approximate the acoustic forward and adjoint operators on a uniform rectilinear grid staggered in space and time \cite{Tabei}. We denote the position of a given grid point in Cartesian coordinates  by $r_\zeta$, where $\zeta = (\zeta_1,  ...  , \zeta_d) \in  \{1, \ ... \ , N_1\}  \times  ... \times \{1,  ...  , N_d\} $ with $N=\prod_{i=1}^d N_i$ the total number of grid points. We denote the grid separation along direction $i$ by $\Delta r_i$.  The time accessible to detectors, i.e., $t \in (0,T)$, is sampled with a temporal separation of $\Delta t$ so that the time step $n$ corresponds to the time instant $t_n=n \Delta t$.
To accommodate a staggered spatial grid, we introduce operator $\mathbb{T}_i$, which shifts the point $r$ by $\Delta r_i/2 $ \Scom{in coordinate $i$, and acts on a field as}
\begin{align}
	\mathbb{T}_i f(r)= f(\mathbb{T}_i r).
\end{align}
To accommodate a staggered temporal grid, we shift the field at time $t_n$ by $-\Delta t/2$, i.e., $t_{n-1/2}=n \Delta t- \Delta t/2$.
Based on the above, the discretised particle velocity vector at time step $n$ is denoted by $v_{(i;\zeta;n-1/2)} \in \mathbb{R}_i^d \times \mathbb{R}_\zeta^N \times \mathbb{R}_n^{N_t+1}$. This approximates the actual velocity on the staggered spatial grid as
\begin{align}
	v_{(i;\zeta;n-1/2)}\approx \mathbb{T}_i v_i(r_\zeta,t_{n-1/2}).
\end{align}
We also approximate the acoustic density as $\rho_{(i;\zeta;n)}\approx  \rho_i(r_\zeta,t_n) \in \mathbb{R}_i^d \times \mathbb{R}_\zeta^N \times \mathbb{R}_n^{N_t+1}$ and the scalar pressure as $p_{(\zeta;n)}\approx  p(r_\zeta,t_n) \in   \mathbb{R}_\zeta^N \times \mathbb{R}_n^{N_t+1}$. (Note that acoustic density is not a physical vector, but it is changed to a vector to accommodate the numerical method.) 

\subsubsection{Discretisation of medium's acoustic properties}
We define a discretised variant of acoustic parameters as diagonal matrices of size $N \times N$ acting on the spatial index $\zeta$ of acoustic fields. We denote a discretised variant of $c_0(r)$ by $\bar{c}$. Also, an approximation of the ambient density $\rho_0(r)$ is denoted by $\bar{\rho}$. We also define $\rho_0(r)$ on a spatial grid staggered in coordinate $i$ by $\bar{\rho}_i$.

Additionally, we define a discretised variant of the absorption and dispersion proportionality coefficients as $\bar{\tau}$ and $\bar{\eta}$, respectively. Also, a discretisation of the associated fractional Laplacian operators $N \times N$ is given by
\begin{align}
	Y_\text{abs}=F^{-1}\left\{ \{k^{y-2}F\{ \cdot\} \right\},   \   \   Y_\text{dis}=F^{-1} \left\{k^{y-1}F\{ \cdot\} \right\},
\end{align}
where $F$ and $F^{-1}$ denote the FFT and its inverse, respectively. Applying the k-space pseudo-spectral method on a staggered spatial grid, the spatial gradient in coordinate $i$ is in the form
\begin{align}  \label{grr}
	\frac{\partial \{\cdot\}}{\partial r_i^{\pm}}= F^{-1}  \left\{ \mathbf{i} k_i  e^{\pm \textbf{i} k_i \Delta r_i /2} \text{sinc}  \left(  \bar{c}_\text{ref}   k \Delta t/2  \right)     F\{\cdot\} \right\},
\end{align}
where $\text{sinc}(\bar{c}_\text{ref}   k \Delta t)/2$ enforces a k-space correction to the spatial gradient using a reference sound speed $\bar{c}_\text{ref}$ in order to minimise the numerical dispersion errors accumulated by the temporal integrations \cite{Tabei}. For further details on the k-space pseudo-spectral method, see \cite{Tabei,Cox}.

Because of the computation of the spatial derivatives via an FFT, wave wrapping may occur, i.e., the acoustic waves leaving one side of the grid reenter the opposite side. To avoid wave wrapping, an absorbing boundary condition, referred to as perfectly matched layer (PML), is added to each side of the grid. Here, the action of a PML on a field in coordinate $i$ is denoted by $A_i \in \mathbb{R}^{N \times N}$. This is defined as \cite{Tabei}
\begin{align} \label{dis-pml}
	A_i= \text{diag} \left( e^{-\alpha_{a,i} \Delta t/2} \right),
\end{align}
where $\alpha_{a,i}$ is the attenuation enforced by the PML in coordinate $i$.

Having defined a discretisation of the optical forward operator, as well as acoustic fields and medium, we now complete the definition of our discretised composite operator.

\begin{definition}
	A discretisation of the composite \textit{opto-acoustic} forward operator $\Lambda$  gives a map from optical coefficients to a set of time series of measured  data in the form
	\begin{align}   \label{com_des}
		\begin{split}
			&\mathbb{H} : \mathbb{R}^{N_e }  \times \mathbb{R}^{N_e } \rightarrow \mathbb{R}^{N_s N_t}\\
			&\hat{P}=\mathbb{H}_a  \mathbb{H}_o [X],
		\end{split}
	\end{align}
	where $\mathbb{H}_o$ and $\mathbb{H}_a$ represent the discretised optical and acoustic forward operators, respectively.
	Also, an operator representing a discretisation of the Fr\'echet derivative of the forward operator $D_X \Lambda$ is defined by
	\begin{align}    \label{J_des}
		\begin{split}
			&\mathbb{J}: \mathbb{R}^{N_e }  \times \mathbb{R}^{N_e } \rightarrow \mathbb{R}^{N_s N_t}\\
			&\delta \hat{P}= \mathbb{H}_a \mathbb{J}_o[X] \delta X.
		\end{split}
	\end{align}
	\Scom{Using \eqref{J_des}, we will use as a discretisation of the adjoint of Fre\'chet derivative operator the map
		\begin{align}
			\begin{split}
				\mathbb{J}^* : \mathbb{R}^{N_s N_t} \rightarrow \mathbb{R}^{N_e} \times \mathbb{R}^{N_e}\\
				\partial X = \mathbb{J}_{o}^*[X] \mathbb{H}_{a}^* \left(\partial \hat{P}\right).
			\end{split}
		\end{align}
	}
\end{definition}

\subsubsection{Numerical computation of the acoustic forward operator $\mathbb{H}_a$}    \label{num-f}
We now explain a discretisation of the acoustic forward operator ($\mathbb{H}_a$). We will use this later for calculation of the acoustic adjoint operator. For a numerical implementation of the acoustic forward operator given by \eqref{mom}, \eqref{mass}, \eqref{state} and \eqref{initial} using a k-space pseudo-spectral method, we used an open-source code, which is freely available on the \textit{k-Wave} website \cite{Treeby-b}.

Using this code, the discretised initial pressure distribution at $t=0$, denoted by $P_0$ is applied as an injection of mass, referred to as \textit{additive source}. To do this, $P_0$ must be split over two time steps $n= \left\{-1/2,+1/2 \right\}$. (cf. \cite{Arridge}, Appendix B, or the \textit{k-Wave} manual \cite{Treeby-b}). This gives a source in coordinate $i$ as \cite{Arridge}
\begin{align}    \label{s1}
	s_{(i;\zeta;n+1/2)} =
	\begin{cases}
		\frac{1}{2 d\Delta t \ \bar{c}^2 } \mathbb{S}P_0  \quad   &n=-1,0 \\
		0                                              \quad   &\text{otherwise},
	\end{cases}
\end{align}
where $\mathbb{S}$ is a symmetric smoothing operator that is applied in order to mitigate unexpected oscillations of $P_0$. (cf. the \textit{k-Wave} manual in \cite{Treeby-b}). Additionally, since $s_i$ is added to the mass equation, a factor $\frac{1}{\Delta t \ \bar{c}^2}$ has been applied in order to account for the conversion of units from pressure to the time rate of density (cf. \cite{Arridge}, Appendix B, or the \textit{k-Wave} manual in \cite{Treeby-b}). Using the definitions given above, the calculation of $\mathbb{H}_a$ proceeds as follows.

\noindent \textit{Start at iterate $n=-1$ with initial conditions $p_{(\zeta;n=-1)} = 0$, $v_{(i;\zeta;n=-3/2)}=0$ and $\rho_{(i;\zeta;n=-1)}=0$, and terminate at iterate $n=N_t-2$.}
\smallskip

\noindent
\textit{1. Update the particle velocity vector field (conservation of momentum \eqref{mom}):}
\begin{align} \label{num-vel}
	\begin{split}
		v_{(i;\zeta;n+\frac{1}{2})}&= A_i   \bigg[A_i \  {v_{(i;\zeta;n-\frac{1}{2})}}
		- \frac{\Delta t}{\bar{\rho}_i} \frac{ \partial}{\partial r_i^+}     p_{(\zeta;n)}  \bigg].
	\end{split}
\end{align}

\noindent
\textit{2. Update the acoustic density field (conservation of mass \eqref{mass}) also adding source:}
\begin{align} \label{num-rho}
	\begin{split}
		\rho_{(i;\zeta;n+1)}&= A_i   \bigg[A_i \  {\rho_{(i;\zeta;n)}}
		- \Delta t  \bar{\rho} \frac{ \partial}{\partial r_i^-}     v_{(i;\zeta;n+\frac{1}{2})}    \bigg] + \ \Delta t \ s_{(i;\zeta;n+1/2)}.
	\end{split}
\end{align}


\noindent
{\textit{3. Update the scalar pressure field (equation of state \eqref{state}):}
	\begin{align} \label{num-p}
		\begin{split}
			p_{(\zeta;n+1)}&= \bar{c}^2   \bigg[   \left( I_N-
			\bar{\eta} Y_\text{dis}\right) \sum_{i=1}^d \rho_{(i;\zeta;n+1)} +  \bar{\tau} Y_\text{abs} \sum_{i=1}^d   A_i \bar{\rho} \frac{ \partial}{\partial r_i^-}     v_{(i;\zeta;n+\frac{1}{2}) }  \bigg]
		\end{split}
	\end{align}
	where $I_N$ denotes an identity matrix of size $N \times N$, and the last term is actually the action of $ \bar{\tau} Y_\text{abs}$ on $ \sum_{i=1}^d  \frac{\partial }{\partial t } \rho_{(i;\zeta;n+1)}$, and is derived from \eqref{num-rho}.
	
	\noindent
	\textit{4. Compute the measured pressure at detectors:}
	This is defined using
	\begin{align} \label{num-m}
		\hat{P}_{n+1}= \mathbb{M} p_{(\zeta;n+1)},
	\end{align}
	where $\mathbb{M}$ denotes a discretised variant of $\mathcal{M}$ (cf. section \ref{skp}), and includes an interpolation operator for mapping the acoustic pressure field from grid points to position of ultrasound detectors. It is worth mentioning that $\mathbb{M}$ also depends on the size and properties of detectors. These effects are neglected in our study by assuming the detectors sample the pressure pointwise.
	
	\subsubsection{A \textit{discretise-then-adjoint} method for derivation of the acoustic adjoint operator}
	Having defined a discretisation of the acoustic forward operator ($\mathbb{H}_a$), we now explain how to calculate the acoustic adjoint operator ($\mathbb{H}_a^*$) using a \textit{discretise-then-adjoint} method. The derivation of $\mathbb{H}_a^*$ in this section is a modification to the study of \cite{Huang} in the sense that the effects of the PMLs and an \textit{additive source} are incorporated in calculation of the adjoint. Additionally, in contrast to \cite{Huang}, we will represent $\mathbb{H}^*_a$ as a discretised linear system of PDEs, the same as our representation for $\mathbb{H}_a$ \cite{Huang}.

	To derive this adjoint, a matrix form of $\mathbb{H}_a$  must be derived using the details given in section \ref{num-f}. To do this, we start with definition of diagonal matrices $C=\bar{c}^2 \in \mathbb{R}^{ N\times N}$ and $Q \in \mathbb{R}^{dN \times dN}$ with diagonal a \Scom{$d$-times} stack of diagonals of $\bar{\rho} \in \mathbb{R}^{N \times N}$. Also, we define a diagonal matrix $Q_s \in \mathbb{R}^{dN \times dN}$ with diagonal a stack of diagonals of $\bar{\rho}_i, \ (i \in \left\{ 1,...,d\right\})$. We will also use $A \in \mathbb{R}^{dN \times dN}$ as a diagonal matrix with diagonal a stack of diagonals of $A_i \  \  (i \in \{1,...,d \})$.
	
	We now define a stack of coordinate-dependent particle velocity and  acoustic density fields $v_{(i;\zeta,n)}\in \mathbb{R}^N , \ \rho_{(i;\zeta,n)}  \in \mathbb{R}^N \  (i \in \{1,...,d \})$ as $\bar{v}_{n-1/2} \in \mathbb{R}^{dN}$ and $\bar{\rho}_n \in \mathbb{R}^{dN}$, respectively. We also define $\bar{p}_n=p_{(\zeta,n)} \in \mathbb{R}^N $. A stack of all \Scom{these} vector fields yields $\bar{z}_n = \{ \bar{v}_{n-1/2}^T \ \bar{\rho}_n^T \ \bar{p}_n^T \}^T \in \mathbb{R}^{7N}$ at time step $n$.
	
	Also, let $\mathcal{S} \in \mathbb{R}^{7NN_t \times N} $ be the map from \Scom{the discretised sought after initial pressure $P_0$} to an additive source $S$. In particular, at time step $n$, $S_{n+1/2}= \mathcal{S}_n P_0$, where
	\begin{align} \label{s2}
		\mathcal{S}_n P_0 = \Delta t   \  \mathcal{I}_{s}^{z}   s_{(i;\zeta;n+1/2)}=  \mathcal{I}_{s}^{z}  1_i \otimes \left( \frac{1}{2d C}\right)   \mathbb{S}P_0 \quad \mbox{for $n = -1, \ 0$},
	\end{align}
	and $\mathcal{S}_n P_0 = 0$ for $n \neq -1, \ 0$. Here, $ s_{(i;\zeta;n+1/2)}$ has been derived from \eqref{s1}, and $1_i \otimes  \in \mathbb{R}^{dN \times N}$ is the adjoint of $\sum_{i=1}^d \in \mathbb{R}^{N \times dN}$. Also, from \eqref{num-rho} and \eqref{num-p}, $ \mathcal{I}_s^{z} \in \mathbb{R}^{7N \times dN}$ is the map from  $s_{(i;\zeta;n+1/2)}$ to $\bar{z}_n$ in the form
	\begin{align} \label{s3}
		\begin{split}
			\mathcal{I}_s^{z}=
			\begin{bmatrix}
				0_{dN \times dN}\\
				I_{dN \times dN} \\
				C\left(\left(I_N - \bar{\eta} Y_\text{dis}\right) \sum_{i=1}^d \right)_{N \times dN}
			\end{bmatrix}.
		\end{split}
	\end{align}
	
	\noindent
	We now introduce an operator $T \in \mathbb{R}^{7N \times 7N} $ for defining a discretised formula for the time sequence of fields as
	\begin{align}\label{t}
		\bar{z}_{n+1}=T \bar{z}_n+ S_{n+1/2}.
	\end{align}
	Additionally, we introduce a measurement matrix $\boldsymbol{M} =  \mathbb{M}  \mathcal{I}_{z}^{p} $ with $ \mathcal{I}_{z}^{p} \in \mathbb{R}^{N \times 7N}$ the projection from the space of $\bar{z}_n$ to the space of $\bar{p}_n$. \Scom{Based on this}, we now define our discretised acoustic forward operator.
	\begin{definition}\label{Hdef}
		The discretised acoustic forward operator is defined by
		\begin{align}
			\mathbb{H}_{a}  &: \mathbb{R}^{N} \rightarrow \mathbb{R}^{N_s N_t}\\
			\mathbb{H}_{a}  &= \mathbb{H}_T \mathcal{S},
		\end{align}
		where $\mathbb{H}_T: \mathbb{R}^{7NN_t} \rightarrow \mathbb{R}^{N_s N_t}$ is in the form
		\begin{align} \label{for_map_dis}
			\hat{P}=\mathbb{H}_T S, \quad \hat{P}_{n}            =\boldsymbol{M}\bar{z}_n \quad (n \in \{0,...,N_t-1\}), \quad
			\hat{P}              =\Big[\hat{P}_{n+1} \Big]_{n=-1}^{N_t-2}.
		\end{align}
		Here, $\bar{z}_n$ is determined by \eqref{t} with initial condition $\bar{z}_{-1} = 0$, and $\hat{P} \in \mathbb{R}^{N_s N_t}$ is a time-series stack of measured data at iterates $n \in \{0,..., N_t-1\}$.
	\end{definition}

	\begin{lemma} \label{adjcor}
		The action of the adjoint operator $\mathbb{H}_{a}^*$ on $\hat{P}$ is given by
		\begin{align} \label{S*H*}
			\mathbb{H}_{a}^*:\mathcal{S}^* \mathbb{H}_T^*  \hat{P} = \sum_{n=-1}^{N_t-2} \mathcal{S}^*_{n} \bar{z}^*_{N_t-2-n}
		\end{align}
		where $\bar{z}^*_{n}$ is determined by
		\begin{align} \label{col}
			\bar{z}^*_{-1} = 0, \quad \bar{z}^*_{n+1} = T^* \bar{z}_n^* + \boldsymbol{M}^* \hat{P}_{N_t - 2 - n} \quad (n \in \{-1, \ ... \ , N_t - 2\}).
		\end{align}
		\begin{proof}
			We have given the proof in \cite{Javaherian-e}, Lemma 2 and and Corollary 1.
		\end{proof}
	\end{lemma}
	\noindent
	Note that in the case we are considering the sum in \eqref{S*H*} is actually just two terms. \Scom{By commuting $n=-1,0$ in \eqref{s2} to $\boldsymbol{M}$ using the fact the the forward operator is linear, we define a source for our acoustic adjoint operator as
		\begin{align}      \label{adj_data}
			\begin{split}
				s_{n+1/2}^\text{adj}=   \mathbb{M}^*
				\begin{cases}
					\hat{P}_{N_t-1},                            \quad &n=-1	\\
					\hat{P}_{N_t-n-1}+\hat{P}_{N_t-n-2},        \quad &n=0,...,N_t-2     \\
					\hat{P}_0,                    \quad & n=N_t-1.
				\end{cases}
			\end{split}
		\end{align}
		From \eqref{adj_data}, equation \eqref{col} in Lemma \ref{adjcor} can be modified as
		\begin{align}   \label{col1}
			\bar{z}^*_{-1} = 0, \quad \bar{z}^*_{n+1} = T^* \bar{z}_n^* + S_{n+1/2}^\text{adj}  \quad (n \in \{-1, \ ... \ , N_t - 2\}),
		\end{align}
		where
		\begin{align}  \label{col1s}
			S_{n+1/2}^\text{adj}= \left(\mathcal{I}_z^{p}\right)^* s_{n+1/2}^\text{adj}.
		\end{align}
	}
	
	\noindent
	Now, we derive a matrix form of $T$ and $T^*$ using the forward model presented in section \ref{num-f} to show how multiplication by each of them may be computed.
	\noindent
	We also denote a matrix form of the k-space spatial gradient in \eqref{grr} as
	\begin{align} \label{dgr}
		\nabla_i^{\pm}=\frac{\partial}{\partial r_i^{\pm}} \in \mathbb{R}^{N \times N}.
	\end{align}
	\Scom{Using this, we will use $\Phi \in \mathbb{R}^{dN \times dN}$ and $\Psi \in \mathbb{R}^{dN \times dN}$ composed of submatrices
		\begin{align} \label{phi}
			(\Phi)_{ij}= -\Delta t  \ \delta_{ij}  \nabla_j^+
		\end{align}
		and
		\begin{align} \label{psi}
			(\Psi )_{ij}= -\Delta t  \ \delta_{ij}   \nabla_{j}^-.
		\end{align}
	}
	To make the notation more compact we also introduce matrices $D,E \in \mathbb{R}^{N \times dN}$ and $G \in \mathbb{R}^{N \times N}$ whose actions are given by
	\begin{align}  \label{D4}
		D \bar{v}=  C \left(I_N-  \bar{\eta} Y_\text{dis}-\bar{\tau} \frac{Y_\text{abs}}{\Delta t}\right) \sum_{i=1}^d  (A Q \Psi A^2 \bar{v})_i
	\end{align}
	\begin{align}\label{E}
		E \bar{\rho}= C  \left(  I_N-\bar{\eta} Y_\text{dis}\right)  \sum_{i=1}^d  (A^2  \bar{\rho})_i
	\end{align}
	\begin{align}\label{G4}
		G \bar{p}=   C \left(I_N-  \bar{\eta} Y_\text{dis}-\bar{\tau} \frac{Y_\text{abs}}{\Delta t}\right) \sum_{i=1}^d   (A Q \Psi A Q_s^{-1} \Phi 1_i \otimes  \bar{p})_i
	\end{align}
	
	\noindent
	The matrix $T$ is then given, in block form, by
	\begin{align} \label{Tprod}
		T =
		\left (
		\begin{matrix}
			A^2                                 & 0                & A Q_s^{-1}\Phi   1_i \otimes                 \\
			A Q \Psi A^2                   &A^2            &A  Q \Psi A Q_s^{-1} \Phi 1_i \otimes   \\
			D                                     &E                 &G
		\end{matrix}
		\right ).
	\end{align}
	From \eqref{Tprod}, the adjoint operator $T^*$ will be in the form
	\begin{align} \label{T*prod}
		T^* =
		\left (
		\begin{matrix}
			A^2                                                     &A^2 \Psi^* Q A                                                              & D^*\\
			0                                                          &A^2                                                                                & E^*\\
			\sum _{i=1}^d  \Phi^* Q_s^{-1} A       & \sum _{i=1}^d \Phi^* Q_s^{-1} A \Psi^* Q A         & G^*\\
		\end{matrix}
		\right).
	\end{align}
	\noindent
	From \eqref{T*prod} and using \eqref{s3}, \eqref{col1} and \eqref{col1s}, as well as calculating the adjoints of $D$, $E$ and $G$ from \eqref{D4}, \eqref{E}, and \eqref{G4}, the time sequence of adjoint fields is given iteratively by
	\begin{align}  \label{t1}
		\begin{split}
			\bar{\rho}_{n+1}&= A^2 \left[ \bar{\rho}_n +1_i \otimes \left( I_N-Y_\text{dis} \bar{\eta}\right) C \bar{p}_n \right]\\
			\bar{v}_{n+1/2} &= A^2 \left[   \bar{v}_{n-1/2}+ \Psi^* Q A^{-1}  \left( \bar{\rho}_{n+1}  -A^2 1_i \otimes \frac{Y_\text{abs}}{\Delta t}\bar{\tau}  C \bar{p}_n   \right)    \right]\\
			\bar{p}_{n+1}   &= \left(\sum_{i=1}^d \Phi^* Q_s^{-1}  A^{-1}\bar{v}_{n+1/2}\right)  +s_{n+1/2}^\text{adj}\\
			\bar{p}^\text{sol}   &=  \boldsymbol{\vartheta} \left( \sum_{i=1}^d  1_i \otimes  \left(I_N- Y_\text{dis} \bar{\eta} \right)  C \bar{p}_{n+1}+  \bar{\rho}_{n+1}\right) \ (n=N_t-2),
		\end{split}
	\end{align}
	where $ \bar{p}^\text{sol}:=\mathbb{H}_a^* (\hat{P}) $, {\textcolor{red}{and $\boldsymbol{\vartheta} =  \mathbb{S} \left( \frac{1}{2d C}\right) $ using \eqref{s2}}}. Now, applying the replacements $\boldsymbol{\hat{\rho}}_{n+1}= QA^{-1} \bar{\rho}_{n+1}$ and $\boldsymbol{\hat{v}}_{n+1}= (AQ_s)^{-1} \bar{v}_{n+1/2}$, together with $\Phi^*=-\Psi$ and $\Psi^*= -\Phi$, gives
	\begin{align}  \label{t41}
		\begin{split}
			\boldsymbol{\hat{\rho}}_{n+1}&= A \left[ A  \boldsymbol{\hat{\rho}}_n +Q1_i \otimes \left( I_N-Y_\text{dis} \bar{\eta}\right) C \bar{p}_n \right]\\
			\boldsymbol{\hat{v}}_{n+1/2} &= A \left[  A \boldsymbol{\hat{v}}_{n-1/2}- Q_s^{-1} \Phi  \left( \boldsymbol{\hat{\rho}}_{n+1}  - QA 1_i \otimes \frac{Y_\text{abs}}{\Delta t}\bar{\tau}  C \bar{p}_n   \right)    \right]\\
			\bar{p}_{n+1}   &= \left(\sum_{i=1}^d - \Psi  \boldsymbol{\hat{v}}_{n+1/2}\right)  +s_{n+1/2}^\text{adj}\\
			\bar{p}^\text{sol}   &=  \boldsymbol{\vartheta} \left( \sum_{i=1}^d  1_i \otimes  \left(I_N- Y_\text{dis} \bar{\eta} \right)  C \bar{p}_{n+1}+ A Q^{-1} \boldsymbol{\hat{\rho}}_{n+1}\right) \ (n=N_t-2).
		\end{split}
	\end{align}
	Our numerical experiments showed that an operator $\mathbb{H}_{a}^*$ that is calculated using a \textit{discretise-then-adjoint} method \eqref{t41} satisfies an adjoint test with a higher accuracy than using an \textit{adjoint-then-discretise} method {\textcolor{red}{used in \cite{Javaherian}}. Therefore, we used \eqref{t41} for the acoustic adjoint operator. 
		
		\section{Iterative model-based approaches for the direct problem of QPAT}  \label{itr}
		Having defined a discretisation of opto-acoustic forward operator $\mathbb{H}$, the Fr\'echet derivative operator $\mathbb{J}$, and its adjoint $\mathbb{J}^*$, we now explain the iterative approaches we will use for minimisation of  \eqref{obj_couple}.
		
		\Scom{Considering \eqref{h}, the dependence of heating field $h$ on $\kappa, \mu$ is nonlinear. Furthermore, because of high scattering of light in tissue media, simultaneous reconstruction of $\kappa,\mu$ from $h$ can be highly ill-posed.} As a result, a minimisation of \eqref{obj_couple} is a non-convex, nonlinear and ill-posed inverse problem. It  has been shown that simultaneous reconstruction of $\kappa$ and $\mu$  using a single optical excitation, i.e., $N_q=1$, does not have a unique solution \cite{Bal-c}. However, the uniqueness and stability, in appropriate norms, of this inverse problem using $N_q>1$ optical excitations under some geometric constraints has been established \cite{Bal-c}. \Scom{It is also worth mentioning that using our forward acoustic operator, which can be adapted to acoustically heterogeneous and lossy media, the direct problem of QPAT is more ill-posed than existing studies, for which spherical mean Radon transform or Green's function techniques have been used for solving the acoustic portion of the forward operator using an acoustically homogeneous and lossless medium, which does not hold in practice \cite{Haltmeier,Gao-b,Pulkkinen-b}.}
		\begin{remark}
			Since the magnitude of $\hat{\kappa}$ is often 1 to 2 orders of magnitude greater than $\hat{\mu}$, we follow \cite{Gao,Gao-b} and minimise $\epsilon$ (cf. \eqref{obj_couple}) with respect to scaled coefficients 
			\begin{equation} \label{scaling}
				\bar{X} =[\bar{\kappa}(\hat{\kappa}),\bar{\mu}(\hat{\kappa})]^T. 
			\end{equation}
			We will use two types of scaling. For the method described in section \ref{sec:ADMM} we use a linear scaling (see \eqref{linscale}) while for the methods described in sections \ref{sec:LD} and \ref{sec:PDIPM} we use logarithmic scaling (see \eqref{logscale}). In our numerical experiments the method of section \ref{sec:ADMM} did not converge with logarithmic scaling.
		\end{remark}
		\noindent
		\Scom{Using these, we now consider a minimisation problem with respect to a scaled vector of optical coefficients $\bar{X}$ in the form}
		\begin{align}  \label{varr}
			\bar{X}_* =\operatornamewithlimits{argmin}_{\bar{X}_l \leqslant    \bar{X}   \leqslant \bar{X}_u}    \    \epsilon(\bar{X}) =  \operatornamewithlimits{argmin}_{\bar{X}_l \leqslant    \bar{X}   \leqslant \bar{X}_u}    \frac{1}{2} \sum_{q=1}^{N_q} \norm{ \mathbb{H}_q[X[\bar{X}]]- \hat{P}_q }^2,
		\end{align}
		where $\epsilon$ is a non-convex, nonlinear and smooth function. \Scom{Also, $\bar{X}$ has been constrained by a lower bound $\bar{X}_l$ and an upper bound $\bar{X}_u$.} We will take two main approaches for solving \eqref{varr}, the first of which is a direct minimisation of the nonlinear objective function using a Quasi-Newton approach, and the second is to solve the minimsation problem as a sequence of convex and linearised subproblems using a matrix-free Jacobian-based method. We now explain these approaches.
		
		\noindent
		A fixed point iteration arising from the optimality conditions of $\epsilon$ gives a sequence
		\begin{align}  \label{sdm}	
			\bar{X}_{k+1}= \bar{X}_k + \alpha _k d_k,
		\end{align}
		where $d_k$ is a search direction for iteration $k$, and $\alpha_k$ is a step size along search direction $d_k$. (Throughout this manuscript, a subscript (resp. superscript) $k$ indicates an iteration for an inner (resp. outer) loop.)
		
		\subsection{Newton's methods}  \label{sm}
		From the second-order optimality condition for minimising $\epsilon$, the \textit{Newton search direction}, is derived using
		\begin{align} \label{sd}
			d_k=  -\mathcal{H}_k^{-1} \nabla \epsilon_k.
		\end{align}
		Here, $\nabla \epsilon_k$ is the first-order derivative of $\epsilon$ at  $\bar{X}_k$,
		and is computed using
		\begin{align} \label{gr}
			\nabla  \epsilon_k =  \sum_{q=1}^{N_q} \frac{\partial X}{\partial \bar{X}}[\bar{X}_k] \  \mathbb{J}_q^*  [X_k]  \left( \mathbb{H}_q[ X_k]- \hat{P}_q \right)   ,
		\end{align}
		and $\mathcal{H}_k$, is the second-order derivative of $\epsilon$ (the Hessian matrix). Using a \textit{Gauss-Newton} method, the Hessian matrix is approximated using
		\begin{align} \label{apj}
			\mathcal{H}_k= \sum_{q=1}^{N_q} \frac{\partial X}{\partial \bar{X}}[\bar{X}_k]  \  \mathbb{J}_q^* [X_k] \mathbb{J}_q[X_k]  \frac{\partial X}{\partial \bar{X}}[\bar{X}_k],
		\end{align}
		\Scom{where a term including the second-order derivatives of $\mathbb{H}_q$ has been neglected as compared to the exact Hessian.} 
		A class of approaches that utilise \eqref{sd} for minimisation of $\epsilon$ \Scom{are} called \textit{Newton's method}. Although \textit{Newton's methods} benefit from a quadratic (optimal) rate of convergence, they \Scom{pose} some practical limitations for QPAT, i.e.,
		
		a) The computation, storage and inversion of the Hessian matrix is expensive. To address this problem, an implicit inversion of $\mathcal{H}_k$ using an explicit form of $\mathbb{J}$ and $\mathbb{J}^*$ has been used in the context of diffuse optical tomography \cite{Schweiger}.
		
		b) The size of the discretised heating field $H$, which represents data for the optical portion of the forward operator, is on the same order as the size of unknown parameters (cf. \eqref{opt_des}). It has been shown that an explicit computation of the matrix $\mathbb{J}_o$  (resp. $\mathbb{J}$) requires at least $N_s+N_e$ times implementation of $\mathbb{H}_o$ (resp. $\mathbb{H}$) (See \cite{Gao}, page 11).
		
		c) Considering the time series of measured data $\hat{P}$, \Scom{which is of size $N_s N_t$}, the Jacobian matrix $\mathbb{J}$ is 
		dense, and thus a storage of $\mathbb{J}$ is impractical.
		
		We will later explain how we have addressed these challenges using an inexact Newton method.
		
		\subsection{Nonlinear gradient-based methods}  \label{sm1}
		An alternative to Newton's method is using gradient-based Quasi-Newton approaches, for which $\mathcal{H}_k$ is not computed explicitly, but is approximated using solely information included in the first-order gradients $\nabla{\epsilon}$ possibly at previous steps. Additionally, an inversion of $\mathcal{H}_k$ can be avoided using a direct approximation of $\mathcal{H}^{-1}_k$, for which the so-called BFGS method is often used. Since $\mathcal{H}^{-1}_k$ is a dense matrix, and BFGS method poses challenges regarding memory, a limited-memory variant of BFGS (L-BFGS) method is used \cite{Gao,Gao-b}.
		Using L-BFGS, $\mathcal{H}^{-1}_k$ is updated using the most recent $m$  pairs of $(\boldsymbol{s},\boldsymbol{y})$ given by
		\begin{align}
			\begin{split}
				\boldsymbol{s}_k &= \bar{X}_{k}-\bar{X}_{k-1}\\
				\boldsymbol{y}_k &= \nabla \epsilon_{k} - \nabla \epsilon_{k-1}.
			\end{split}
		\end{align}
		(In our study, we empirically use $m=5$.) We will also use $\boldsymbol{\rho}_k= 1/y_k^T s_k$, and an initial guess for the Hessian matrix in the form
		\begin{align}
			\mathcal{H}_{k,0}= \frac{\boldsymbol{s}_{k-1}^T \boldsymbol{y}_{k-1}} {\boldsymbol{y}_{k-1}^T  \boldsymbol{y}_{k-1}} I,
		\end{align}
		where $I$ is the identity matrix. By applying L-BFGS method to the constrained minimisation problem \eqref{varr}, the search direction $d_k$ is computed using Algorithm~\ref{alg:41}. (See \cite{Gao} and \cite{Gao-b} for applications on QPAT and direct QPAT.)
		\begin{algorithm}   
			\caption{L-BFGS (search direction): inner iteration $k$ }
			\label{alg:41}
			\begin{algorithmic}[1]
				\State Input: $ \nabla \varepsilon_k, \mathcal{H}_{k,0}^{-1}$
				\State Initialise: $\boldsymbol{q}= \nabla \varepsilon_k$
				\For {$i=k-1,k-2,...,k-m$}
				\State $\boldsymbol{\alpha}_i= \boldsymbol{\rho}_i \boldsymbol{s}_i^T \boldsymbol{q}$
				\State $\boldsymbol{q}=\boldsymbol{q}-\boldsymbol{\alpha}_i \boldsymbol{y}_i$
				\EndFor
				\State $\boldsymbol{r}=\mathcal{H}_{k,0}^{-1} \boldsymbol{q}$
				\For {$i=k-m,k-m+1,...,k-1$}
				\State $\boldsymbol{r}=\boldsymbol{r}+\boldsymbol{s}_i\left(\boldsymbol{\alpha}_i-\boldsymbol{\rho}_i \boldsymbol{y}_i^T \boldsymbol{r} \right)$
				\EndFor
				\State $\boldsymbol{d}_k=-\boldsymbol{r}$.
			\end{algorithmic}
		\end{algorithm}

		From the First-order Karush Kuhn Tucker (KKT) conditions associated with the constraint on $\bar{X}$, $d_k$ is projected onto the feasible region using (See \cite{Gao-b}, equation (37))
		
		\begin{align}\label{LBFGSsearch}
			\begin{split}
				d_k =
				\begin{cases}
					-\bar{X}_l, \  &\text{if} \  \bar{X}_k+ \boldsymbol{d}_k \leqslant \bar{X}_l\\
					\boldsymbol{d}_k, \   &\text{if} \  \bar{X}_l  < \bar{X}_k+ \boldsymbol{d}_k < \bar{X}_u\\
					-\bar{X}_u, \  &\text{if} \  \bar{X}_k+\boldsymbol{d}_k \geqslant \bar{X}_u.
				\end{cases}
			\end{split}
		\end{align}
		
		\noindent
		Using $d_k$ given by \eqref{LBFGSsearch}, the step size in \eqref{sdm} is chosen by a standard backtracking line search satisfying the \textit{Wolfe} conditions as well as the constraints
		\begin{align}  \label{wo}
			\begin{split}
				&\epsilon (\bar{X}_k +\alpha d_k) \leqslant \epsilon (\bar{X}_k ) +  c_1   \alpha_k   d_k^T  \nabla \epsilon (\bar{X}_k)   \\
				& d_k^T   \nabla\epsilon (\bar{X}_k + \alpha d_k) d_k    \geqslant   c_2   d_k^T  \nabla \epsilon (\bar{X}_k),
			\end{split}
		\end{align}
		together with enforcing the bounds associated with the constraint, i.e.,
		\begin{align} \label{cons}
			\bar{X}_l \leqslant    \bar{X}_k + \alpha d_k   \leqslant \bar{X}_u.
		\end{align}
		(See \cite{Gao,Gao-b}.)
		In \eqref{wo}, $ 0 < c_1 < c_2 <1  $ are user-defined parameters. Applying these conditions, $\alpha_k$ is chosen using a backtracking line search, as given in Algorithm~\ref{alg:42}.
		\begin{algorithm}
			\caption{Backtracking Line Search: inner iteration $k$}
			\label{alg:42}
			\begin{algorithmic}[1]
				\State Input: $c_1,c_2, \boldsymbol{\tau}, \bar{X}_l, \bar{X}_u$
				\State Initialise: $\alpha_0=1$
				\While {\eqref{wo} or \eqref{cons} are not satisfied}
				\State $\alpha=\boldsymbol{\tau} \alpha$
				\EndWhile
				\State $\alpha_k=\alpha_*$.
			\end{algorithmic}
		\end{algorithm}
		\noindent
		Here, $\boldsymbol{\tau}<1$ is a user-defined parameter.
		
		\section{Total Variation (TV) regularisation}  \label{tvr}
		To mitigate the ill-posedness of the problem, a regularisation functional must be added to the data fidelity term in \eqref{varr} \cite{Gao,Gao-b}. This results in a minimisation problem in the form
		\begin{align}\label{obj_dis}
			\bar{X}_*=\operatornamewithlimits{arg \  min}_{\bar{X}_l\leqslant \bar{X} \leqslant \bar{X}_u}  \left\{ \mathcal{F}:= \epsilon[\bar{X}]+ \lambda \mathcal{J}[\bar{X}]\right\},
		\end{align}
		where $\mathcal{J}[\bar{X}]$ and $\lambda$, respectively denote \Scom{the} regularisation functional and the regularisation parameter, \Scom{the latter of which makes a balance} between a fidelity to the measured data $\hat{P}$ and \Scom{to} \textit{a priori} knowledge about the \Scom{true solution}. In our study, based on \Scom{an} assumption that the optical coefficients are piecewise constant with sharp edges, we use $\mathcal{J}[\bar{X}]:= \mathcal{R}[\bar{\kappa}]+  \mathcal{R}[\bar{\mu}]$, where $\mathcal{R}[u]$ is a discretisation of the \textit{Total-Variation (TV)} functional
		\begin{align} \label{tvc}
			\int_\Omega | \nabla u | dr
		\end{align}
		with $u$ either $\bar{\kappa}$ or $\bar{\mu}$. Using \eqref{tvc},
		\begin{align} \label{tvm}
			\mathcal{J} [\bar{X}] = \| D \bar{X} \|_1
	\end{align}}
	with
	\begin{align} \label{dd}
		\begin{split}
			D= \begin{bmatrix}
				&D^{\bar{\kappa}}                  &0_{N_l \times N_e}\\
				&0_{N_l\times N_e}                    &  D^{\bar{\mu}}
			\end{bmatrix}.
		\end{split}
	\end{align}
	Here, $D^u \in \mathbb{R}^{N_l \times N_e}$ is a sparse matrix with $N_e$ and $N_l$ the total number of elements and the total number of internal edges between elements, respectively. Each row $D_l^{u}   \in \mathbb{R}^{1 \times N_e} \ (l \in \{1,...,N_l\})$ has two nonzero \Scom{components} at indices $j_\text{1}$ and $j_\text{2}$, which correspond to two elements connected by the internal edge $l$. These have values $a_l$ and $-a_l$ with $a_l$ the length (or area) of internal edge $l$ \cite{Borsic}. Using this, the gradient of $\mathcal{J}(\bar{X})$ is a nonlinear operator in the form
	\begin{align} \label{jm}
		\begin{split}
			&\nabla \mathcal{J}: \mathbb{R}^{2N_e} \rightarrow \mathbb{R}^{2N_e}\\
			&\nabla \mathcal{J}[\bar{X}] =M[\bar{X}] \bar{X},
		\end{split}
	\end{align}
	where $M(\bar{X})$ is given by
	\begin{align}  \label{m}
		M[\bar{X}] = D^T C [\bar{X}] D
	\end{align}
	with $C  (\bar{X}) $ a diagonal matrix
	\begin{align}  \label{m1}
		C  [ \bar{X}]  =   \text{diag}\left(   ( |D \bar{X}|^2+\beta )^{-1/2}   \right).
	\end{align}
	Here, the smoothing parameter $\beta$ is added in order to make $\nabla \mathcal {J}[\bar{X}] $ differentiable. Having defined our regularisation functional, we now explain the minimisation approaches we use for solving the direct problem of QPAT.
	
	\subsection{TV regularisation using Alternating Direction Method of Multipliers (ADMM)} \label{sec:ADMM}
	Two major issues for minimisation of $\mathcal{F}$ is the nonlinearity of $\nabla \mathcal{J}[\bar{X}]$ and a loss of accuracy due to the smoothing parameter $\beta$. Note that a small value for $\beta$ may deteriorate the convergence \cite{Vogel}. One way for addressing these difficulties is to use a slack variable for shifting the gradient of $\|D \bar{X} \|_1$ out of the non-differentiable region and penalising the applied shift. To do this, the Augmented Lagrangian is introduced which,
	following \cite{Gao,Gao-b}, may be further rewritten as
	\begin{align} \label{FAm}
		\mathcal{F}_{A} (W, \bar{X})  =  \varrho \left(\frac{1}{2}\| D\bar{X}-W +  U_w   \|_2^2+ \nu  \|W \|_1    \right)+  \sum_{q=1}^{N_q}  \frac{1}{2}\|  \mathbb{H}_q[X(\bar{X})]- \hat{P}_q +U_{p,q}   \|_2^2,
	\end{align}
	where $\nu$ and $\varrho$ are constants and $U_w$ and $U_{p,q}$ are rescaled Lagrange multipliers. Minimisation of $\mathcal{F}$ is then accomplished by alternating minimisation of \eqref{FAm} in $W$ and $\bar{X}$, and updating of the Lagrange multipliers, using Algorithm~\ref{alg:43}.
	\begin{algorithm}
		\caption{Gradient-based Quasi-Newton method using ADMM}
		\label{alg:43}
		\begin{algorithmic}[1]
			\State Input: $\hat{P}_q,   \   (q\in \left\{ 1,..., N_q\right\})$
			\State Initialise: $X^0,W^0=0$
			\While {$ \| \nabla_{\bar{X}} \mathcal{F}_A \left( \bar{X}^k, W^k\right)  \| > Tol_\text{out} $}
			\State $W^{k+1}= \operatornamewithlimits{arg \ min}\limits_{W} \mathcal{F}_{A}  \left(  \bar{X}^k, W \right)    $  \label{al3l1}
			\State $\bar{X}^{k+1} = \operatornamewithlimits{arg \ min}\limits_{\bar{X}_l  \leqslant  \bar{X}   \leqslant \bar{X}_u} \mathcal{F}_{A}  \left(  \bar{X}, W^{k+1} \right)$  \label{al3l2}
			\State $U_w^{k+1}=U_w^k +D \bar{X}^{k+1}-W^{k+1}$
			\State $U_{p,q}^{k+1}=U_{p,q}^k +\mathbb{H}_q[ X^{k+1}]-\hat{P}_q     \     (\forall    q)  $
			\State Output: $X^*$
			\EndWhile
		\end{algorithmic}
	\end{algorithm}
	In Algorithm~\ref{alg:43}, $Tol_\text{out} $ is a terminating threshold, and $X^*$ denotes an optimal solution. The line \ref{al3l1} in Algorithm~\ref{alg:43} is a minimisation of the first term in \eqref{FAm} with respect to $W$, and can be calculated exactly using a scalar-wise \textit{Shrinkage formula} of the form
	\begin{align} \label{thre}
		W^{k+1}= \text{max} \left\{  |  D\bar{X}^k+U_w^k |- \nu, 0 \right\} \text{sgn}( D\bar{X}^k+U_w^k).
	\end{align}
	Additionally, the line \ref{al3l2} in this algorithm is a minimisation of $\mathcal{F}_{A}$ at $W^{k+1}$ with respect to $\bar{X}$, and is done using L-BFGS algorithm, as explained in section \ref{sm1}. \Scom{Note that we should replace $\epsilon$ by $\mathcal{F}_{A}$, and we hope this is not confusing for the reader}. To do this, the first-order derivative of $\mathcal{F}_{A}$ with respect to $\bar{X}$ at inner iteration $k$ is computed using
	\begin{align}
		\nabla_{\bar{X}} \mathcal{F}_{A}  = \varrho D^T \left( D\bar{X} - W+ U_w\right)  +  \sum_{q=1}^{N_q}  \frac{\partial X}{\partial \bar{X}}[\bar{X}]  \mathbb{J}_q^*[X]  \left(   \mathbb{H}_q[X]- \hat{P}_q  + U_{P_q} \right).
	\end{align}
	For this method we use a linear scaling
	\begin{equation} \label{linscale}
		\bar{X} = 
		\left [
		\begin{matrix}
			\bar{\kappa}\\\bar{\mu}
		\end{matrix}
		\right ]
		=
		\left [
		\begin{matrix}
			\frac{\hat{\kappa}}{\text{mean} (\hat{\kappa}^0)}\\
			\frac{\hat{\mu}}{\text{mean} (\hat{\mu}^0)}
		\end{matrix}
		\right ]
	\end{equation}
	where the dominators are the mean value of initial guesses $\hat{\kappa}^0$ and $\hat{\mu}^0$.
	
	\subsection{Linearised matrix-free Jacobian-based method} \label{secj}
	In this section, we explain two methods that we use for solving the problem \eqref{varr} as a sequence of linearised subproblems. Here, we set $X_l=- \infty$ and $X_u=+\infty$. This gives a non-constrained form of the problem. Our numerical results showed that enforcing bounds on the solutions is not required because of a good stability provided by these approaches. Additionally, we will make a balance between reconstruction of $\hat{\kappa}$ and $\hat{\mu}$ by using the logarithmic rescaling
	\begin{equation} \label{logscale}
		\bar{X}=  \log{\frac{X}{X^0}}
	\end{equation}
	which also implicitly enforces positivity on $\hat{\kappa}$ and $\hat{\mu}$. Here $X^0$ denotes an initial guess.

	Accordingly, given an iterate $\bar{X}^k$ and a point $\bar{X}$ in a neighborhood of $\bar{X}^k$, the forward operator $\mathbb{H}$ is linearised using an approximation
	\begin{align}   \label{jl}
		\mathbb{H}(\bar{X}) \approx \mathbb{H}(\bar{X}^k)+  \mathbb{J}[X^k] \frac{\partial X}{\partial \bar{X}}[\bar{X}^k] ( \bar{X} - \bar{X}^k).
	\end{align}
	Applying the approximation \eqref{jl} on the problem \eqref{varr} yields the minimisation problem
	\begin{align} \label{dsj}
		d^k = \operatornamewithlimits{argmin}_ {d} \ \frac{1}{2} \sum_q^{N_q} \norm{ \mathbb{J}_q [X^k] \frac{\partial X}{\partial \bar{X}}[\bar{X}^k]  d -\left (   \hat{P}_q- \mathbb{H}_q [X^k] \right) }^2,
	\end{align}
	where we have changed from $\bar{X}$ to $d=\bar{X} - \bar{X}^k $ in the minimisation.
	Note that we have used $k$ as a superscript in order to indicate a linearised subproblem (outer iteration), as opposed to a subscript in \eqref{sd} that indicates an inner iteration. The $k$-th linearised subproblem \eqref{dsj} gives a \textit{normal equation} in the form of
	\begin{align} \label{no}
		\mathcal{H}^k d^k =-\nabla \epsilon^k,
	\end{align}
	where $\nabla \epsilon^k$ and $\mathcal{H}^k$ are obtained from \eqref{gr} and \eqref{apj}, respectively. The normal equation \eqref{no} is a variant of \eqref{sd}, for which $\mathcal{H}^k$ is approximated using a \textit{Gauss-Newton} method. Here, to avoid a storage of $\mathbb{J}$, we solve each linearised subproblem \eqref{no} using a \textit{Krylov subspace} method in a matrix-free manner, for which implicit forms of operators $\mathbb{J}$ and $\mathbb{J}^*$ are used.
	
	As discussed in section \ref{sm}, \textit{Newton's methods} converge rapidly, but solving a normal equation with a high accuracy for each linearisation is very expensive. From a theoretical point of view, using \textit{Krylov} methods, the total number of iterations for reaching a minimiser is on the same order of the number of unknowns. Therefore, we solve \eqref{no} roughly using a loose stopping tolerance, i.e.,
	\begin{align}\label{mf}
		\mathcal{H}^k    \tilde{d}^k =- \nabla  \epsilon^k+ v_k, \  \ \|v^k\| /  \| \nabla  \epsilon^k \|   \leqslant  \eta^k,
	\end{align}
	where $\tilde{d}^k$ denotes a rough solution.
	It has been shown that under assumptions that $\mathcal{H}^k$ is symmetric and positive definite, the solutions $\tilde{d}^k$ are sufficiently small and $\eta^k<1$, the local convergence is guaranteed using a step size $\alpha^k=1$ (cf. \cite{Dembo}, section 2). A class of approaches that use \eqref{mf} for minimisation of \eqref{varr} are called \textit{Inexact Newton's methods} \cite{Dembo}. In the sequel, we use $d^k$, rather than $\tilde{d}^k$, for indicating a rough solution of \eqref{mf}, and we hope this is not confusing for the readers.
	
	\subsubsection{Lagged diffusivity (LD) method with priorconditioning} \label{sec:LD}
	As discussed above, for our direct QPAT problem, the Hessian matrices $\mathcal{H}^k$ are ill-conditioned. Therefore, a regularisation functional must be added in order to stabilise the problem. Our first approach for an inclusion of TV regularisation in \eqref{mf} is based on solving \eqref{obj_dis} via an iterative linearisation of an associated objective function $\mathcal{F}$. Strictly speaking, we first add the regularisation functional to an original nonlinear problem, and then the linearisations are applied to a regularised form of a nonlinear objective function. The k-th linearisation of the data fidelity term $\epsilon$ using \eqref{dsj}, together with  $\nabla \mathcal{J}[\bar{X}]$ defined by \eqref{jm}, \Scom{gives a TV regularised variant of \eqref{mf} in the form}
	\begin{align}  \label{unprior}
		\mathcal{H}^k d^k + \lambda M[\bar{X}]\bar{X}  =  - \nabla \epsilon^k .
	\end{align}
	Let us denote an initial guess for k-th subproblem by $\bar{X}_0^k$, which is calculated using the previous linearised subproblem $k-1$. One way for addressing the nonlinearity of $M[\bar{X}] \bar{X} $ (cf. \eqref{m})  is replacing $M[\bar{X}]$ by $M^k=M[\bar{X}_0^k]$. This gives a normal equation in the form
	\begin{align}  \label{unprior1}
		\left(	\mathcal{H}^k+ \lambda M^k \right) d^k =  - \nabla \epsilon^k - \lambda M^k \bar{X}_0^k.
	\end{align}
	Linearisation of $M$ using the above equation is called the \textit{Lagged Diffusivity} (LD) method \cite{Vogel}. (See \cite{Hannukainen} for an application of LD
	on the purely optical problem of QPAT.) Our method now is to follow \cite{Arridge-c} in order to convert \eqref{unprior1} into a similar problem in which the regularisation is obtained by early termination of an iterative method rather than tuning of parameter $\lambda$. First, since $M^k$ may only be positive semi-definite we approximate it by $M^k_\gamma = M^k + \gamma I$. Replacing $M^k$ by $M^k_\gamma$ in \eqref{unprior1} we next multiply by $(M^k_\gamma)^{-1}$ so that \eqref{unprior1} becomes
	\begin{equation}\label{unprior2}
		\left ( (M^k_\gamma)^{-1} \mathcal{H}^k + \lambda I \right ) d^k = - (M^k_\gamma)^{-1} \nabla \epsilon^k - \lambda\bar{X}_0^k
	\end{equation}
	Applying \textit{Krylov} methods for solving \eqref{unprior1}, the iterates lie \Scom{in a subspace} \cite{Arridge-c}
	\begin{align}
		\mathcal{K}^{(M^k_\gamma)^{-1} \mathcal{H}^k + \lambda I  }= \mathcal{K}^{(M^k_\gamma)^{-1} \mathcal{H}^k} =   \text{span}\left\{  - \left(  (M^k_\gamma)^{-1}  \mathcal{H}^k  \right)^{i\in    \left\{    0,...,   i_\text{max}^k  -1  \right\}   } \left ((M^k_\gamma)^{-1} \nabla \epsilon^k + \lambda \bar{X}_0^k \right )  \right\},
	\end{align}
	where $i_\text{max}^k$ is the maximum number of iterations for the \textit{Krylov} method. {\textcolor{red}{Our numerical experience shows that using small values for $\lambda$, the term $\lambda \bar{X}_0^k$ will have a small effect, and indeed we drop this in our method by taking $\lambda = 0$.}}
	
	Using this approach for applying regularisation on a normal equation is called \textit{priorconditioning}. Note that in this approach we adjust the regularisation by $i_\text{max}$. In contrast to an empirical choice for the regularisation parameter $\lambda$, which requires a recomputation of the problem, $i_{\text{max}}$ can be implicitly controlled by a stopping tolerance \cite{Arridge-c}.
	\noindent
	Using the above, our subproblem is to solve a \textit{priorconditioned} variant of \eqref{unprior1} in the form
	\begin{align}  \label{pre}
		({M_\gamma^k})^{-1}  \mathcal{H}^k   d^k =  -  ({M_\gamma^k})^{-1} \nabla \epsilon^k.
	\end{align}
	It turns out that \eqref{pre} provides a better convergence than \eqref{unprior1} since the structure of the prior is directly included in the Jacobian matrix \cite{Arridge-c}. Here, we solve \eqref{pre} using the \textit{Preconditioned Conjugate Gradient} (PCG) method, as outlined in Algorithm~\ref{alg:44}.
	
	\begin{algorithm}
		\caption{PCG algorithm for solving linearised subproblem $k$}
		\label{alg:44}
		\begin{algorithmic}[1]
			\State Input: ${\nabla \epsilon}^k,\mathcal{H}^k,M_\gamma^k$
			\State  Initialise: $ i=0,\bar{X}_0=0 $
			\State $\boldsymbol{r}_0 = - {\nabla \epsilon}^k$
			\State Solve $  {M_\gamma^k} \boldsymbol{z}_0 =   \boldsymbol{r}_0 $
			\State $\boldsymbol{d}_0 = \boldsymbol{z}_0$
			\While {$ {i<i_\text{max}} \ \cap \ \left(   i<i_m    \  \cup \      1- \frac{ \boldsymbol{r}_i^T  \boldsymbol{z}_i }{ \boldsymbol{r}_{i - i_m}^T  \boldsymbol{z}_{i - i_m} }> Tol_\text{in}  \right)  $}
			\State $\boldsymbol{\alpha}_i= \frac{\boldsymbol{r}_i^T  \boldsymbol{z}_i }{\boldsymbol{d}_i^T \mathcal{H}^k \boldsymbol{d}_i}$
			\State $\bar{X}_{i+1}= \bar{X}_i+ \boldsymbol{\alpha}_i \boldsymbol{d}_i$
			\State $\boldsymbol{r}_{i+1}= \boldsymbol{r}_i- \boldsymbol{\alpha}_i   \mathcal{H}^k  \boldsymbol{d}_i$
			\State Solve $ {M_\gamma^k}  \boldsymbol{z}_{i+1} =  \boldsymbol{r}_{i+1} $
			\State $\boldsymbol{\beta}_{i}= \frac{\boldsymbol{r}_{i+1}^T    \boldsymbol{z}_{i+1}} {\boldsymbol{r}_i^T    \boldsymbol{z}_i}$
			\State $\boldsymbol{d}_{i+1}=  \boldsymbol{z}_{i+1}+\boldsymbol{\beta}_{i}  \boldsymbol{d}_i $
			\EndWhile
			\State Output: $\bar{X}_*$
		\end{algorithmic}
	\end{algorithm}
	
	\noindent
	We will terminate Algorithm~\ref{alg:44} if $i>i_m$ with $i_m$ a user-adjusted number of inner iterations, and a relative reduction in $\boldsymbol{r}_i^T \boldsymbol{z}_i $ during $i_m$ inner iterations becomes less than a user-adjusted threshold $Tol_\text{in}$. Also, the PCG algorithm is unconditionally stopped whenever $i> i_\text{max}$.
	
	{\textcolor{red}{
			\begin{remark}
				Each inner iteration of the PCG loop involves an implicit inversion of the sparse matrix $M_\gamma^k$ (See \cite{Arridge-c}). For the direct QPAT, the cost of an inversion of $M_\gamma^k$ is negligible, compared to an implementation of the Jacobian and its transpose. 
	\end{remark}}}

	Using the LD method, together with a logarithmic scaling, our inexact Newton algorithm is outlined in Algorithm~\ref{alg:45}.
	\begin{algorithm}
		\caption{Inexact Newton method using LD}
		\label{alg:45}
		\begin{algorithmic}[1]
			\State Input: $\hat{P}_q  \   (q\in \left\{ 1,..., N_q\right\})$
			\State Initialise: $k=0, X^0$
			\While {$ k=0   \   \cup   \     1- \frac{\epsilon^k}{\epsilon^{k-1}}> Tol_\text{out}  $}
			\State Apply linearisation on \eqref{obj_dis}
			\State Compute $d^k$ from \eqref{pre} using Algorithm~\ref{alg:44}
			\State Compute $\bar{X}^{k+1}$ using \eqref{sdm} and $\alpha^k=1$
			\State $X^{k+1} = X^{0} e^{\bar{X}^{k+1}}$
			\EndWhile
			\State Output: $X^*$
		\end{algorithmic}
	\end{algorithm}

	\subsubsection{Primal-Dual Interior-Point-Method (PD-IPM)} \label{sec:PDIPM}
	A technique for linearisation of $M$ was developed using a primal-dual method, and was shown to give better convergence than the LD method, especially  for small values of $\beta$ \cite{Chan}. In \cite{Chan}, the PD-IPM technique was used for enforcing TV regularisation when inverting a linear blurring operator. \Scom{In contrast to our first approach using the LD method, here we first linearise the data fidelity function, and then add a TV regularisation function to each linearised subproblem using the PD-IPM approach.} Using this, we iteratively solve a TV regularised variant of the linearised subproblem \eqref{dsj} in the form
	\begin{align} \label{dsjj}
		d^k = \operatornamewithlimits{argmin}_ {d} \ \frac{1}{2} \sum_{q=1}^{N_q} \norm{ \mathbb{J}_q [X_0^k]  \frac{\partial X}{\partial \bar{X}}[\bar{X}_0^k] d -\left (   \hat{P}_q- \mathbb{H}_q [X_0^k] \right) }_{\mathcal{L}^2}^2  +   \lambda \|D d \|_{\mathcal{L}^1}.
	\end{align}
	The main idea for linearisation of the nonlinear $M$ (cf. \eqref{m} and \eqref{m1}) using the PD-IPM approach is introducing a dual parameter
	$\chi= C[d] D d \in \mathbb{R}^{2N_l}$. This gives a system of coupled nonlinear PDEs for subproblem $k$ in the form
	\begin{align}
		\begin{split}
			\lambda D^T \chi + \mathcal{H} d - \sum_{q=1}^{N_q} \frac{\partial X}{\partial \bar{X}}[\bar{X}_0] \mathbb{J}_q^* [X_0] \left (   \hat{P}_q- \mathbb{H}_q [X_0] \right)   =g(\chi,d)=0\\
			C^{-1}[d]\chi -Dd = f(\chi,d)=0,
		\end{split}
	\end{align}
	where we have removed the superscripts indicating subproblem $k$ for brevity.
	A linearisation of this system with respect to $(\chi,d)$ gives
	\begin{align}  \label{nsw}
		\begin{split}
			\begin{bmatrix}
				C^{-1}[d]   &-\Big(I-C[d] \chi {(Dd)}^T \Big) D \\
				\lambda D^T  &\mathcal{H}
			\end{bmatrix}
			\begin{bmatrix}
				\delta \chi \\
				\delta d
			\end{bmatrix}=
			\begin{bmatrix}
				-f(\chi,d)\\
				-g(\chi,d)
			\end{bmatrix}
		\end{split}
	\end{align}
	If we make the replacement $\xi = C[d] D d$, then above linearised system gives decoupled equations
	\begin{align}  \label{du}
		\begin{split}
			&\bigg[  \lambda   D^T C[d_{k'}] \Big(I-C[d_{k'}] \chi_{k'} ( Dd_{k'} )^T \Big) D + \mathcal{H} \bigg]\delta d_{k'}	\\
			=& - \lambda M(d_{k'}) + \sum_{q=1}^{N_q} \frac{\partial X}{\partial \bar{X}}[\bar{X}_0] \mathbb{J}_q^* [X_0] \left (   \hat{P}_q- \mathbb{H}_q [X_0] - \mathbb{J}_q [X_0]  \frac{\partial X}{\partial \bar{X}}[\bar{X}_0] d_{k'} \right)
		\end{split}
	\end{align}
	and
	\begin{align}   \label{dw}
		\delta \chi_{k'} = C[d_{k'}] \Big( I-C[d_{k'}] \chi_{k'} ( Dd_{k'} )^T \Big) D  \delta d_{k'}  -\chi_{k'} + C[d_{k'}]D  d_{k'},
	\end{align}
	where the subscript $k'$ indicates an inner \Scom{sub-subproblem}.
	The same as the LD method, we solve a priorconditioned form of \eqref{du} with $\lambda=0$, i.e.,
	\begin{align}  \label{pre1}
		\tilde{M}_{\gamma,k'}^{-1}  \mathcal{H}   \delta d_{k'}  =  -\tilde{M}_{\gamma,k'}^{-1} \nabla \tilde \epsilon_{k'},
	\end{align}
	where
	\begin{align}
		\tilde{M}_{\gamma,k'}= D^T C[d_{k'}] \Big(I-C[d_{k'}] \chi_{k'} (Dd_{k'})^T \Big)D+\gamma I,
	\end{align}
	{\textcolor{red}{and $-\nabla \tilde \epsilon_{k'}$ is actually the second term in the right-hand-side of \eqref{du}.}}
	Here, the subscript $k'$ indicates the fact that for each linearised subproblem $k$ (superscript), we solve a sequence of sub-subproblems \eqref{dw} and \eqref{pre1} using an update of $d_{k'}$, $M_{\gamma,k'}$ and $\nabla \tilde{\epsilon}_{k'}$. 
	The developed inexact Newton method using a TV regularisation based on the PD-IPM approach is outlined in Algorithm~\ref{alg:46}. We use step sizes $\alpha_{k'}=1$ and $\alpha^{k}=1$ for all $k'$ and $k$, respectively. The step size $s_{k'}$ is described below the algorithm (see \eqref{sstep}).
	
	\begin{algorithm}
		\caption{Inexact Newton method using PD-IPM}
		\label{alg:46}
		\begin{algorithmic}[1]
			\State Input: $\hat{P}_q  \   (q\in \left\{ 1,..., N_q\right\})$
			\State Initialise: $k=0, X^0$
			\While {$ k=0 \  \cup  \  \ 1- \frac{\epsilon^k}{\epsilon^{k-1}}> Tol_\text{out}     $}
			\State Apply linearisation and derive the objective function in \eqref{dsjj}
			\State Initialise: $\bar{X}_0^k, d_0=0, \chi_0=0$
			\While{$  k'< k'_{\text{max}}   \  \cap \  \left(    k'=0  \  \cup \  1- \frac{\boldsymbol{r}_{k',*}}{\boldsymbol{r}_{k'-1,*}}> Tol_\text{med} \right) $}  \label{stal6}
			\State Compute $\delta d_{k'}$ from \eqref{pre1} using Algorithm~\ref{alg:44}
			\State $d_{k'+1}=d_{k'}+ \alpha_{k'}\delta d_{k'}$
			\State Compute $\delta \chi_{k'}$ using \eqref{dw}
			\State $\chi_{k'+1}=\chi_{k'}+ s_{k'}\delta \chi_{k'}$
			\EndWhile
			\State Compute $\bar{X}^{k+1}$ using \eqref{sdm} and $\alpha^k=1$
			\State $X^{k+1} = X^{0} e^{ \bar{X}^{k+1} }$
			\EndWhile
			\State Output: $X^*$
		\end{algorithmic}
	\end{algorithm}
	Here, each iteration $k'$ amounts to solving a PCG loop. The optimal $\boldsymbol{r}$ provided by each PCG loop $k'$ is denoted by $\boldsymbol{r}_{k',*}$. Using this, we terminate each outer subproblem $k$ using a stopping criterion given in line \ref{stal6} in Algorithm~\ref{alg:46}. This stopping criterion uses a stopping threshold $Tol_{\text{med}}$.
	
	Additionally, following \cite{Borsic}, we choose $s_{k'}$ using a \textit{step length rule}. Using this approach, we choose
	\begin{align} \label{sstep}
		s_{k'}= \text{min}\left(1,\varphi_*\right) \delta \chi_{k'}.
	\end{align}
	Here, $\varphi_*$ is the largest $\varphi$ that satisfies a feasibility condition
	\begin{align}
		\left(|\chi_{k'}+ \varphi  \delta \chi_{k'}|\right)_j\leqslant 1, \  \forall j=1,...,2N_l,
	\end{align}
	where $j$ denotes the index of components of $\chi$.

	\section{Numerical results}
	The TV regularised minimisation approaches that have been explained in section \ref{tvr}, i.e., ADMM, LD and PD-IPM, were used for a simultaneous reconstruction of images of optical absorption coefficient $\mu$ and diffusion coefficient $\kappa$ for 2D and 3D phantoms.
	
	\subsection{2D phantom}  \label{numr2d}
	The 2D simulation was performed on a square domain $[-5,+4.92] \times [-5,+4.92]$ $\text{mm}^2$. 
	
	\subsubsection{Optical excitation}
	Four different optical excitation patterns were used, i.e., $N_q=4$. For each optical excitation $q$, we used a discretisation of an inward directed diffuse boundary current $I_{s,q} \ (J/\text{mm})$ that obeys
	\begin{align}  \label{isn}
		I_{s,q}(r)= 
		\begin{cases}
			1, & r \in \iota_q\\
			0, & r \in \partial \Omega \backslash \iota_q, 
		\end{cases}	
	\end{align}
	where $\iota_q \in \partial \Omega $ denotes the source position for optical excitation $q$, and was set each side of the square for each optical excitation (cf. the second line in equation \eqref{DA}.) 
	
	\subsubsection{Discretisation for data generation}
	For generation of time series of boundary data, the square domain was discretised using a grid with $128 \times 128$ nodes and an even separation distance of $7.81 \times 10^{-2}$mm along both Cartesian coordiantes. For the optical portion of the problem, a triangulation was applied so that each two finite elements form a pixel, and the centre of the pixel matches an associated node on the acoustic grid. For the acoustic portion of the problem, to mitigate wave wrapping \cite{Tabei}, a perfectly matched layer (PML) having a thickness of $20$ grid points and a maximum attenuation coefficient of $2$ nepers per grid point was added to each side of the grid. The propagated wavefield was detected in $1017$ time steps using $158$ detectors that are equidistantly placed on the left and top sides of the computational grid, as shown in figure~\ref{fig:41}. A $30 \ \text{dB}$ Additive White Gaussian Noise (AWGN) was then added to the simulated data.

	\subsubsection{Discretisation for image reconstruction}
	To avoid an inverse crime for discretisation, the image reconstruction was done using a computational grid made up of $80 \times 80$ nodes with an even separation distance of $1.25 \times 10^{-1}$mm.
	
	\subsubsection{Acoustic properties}  \label{acous2d}
	To the best of our knowledge, our manuscript reports the first results on the direct QPAT for realistic acoustic media, for which acoustic characteristics of tissue media such as heterogeneity and attenuation are taken into account. 
	
	In addition, PAT and QPAT use an assumption that the acoustic properties of the medium are known. This assumption does not hold in practical cases. For example, it has been shown that the acoustic properties of the breast vary up to 15 $\%$.
	These variations are often not exactly known for reconstruction. As a result, using the same acoustic properties for data generation and image reconstruction may be an inverse crime. 
	To avoid this, for data generation, we corrupted the acoustic properties of the medium with $30 \ \text{dB}$ AWGN noise.
	Figures \ref{F1a} and \ref{F1b} show the contaminated distributions of sound speed ($c_0$) and ambient density ($\rho_0$) for data generation, respectively. The minimal and maximal values of these maps are given in the top row of Table \ref{tab1}. Using these values, the computational grid for data generation supports a maximal frequency up to $7.223$ MHz. (See the \textit{k-Wave} manual \cite{Treeby-b}.) Also, in figures \ref{F1a} and \ref{F1b}, the position of detectors is shown by the black circles matching the left and top sides of the grid.

	For image reconstruction, we used the clean maps that are shown in figures \ref{F1c} and \ref{F1d}. The minimal and maximal values of these maps are given in the bottom row of Table \ref{tab1}. Using these acoustic maps,  
	the grid for image reconstruction supports a maximal frequency up to $5.101$ MHz. From Table \ref{tab1}, the incorporated noise has provided a 10-15$\%$ relative discrepancy between the acoustic properties used for data generation and image reconstruction. 
	
	\begin{table}\centering
		\caption {The minimal and maximal values for acoustic properties of the 2D phantom.} \label{tab1} 
		\begin{tabular}{ccccc}
			\hline
			& \multicolumn{2}{c}{\text{$c_0 (\text{m} \text{s}^{-1}$) }}            & \multicolumn{2}{c}{$\rho_0 (\text{kg} \text{m}^{-3}$)}                 \\ \cline{2-5} 
			& min                   & max                    & min                   & max                   \\ \hline
			Data generation      & $1.129 \times 10^3$ & $1.902 \times 10^3$ & $0.620 \times 10^3 $ & $1.390 \times 10^3$ \\ 
			Image reconstruction  & $1.276 \times 10^3$ & $1.725 \times 10^3$ & $0.750 \times 10^3 $ & $1.250 \times 10^3$                                           \\ \hline
		\end{tabular}
	\end{table}

	Furthermore, the acoustic medium was assumed attenuating, where acoustic absorption and dispersion follow a frequency power law \cite{Treeby}. Accordingly, we used a constant attenuation coefficient $\alpha_0= 0.75 \ \text{dB} \hspace{0.1cm}
	\text{MHz}^{-y} \hspace{0.1cm} \text{cm}^{-1}$ and an exponent factor $y= 1.5$ for both data generation and image reconstruction (cf. \eqref{Twelve}). These values were chosen so that they approximately simulate the acoustic attenuation properties of the breast.

	\begin{figure}  \centering
		{\subfigure[]{\includegraphics[width=.49\textwidth]{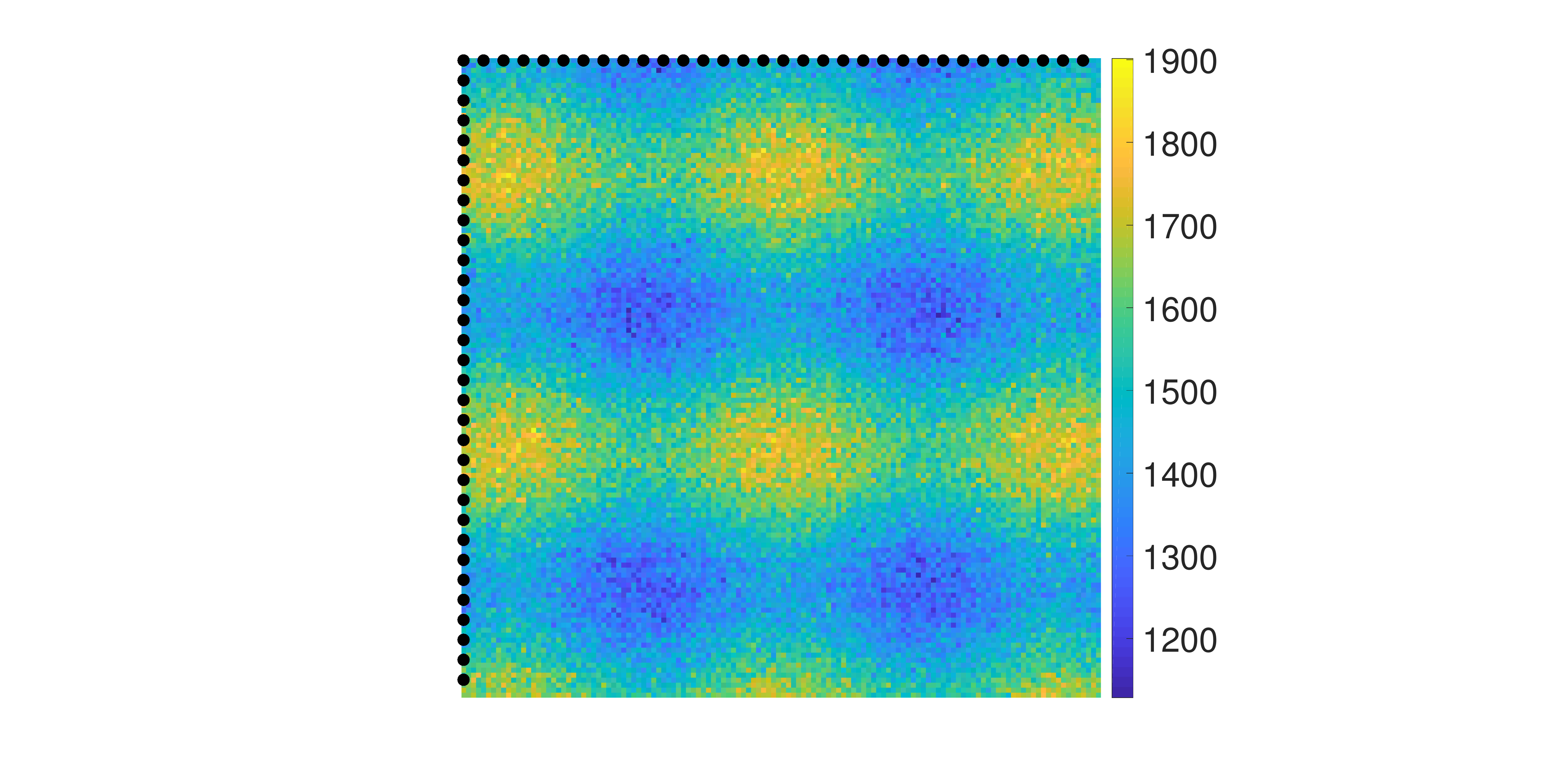}\label{F1a}}
			\subfigure[]{\includegraphics[width=.49\textwidth]{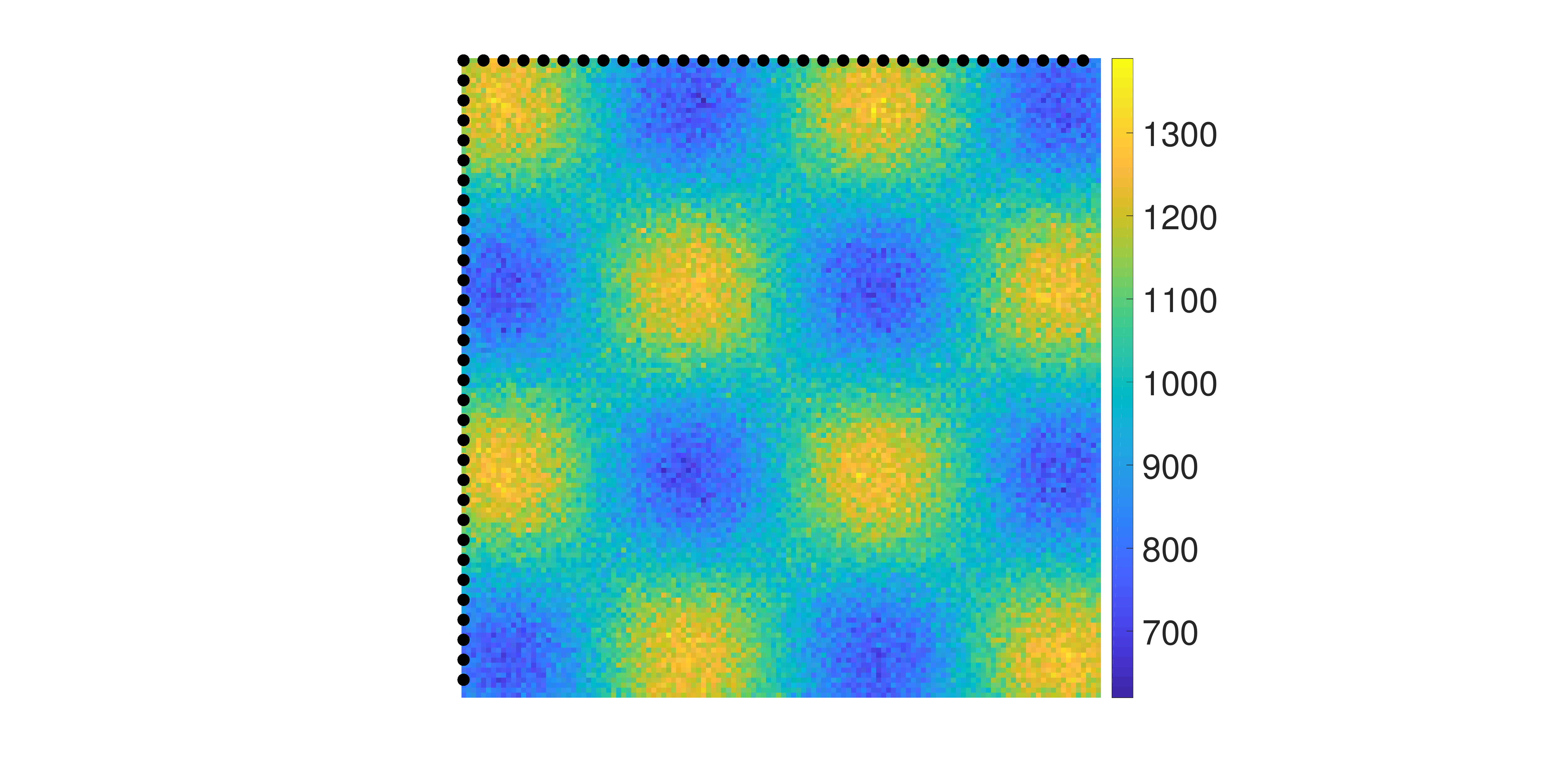}\label{F1b}}
			\subfigure[]{\includegraphics[width=.49\textwidth]{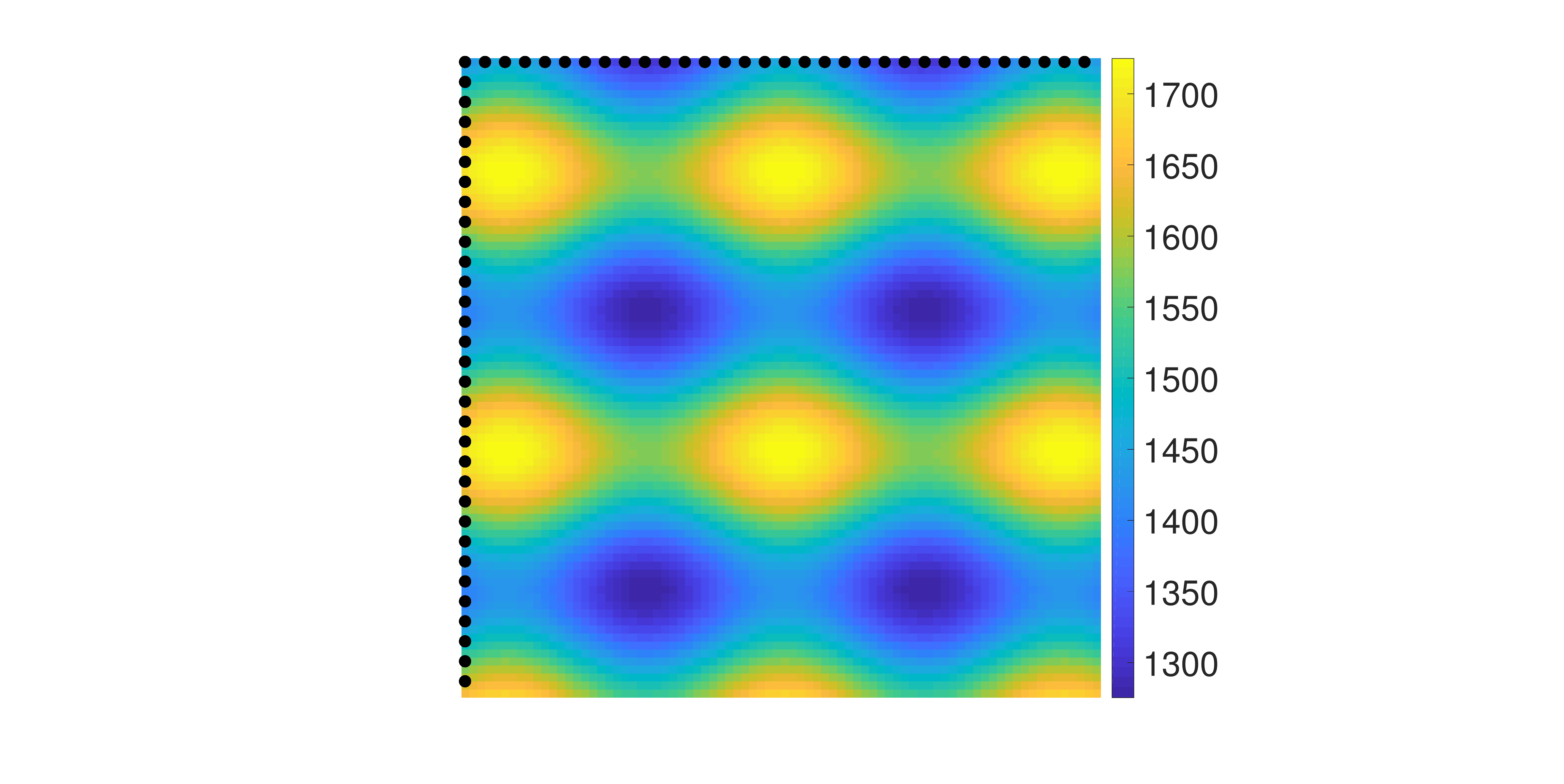}\label{F1c}}
			\subfigure[]{\includegraphics[width=.49\textwidth]{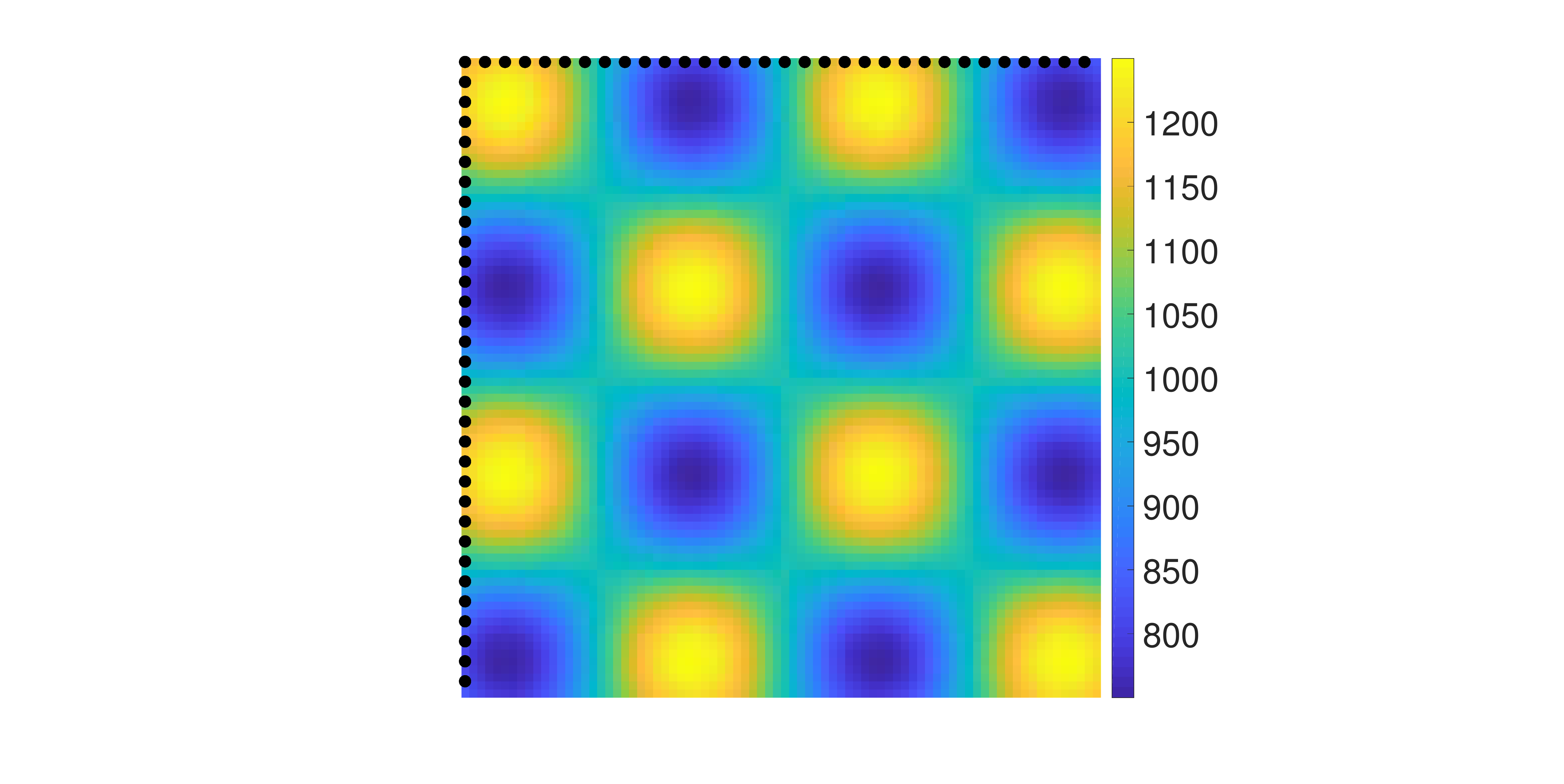}\label{F1d}}
		}
		\caption{Acoustic properties for 2D case. Data generation: (a) $c_0$ (b) $\rho_0$, and image reconstruction: (c) $c_0$ (d) $\rho_0$.}
		\label{fig:41}
	\end{figure}

	\subsubsection{Optical phantom}
	We simulate the distributions of optical coefficients so that they follow the optical properties of soft tissues in the sense that they often possess a broad range of values for the optical coefficients. Our numerical experience showed that this is a challenge for image reconstruction, although this issue has been neglected in many of studies of QPAT. The reader is referred to \cite{Pulkkinen-b} for a study on the direct QPAT, in which this issue has been considered. Accordingly, we simulate a distribution for optical absorption coefficient with 20 values within a range $\mu \in [0.025,0.325] \ \text{mm}^{-1}$ with a background $0.075\ \text{mm}^{-1}$. Also, a distribution for the diffusion coefficient is simulated so that it has 6 values within a range $\kappa \in [0.2,0.4] \ \text{mm}^{-1}$ with a background $0.3\ \text{mm}^{-1}$. Figures \ref{F2a} and \ref{F2b} show the map for $\mu$ and $\kappa$, respectively.

	\subsubsection{Image reconstruction}  \label{imr}
	We initialised all algorithms using values $1.2$ times more than the mean of optical coefficients for the associated phantoms. 
	
	\textit{ADMM.} The ADMM approach (cf. Algorithm~\ref{alg:43}) was applied for a simultaneous reconstruction of $\mu$ and $\kappa$. Since ADMM is our benchmark method, the associated parameters were chosen very carefully in order to obtain the best possible image. The line \ref{al3l2} in Algorithm~\ref{alg:43} is a minimisation of $\mathcal{F}_A$ (cf. \eqref{FAm}) with respect to $\bar{X}$ using the L-BFGS algorithm given in Algorithm~\ref{alg:41}. The L-BFGS algorithm uses a backtracking line search given in Algorithm~\ref{alg:42} with $c_1=1\times 10^{-4}$, $c_2=0.9$ and $\tau=0.25$. Also, we set $\varrho=1$  (cf.\eqref{FAm}), and $\nu=1 \times 10^{-3}$ for the \textit{Shrinkage} operator \eqref{thre}. We terminated the ADMM algorithm using $Tol_{\text{out}}=1 \times 10^{-2}$ (cf. Algorithm~\ref{alg:43}). 
	
	\textit{LD.} The LD method was applied using Algorithm~\ref{alg:45}. Using this algorithm, each linearised subproblem is solved using a PCG loop (cf. Algorithm~\ref{alg:44}) by setting $i_\text{max}=30$, $i_m=5$, and $Tol_\text{in}=0$.
	The latter parameter implies that we terminate each PCG loop, if 
	\begin{align}
		{i>i_\text{max}}  \ \cup  \ \left(  i>i_m   \  \cap \  \boldsymbol{r}_i^T  \boldsymbol{z}_i > \boldsymbol{r}_{i - i_m}^T  \boldsymbol{z}_{i - i_m}   \right). 
	\end{align}
	Note that we observed a nonmonotone convergence for iterates of each PCG loop (inner iterations) in regions close to an optimal solution $X^*$, but the sequence $\epsilon (X^k)$ always monotonically converged to $\epsilon (X^*)$ using $\alpha^k=1$. (We suggest using an Armijo condition for the outer iterations, although we observed that a nonmonotonic reduction in iterates of the PCG loop associated with outer iteration $k$ is sufficient for providing a descent search direction, i.e., $\epsilon(X^{k+1}) <\epsilon(X^k) $).
	The TV preconditioner was applied using $\gamma=1 \times 10^{-6}$ and $\beta= 2 \times 10^{-5}$. Our LD algorithm was stopped using $Tol_\text{out}=1 \times 10^{-3}$ (cf. Algorithm~\ref{alg:45}).
	
	\textit{PD-IPM.} The PD-IPM technique was applied using Algorithm~\ref{alg:46}. Because of applying two layers of linearisation, each of the outer linearised problems are solved using a sequence of inner linearised subproblems. For solving a normal equation associated with each inner linearised subproblem, we terminated each PCG algorithm using the same parameters as in the LD method. The TV preconditioner was applied using $\gamma= 1 \times 10^{-6}$ and $\beta= 1 \times 10^{-6}$.
	
	We terminated each outer linearised suproblem using a threshold $Tol_{med}= 1 \times 10^{-3}$ and $k'_\text{max}=20$ (cf. Algorithm~\ref{alg:46}). Our PD-IPM was terminated using a stopping threshold $Tol_{out}= 1 \times 10^{-3}$.
	
	\subsubsection{Evaluation of image reconstruction}
	The criterion that we use for measuring the convergence of sequence $X^k$ to a ground truth image (phantom) is
	Relative Error (RE), which is calculated as
	\begin{align}
		RE (u^k) = 100 \times \frac{\|u^k -u^{phantom} \|}{\|u^{phantom}\|}.
	\end{align}
	Here, $u^k$ is the solution for either $\kappa$ or $\mu$ at outer iteration $k$ that is interpolated back to the grid for data generation (phantom). Also, the superscript $phantom$ indicates the distribution of optical coeffcients for the phantom.
	
	We also consider $\epsilon(X^k)$ as the second criterion for convergence (cf. \eqref{varr}).

	\subsubsection{Observations}
	
	\textit{ADMM.} Using the parameters given in section \ref{imr}, the ADMM algorithm was stopped after outer iteration 6. The final reconstructed images for the optical absorption coefficient $\mu$ and diffusion coefficient $\kappa$ are shown in figures \ref{F2c} and \ref{F2d}, respectively. Figure \ref{F3a} shows the RE of sequence computed by the ADMM algorithm at outer iterations. 
	
	\textit{LD.} For the LD algorithm, the associated stopping criterion was satisfied after outer iteration 30. Figures \ref{F2e} and \ref{F2f} show the final reconstructed images for $\mu$ and $\kappa$, respectively. Figure \ref{F3b} shows the RE of the iterates provided by LD for outer iterations $k$. The computed values for $\epsilon(X^k)$ are shown in figure \ref{F4a}. This figure is shown from an enlarged view around the optimal solution in figure \ref{F4b}. As shown in these figures, $\epsilon$ monotonically converges to a minimiser for all outer iterations using our choice for the step size $\alpha^k=1 \ (\forall k)$.

	\textit{PD-IPM.} The stopping criterion for the PD-IPM algorithm was satisfied after 4 outer iterations. The final reconstructed images for $\mu$ and $\kappa$ are shown in figures \ref{F2g} and \ref{F2h}, respectively. Figure \ref{F3c} shows the RE of solutions (optical coefficients) computed by the PD-IPM algorithm for outer iterations $k$. Figure \ref{F4c} shows the obtained values for $\epsilon(X^k)$. This figure is shown from an enlarged view around the optimal solution in figure \ref{F4d}. As shown in these figures, our choices for the step sizes associated with outer (resp. inner) subproblems, which are $\alpha^k=1$ (resp. $\alpha_{k'}=1$), provided a monotonic reduction for values of $\epsilon^k$ (resp. $\epsilon_{k'}$) for all iterations.
	
	\begin{remark}
		Using PD-IPM, for both outer and inner linearised subproblems, the first step is a compuation of associated $\epsilon$ and the gradient $\nabla \epsilon$ (cf. Algorithm~\ref{alg:46}). As a result, using a line search for compuation of $\alpha^k$ and $\alpha_{k'}$, e.g., a backtracking line search using Wolfe conditions, is straightforward, and does not impose additional computational cost. However, our numerical experience showed that a reduction in $\boldsymbol{r}_i^T \boldsymbol{z}_i$ provided by the PCG loops (cf. Algorithm~\ref{alg:44}) is sufficient for a monotonic reduction of the objective function without using a line search.
	\end{remark}

	\begin{figure}    \centering
		{\subfigure[]{\includegraphics[width=.49\textwidth]{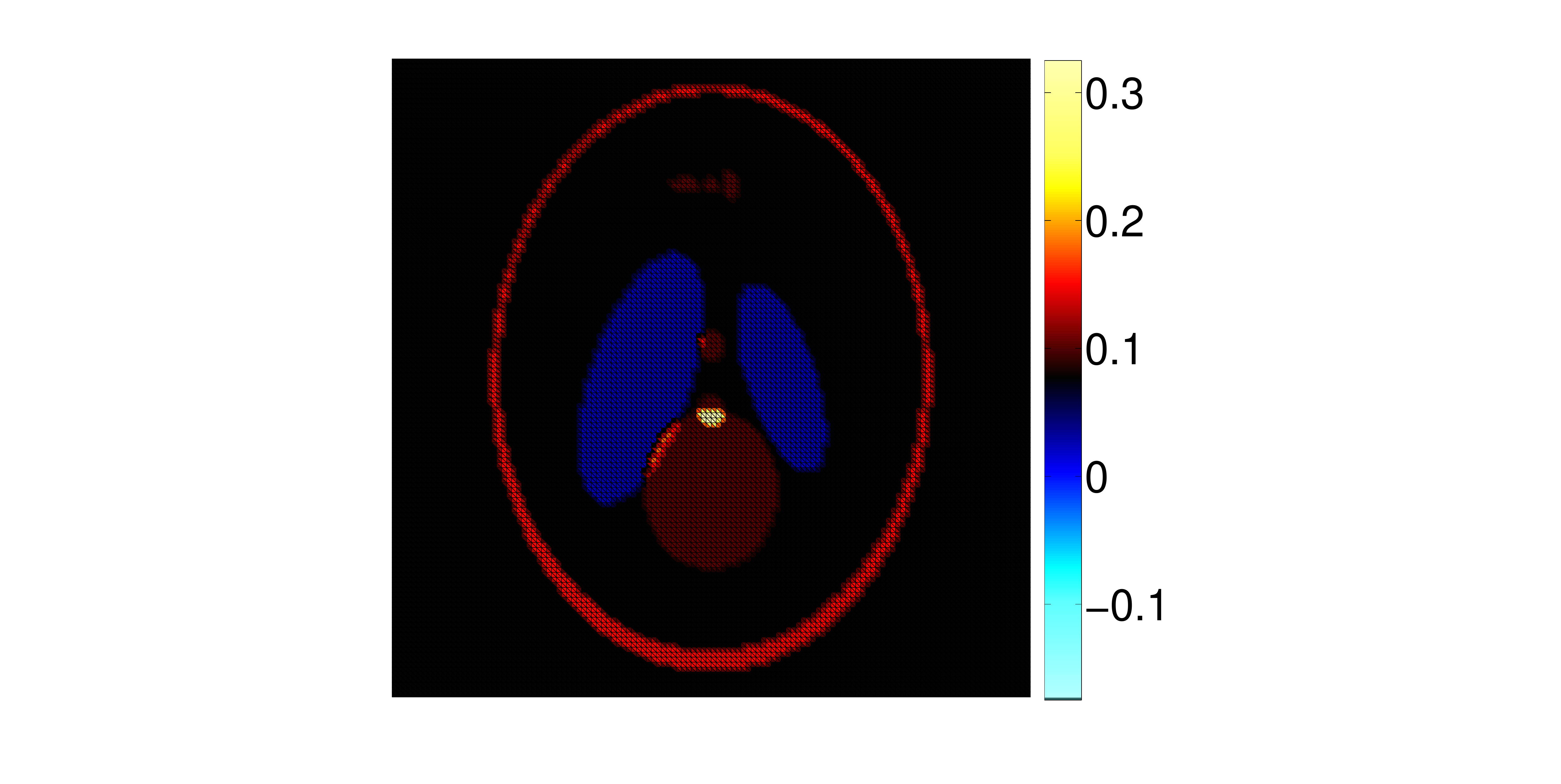}\label{F2a}}
			\subfigure[]{\includegraphics[width=.49\textwidth]{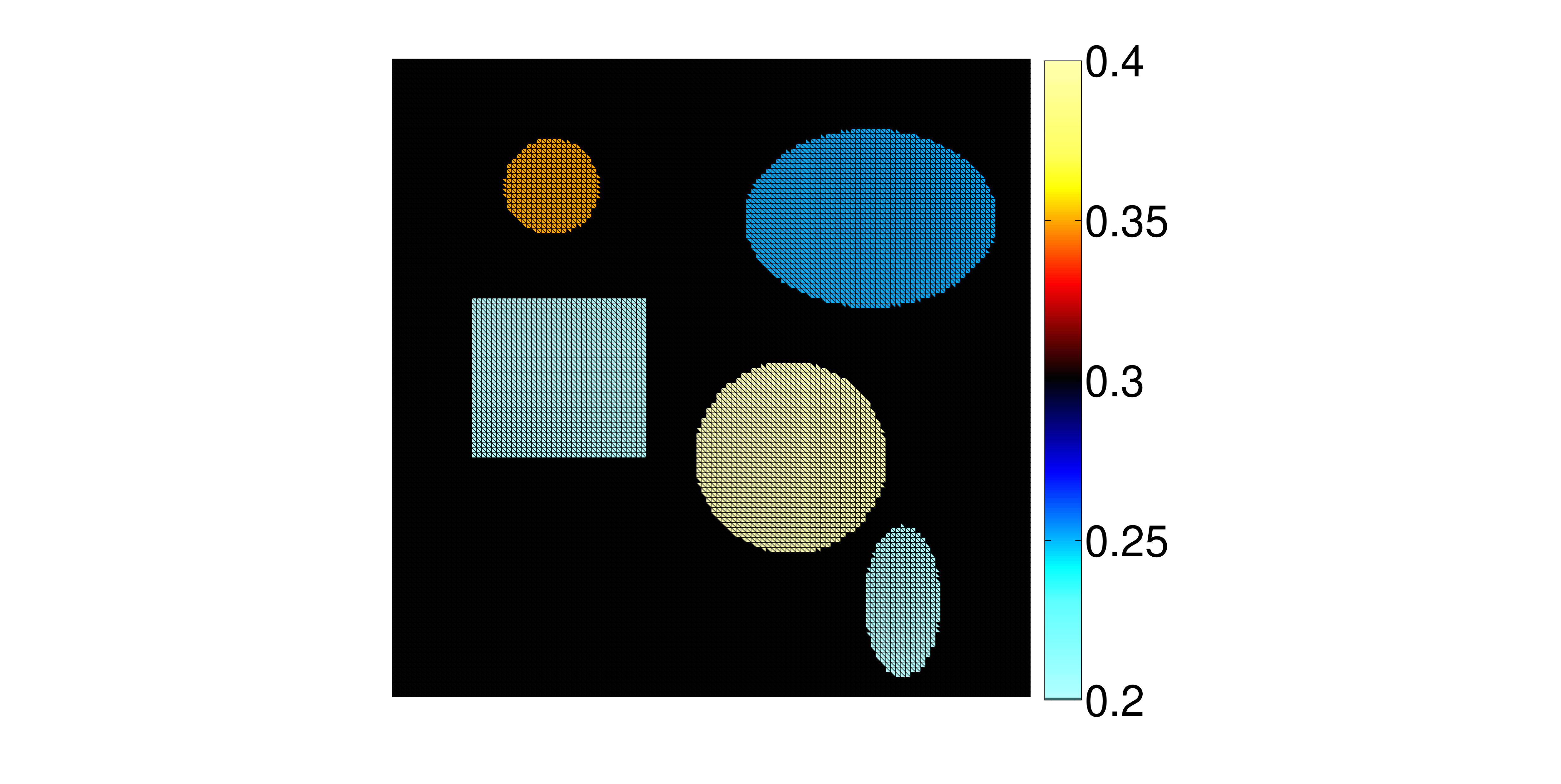}\label{F2b}}
			\subfigure[]{\includegraphics[width=.49\textwidth]{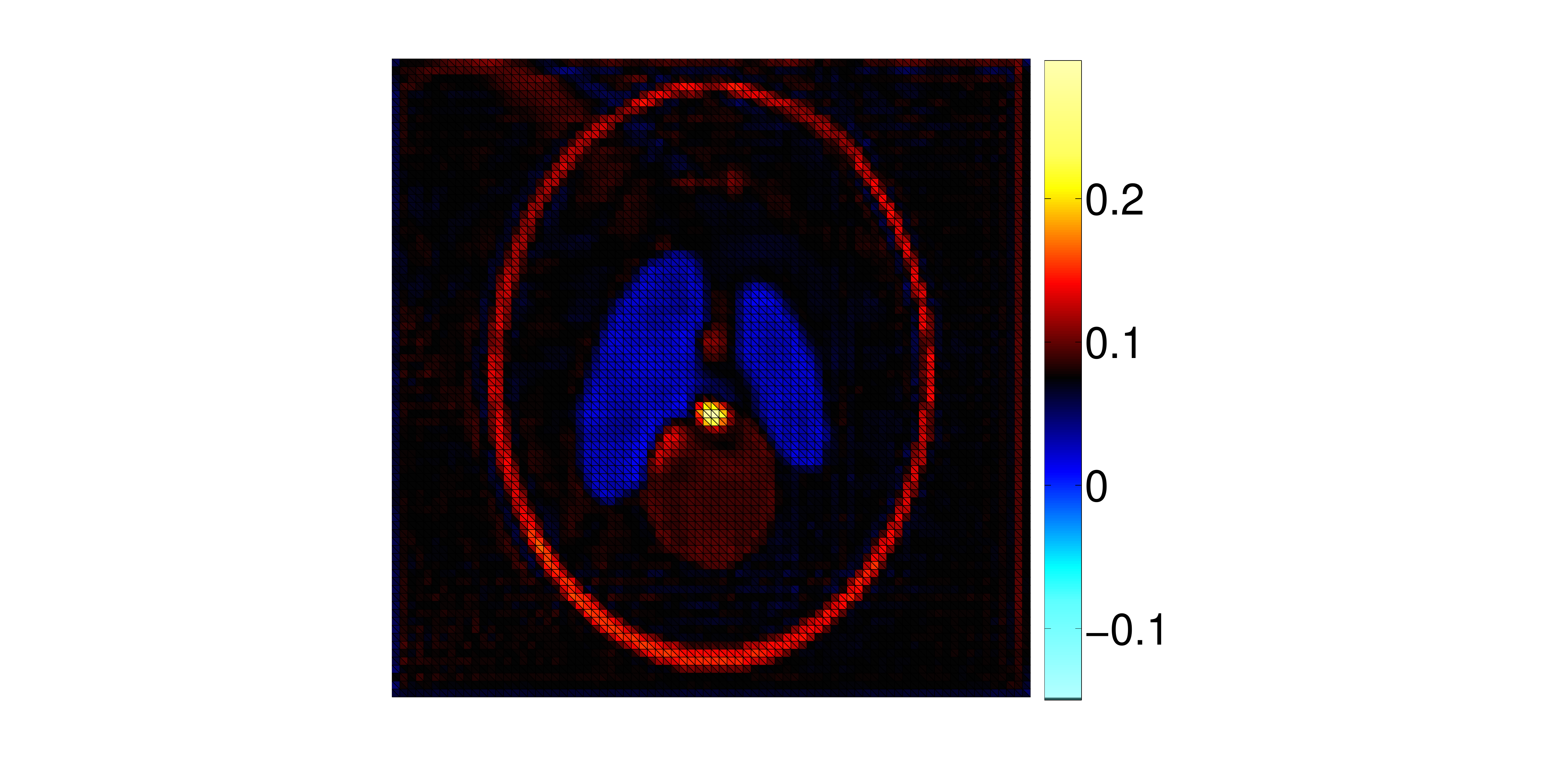}\label{F2c}}
			\subfigure[]{\includegraphics[width=.49\textwidth]{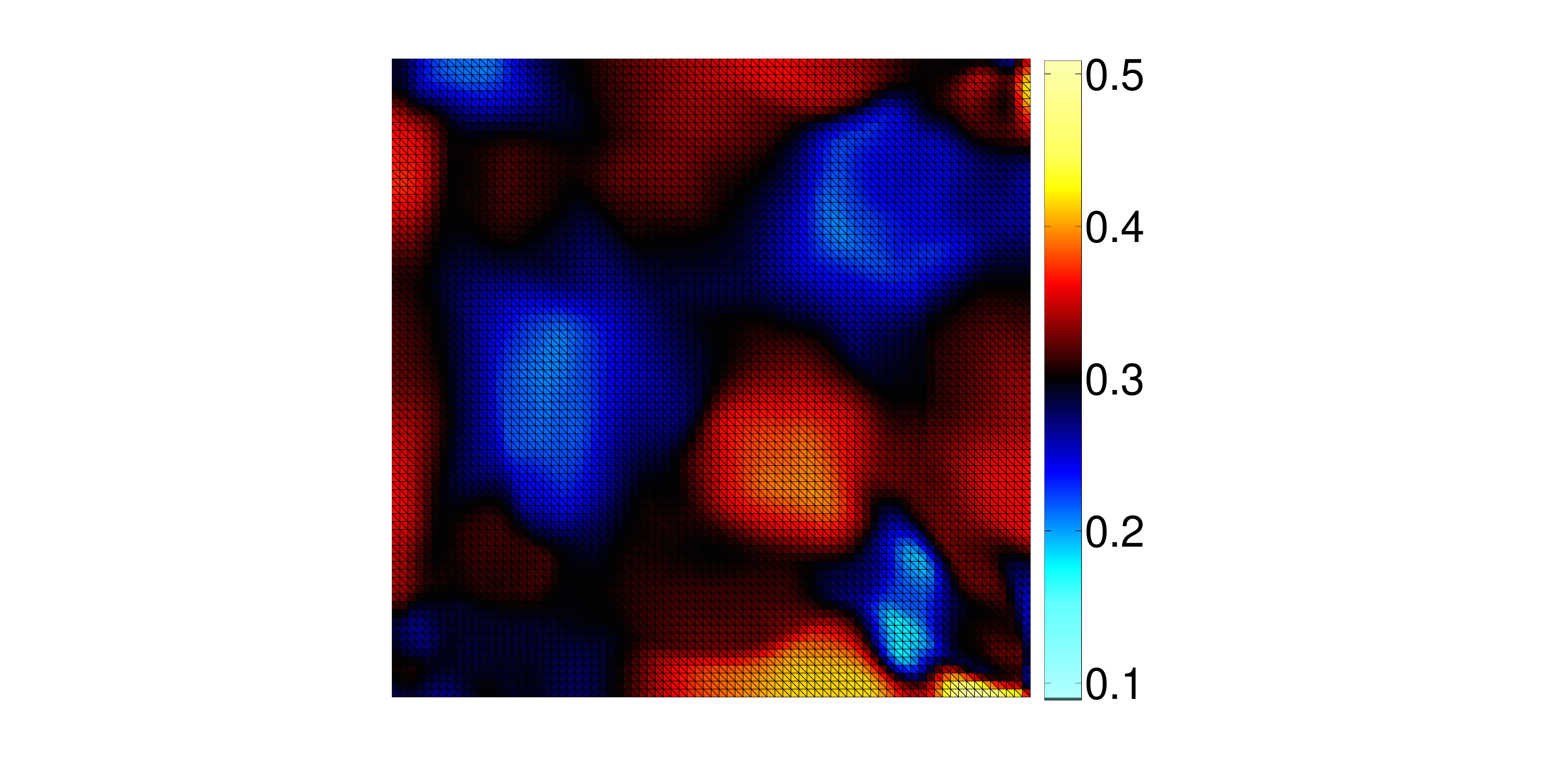}\label{F2d}}
			\subfigure[]{\includegraphics[width=.49\textwidth]{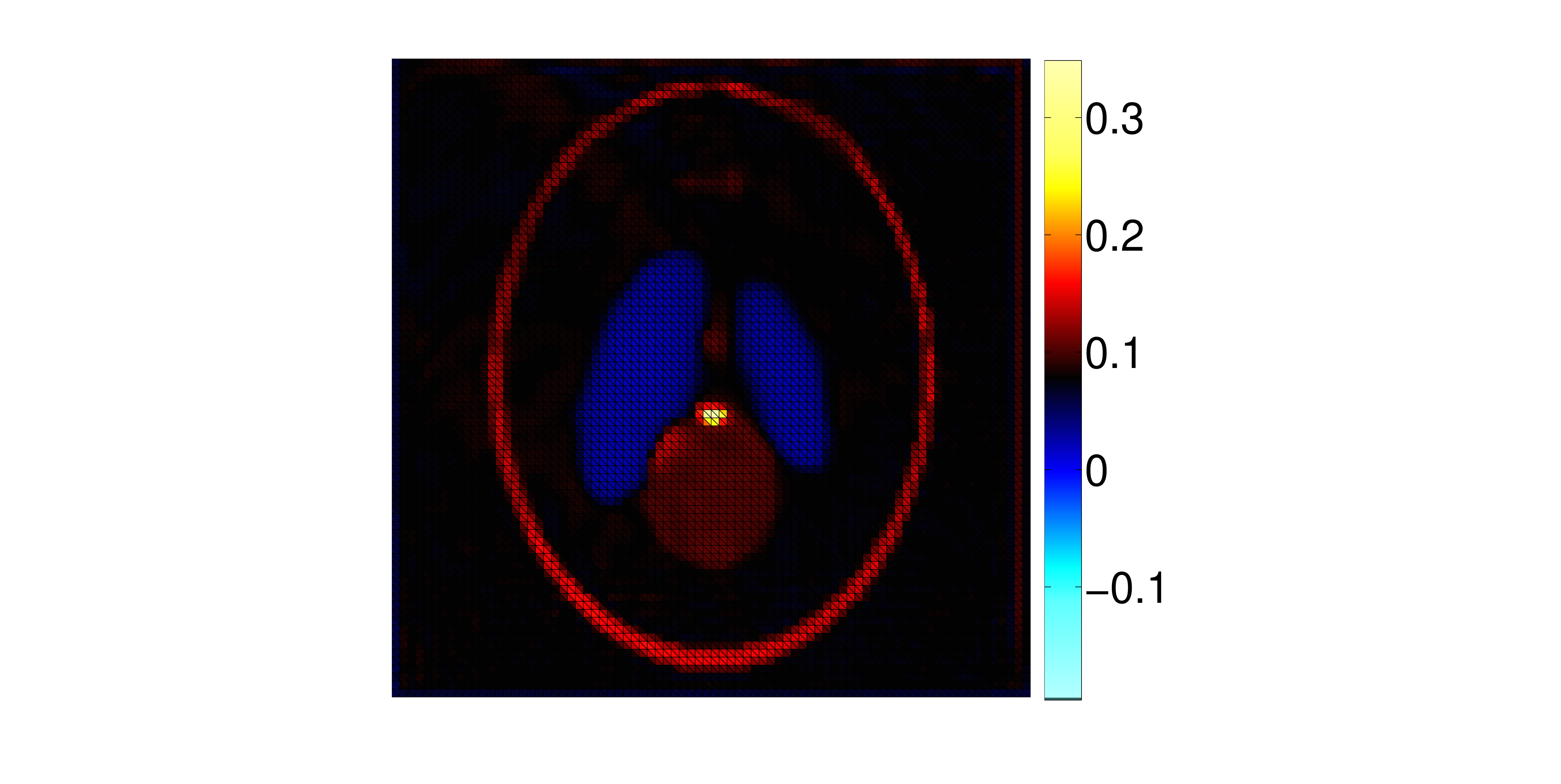}\label{F2e}}
			\subfigure[]{\includegraphics[width=.49\textwidth]{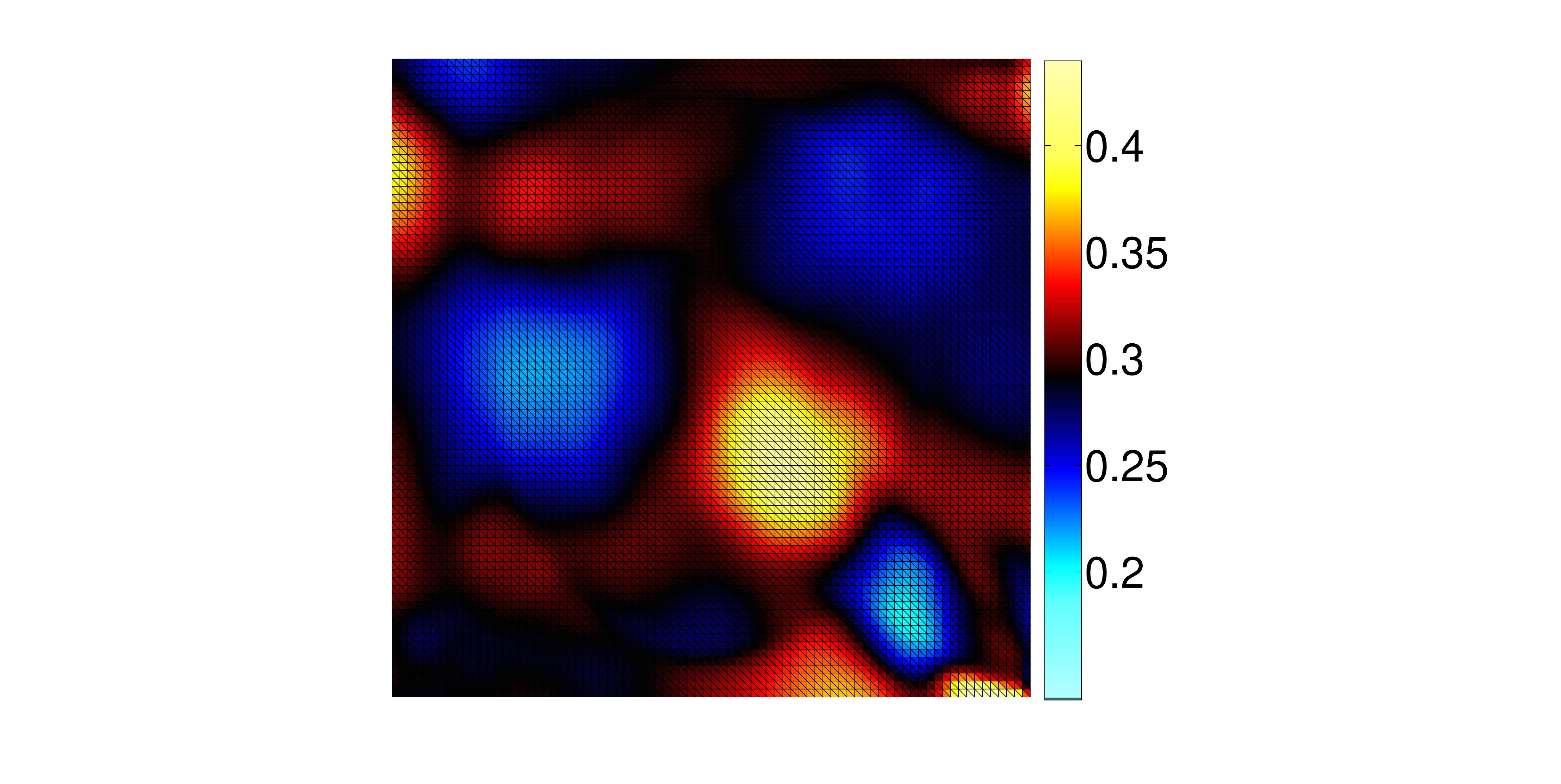}\label{F2f}}
			\subfigure[]{\includegraphics[width=.49\textwidth]{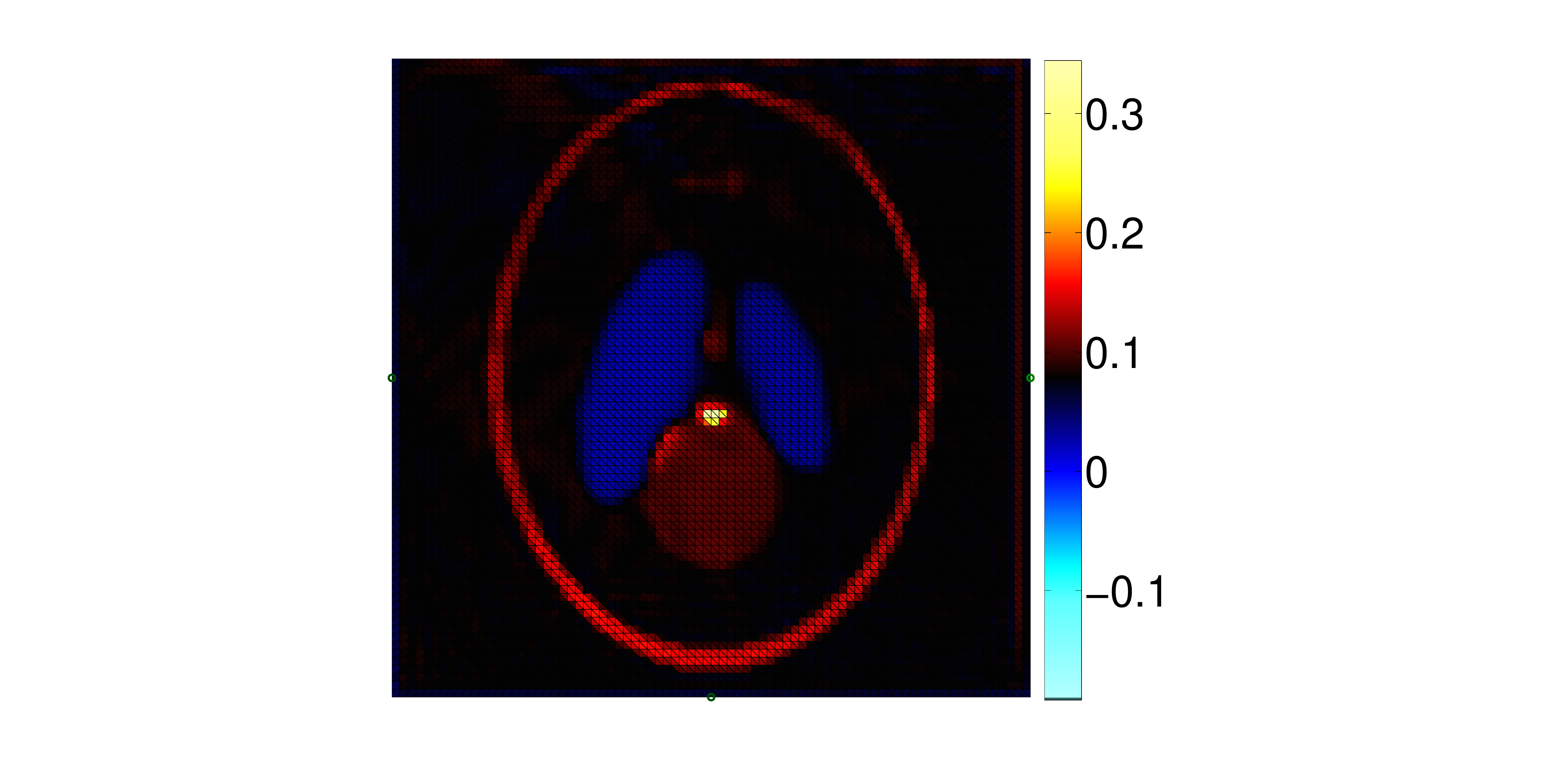}\label{F2g}}
			\subfigure[]{\includegraphics[width=.49\textwidth]{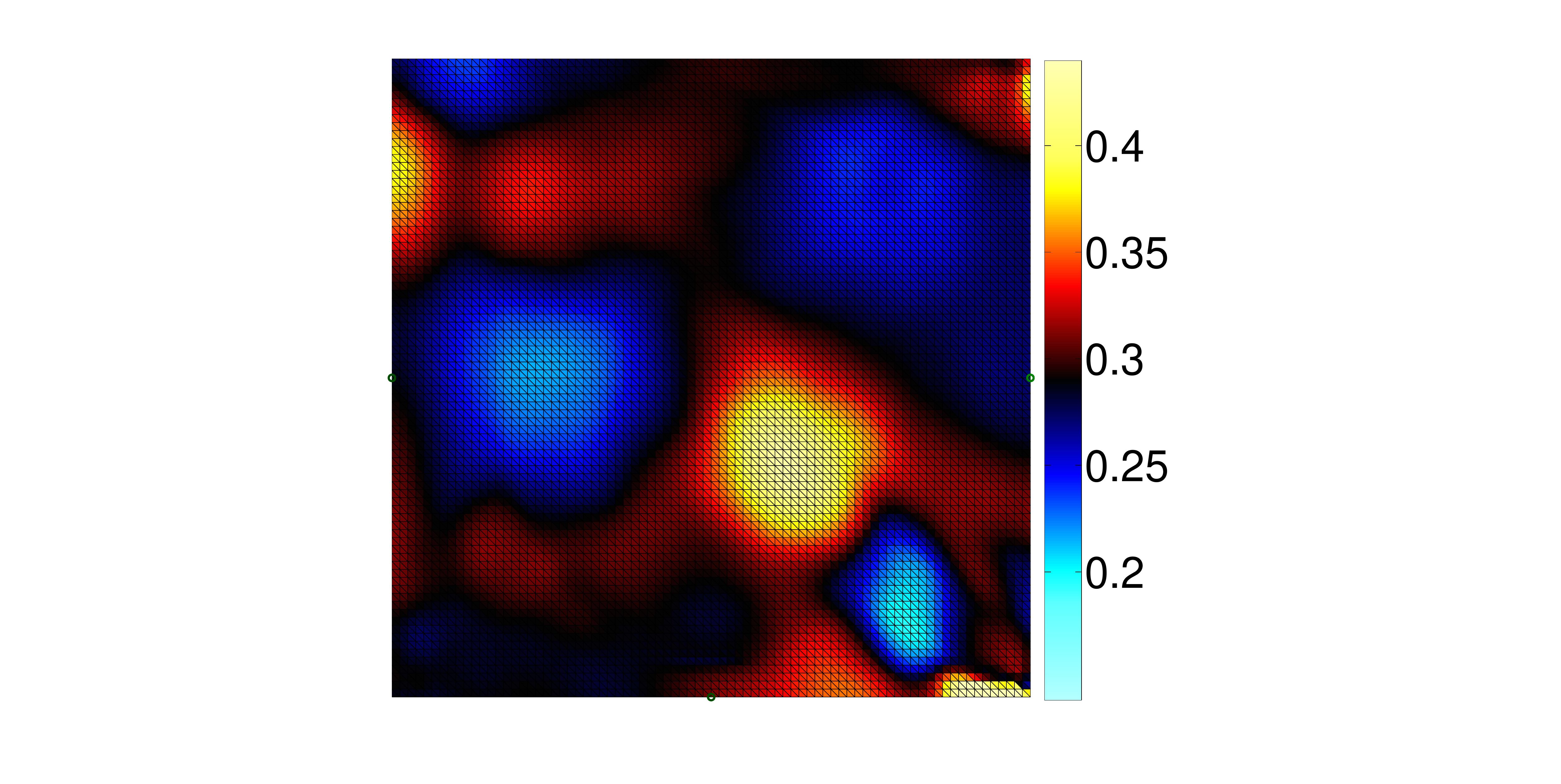}\label{F2h}}
		}
		\caption{Optical coefficients for 2D case. Phantom: (a) absorption coefficient $\mu$ (b) diffusion coefficient $\kappa$. The images reconstructed by ADMM: (c) $\mu$ (d) $\kappa$, LD: (e) $\mu$ (f) $\kappa$, and PD-IPM: (g) $\mu$ (h) $\kappa$.}
		\label{fig:42}
	\end{figure}

	\begin{figure}    \centering
		{\subfigure[]{\includegraphics[width=.49\textwidth]{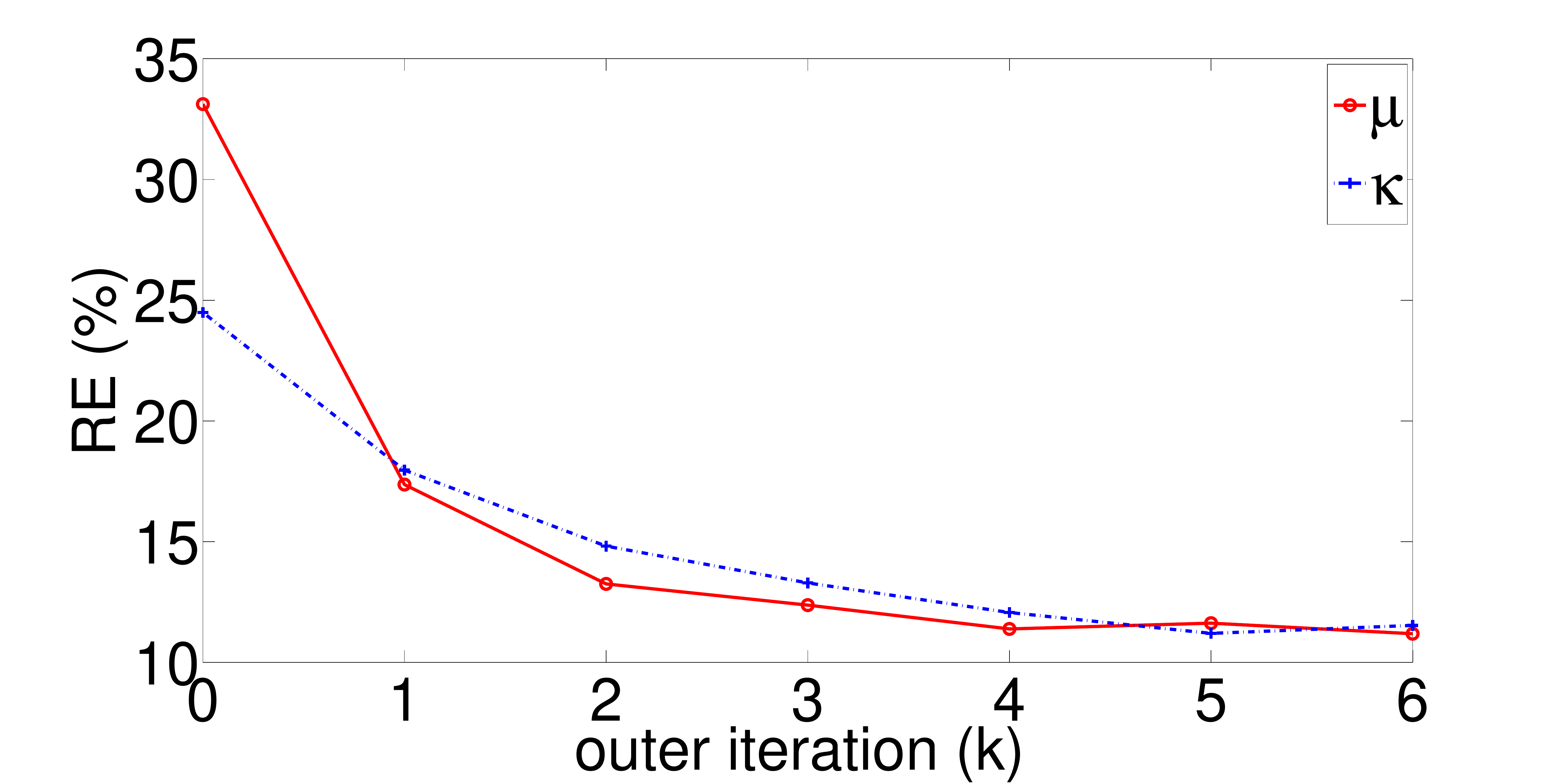}\label{F3a}}
			\subfigure[]{\includegraphics[width=.49\textwidth]{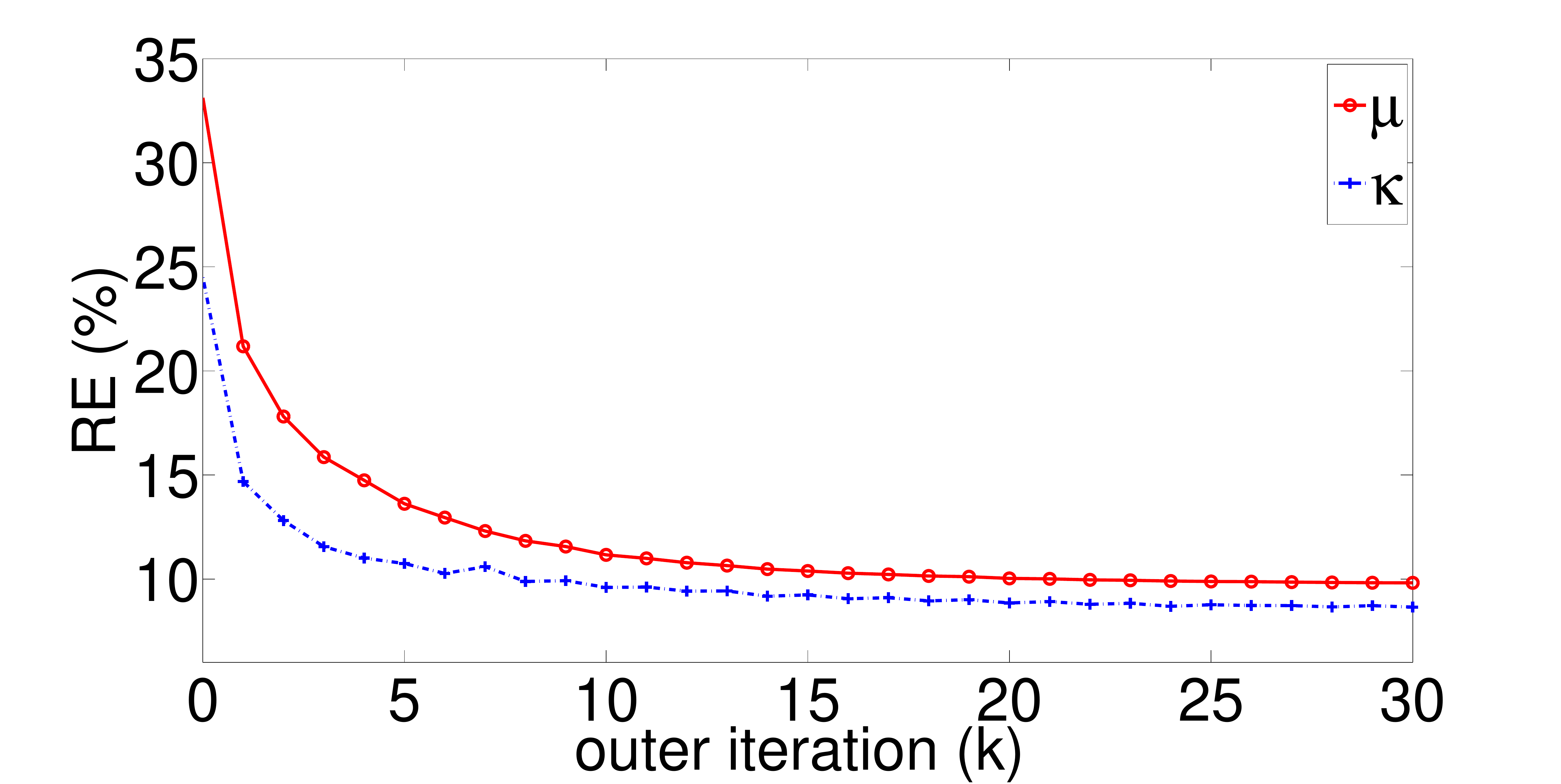}\label{F3b}}
			\subfigure[]{\includegraphics[width=.49\textwidth]{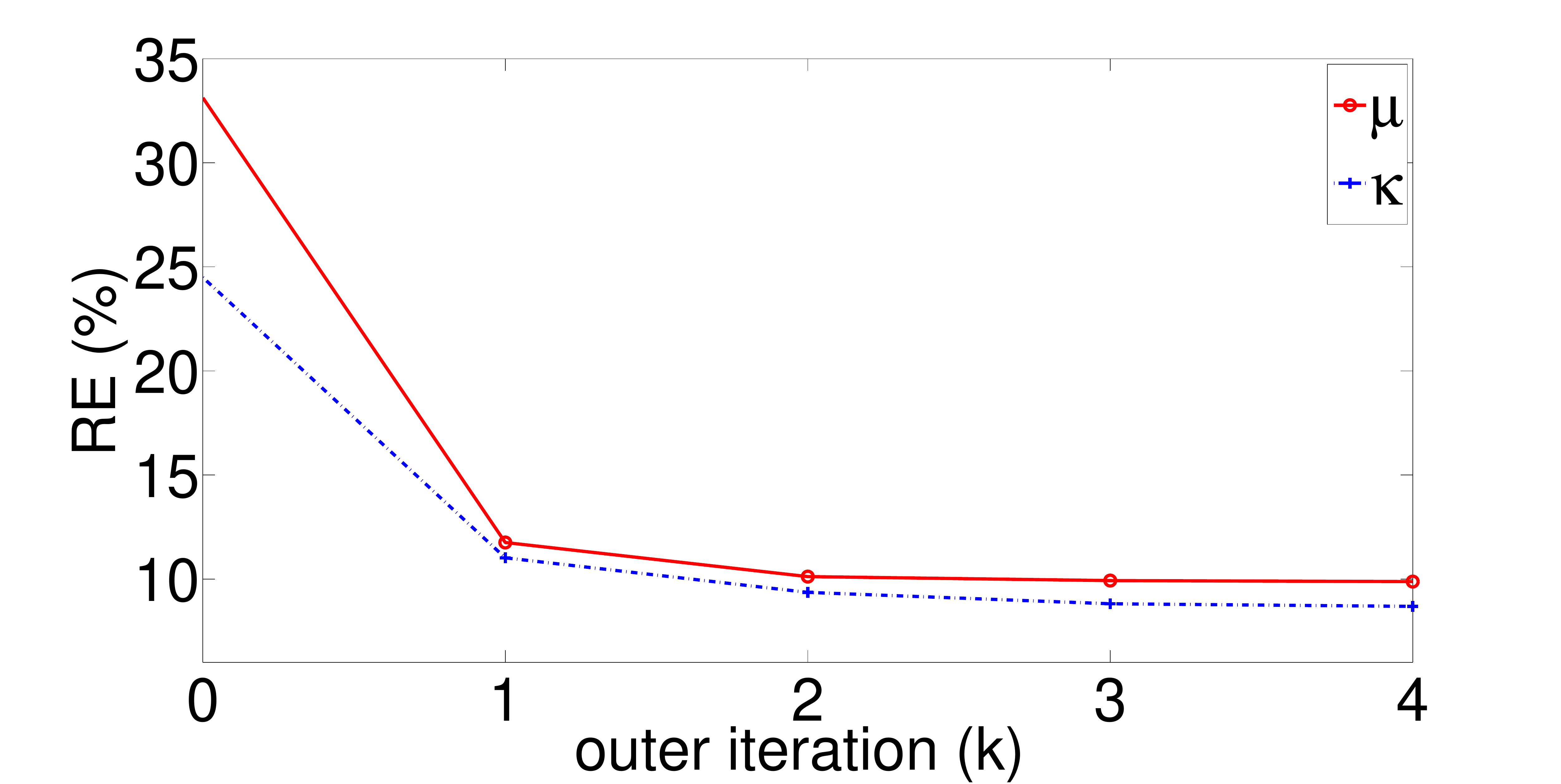}\label{F3c}}
		}
		\caption{RE versus outer iteration $k$ for 2D case: (a) ADMM (b) LD (c) PD-IPM.}
	\end{figure}

	\begin{figure}    \centering
		{\subfigure[]{\includegraphics[width=.35\textwidth]{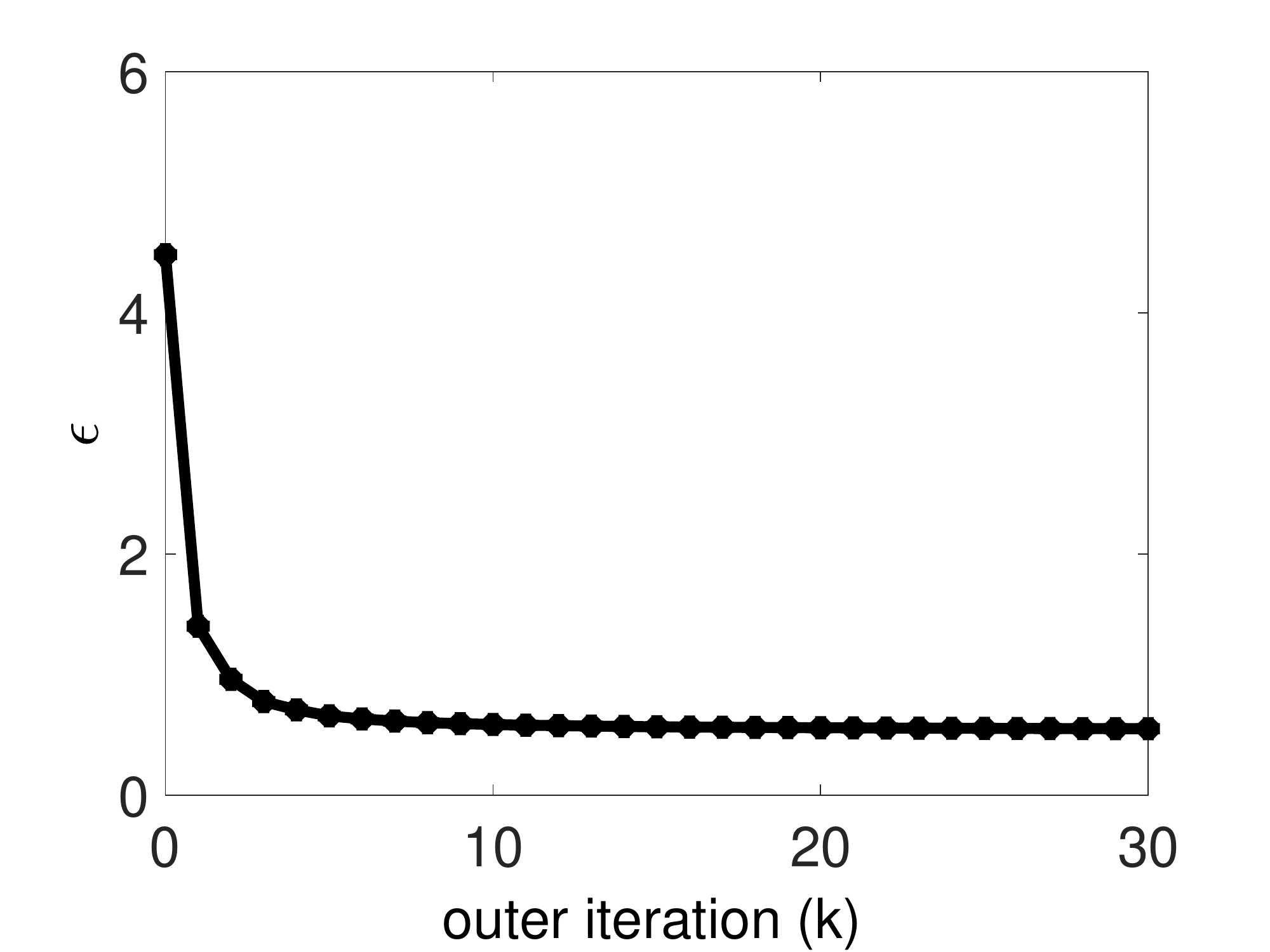}\label{F4a}}
			\subfigure[]{\includegraphics[width=.35\textwidth]{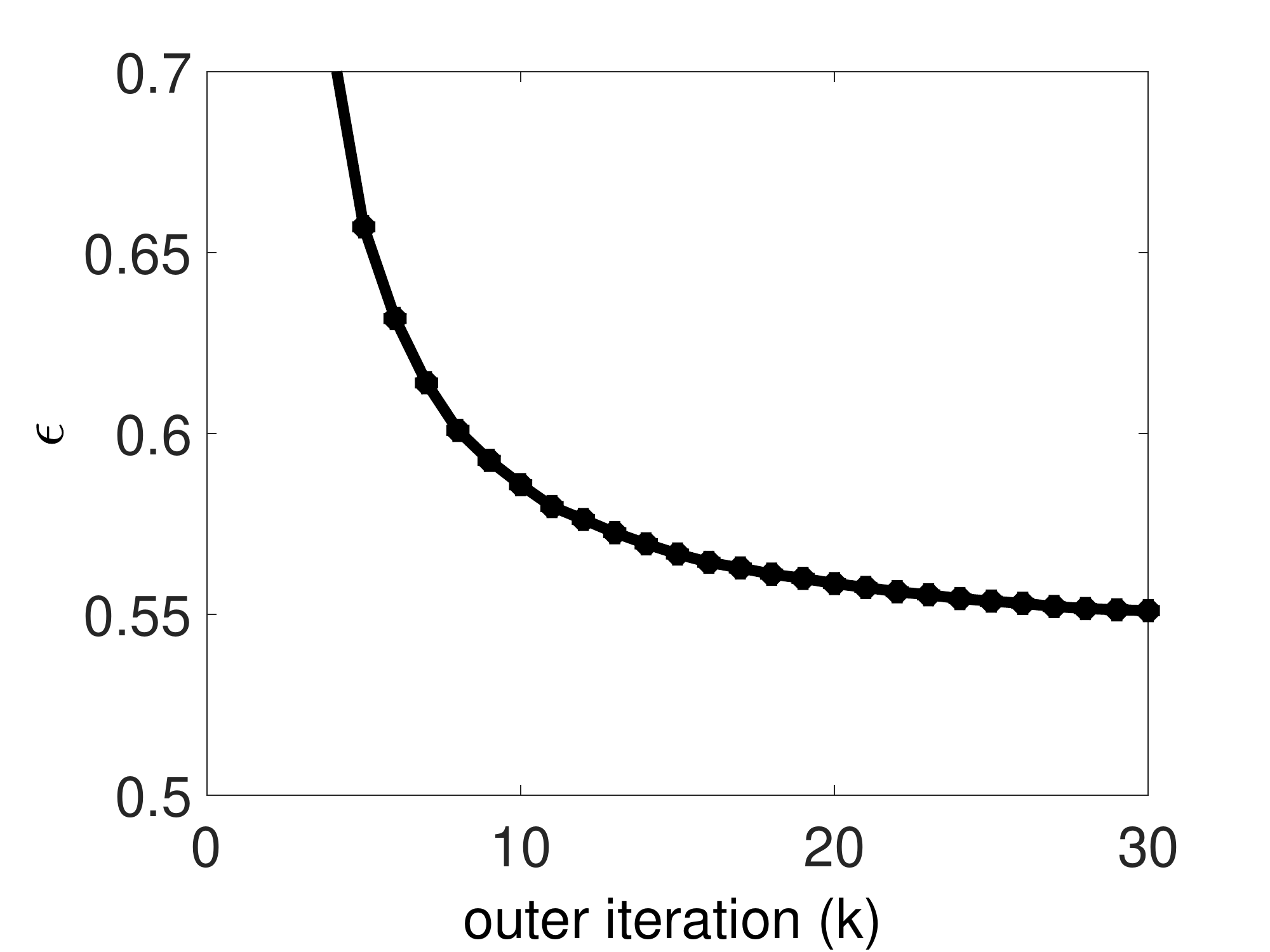}\label{F4b}}
			\subfigure[]{\includegraphics[width=.35\textwidth]{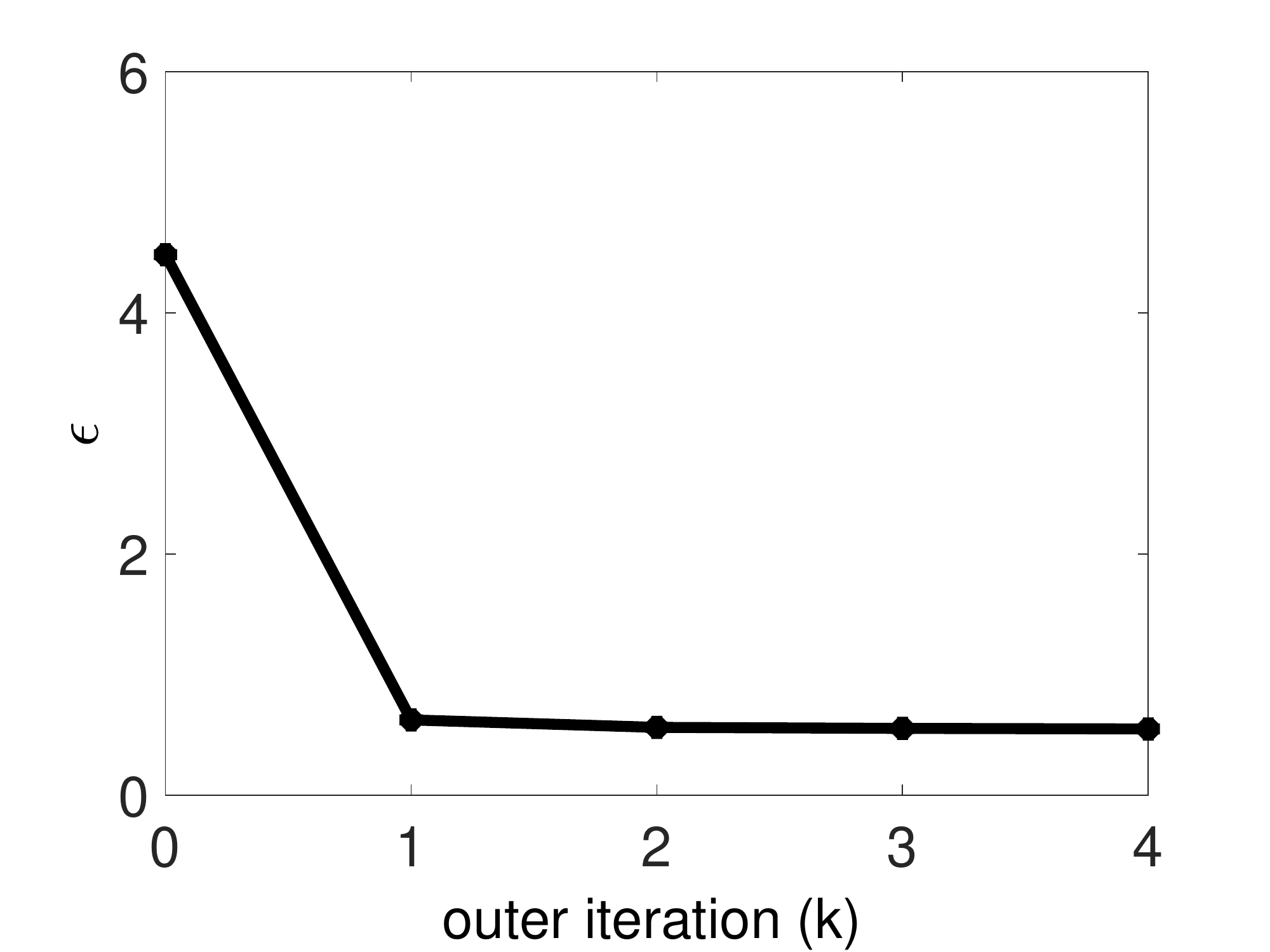}\label{F4c}}
			\subfigure[]{\includegraphics[width=.35\textwidth]{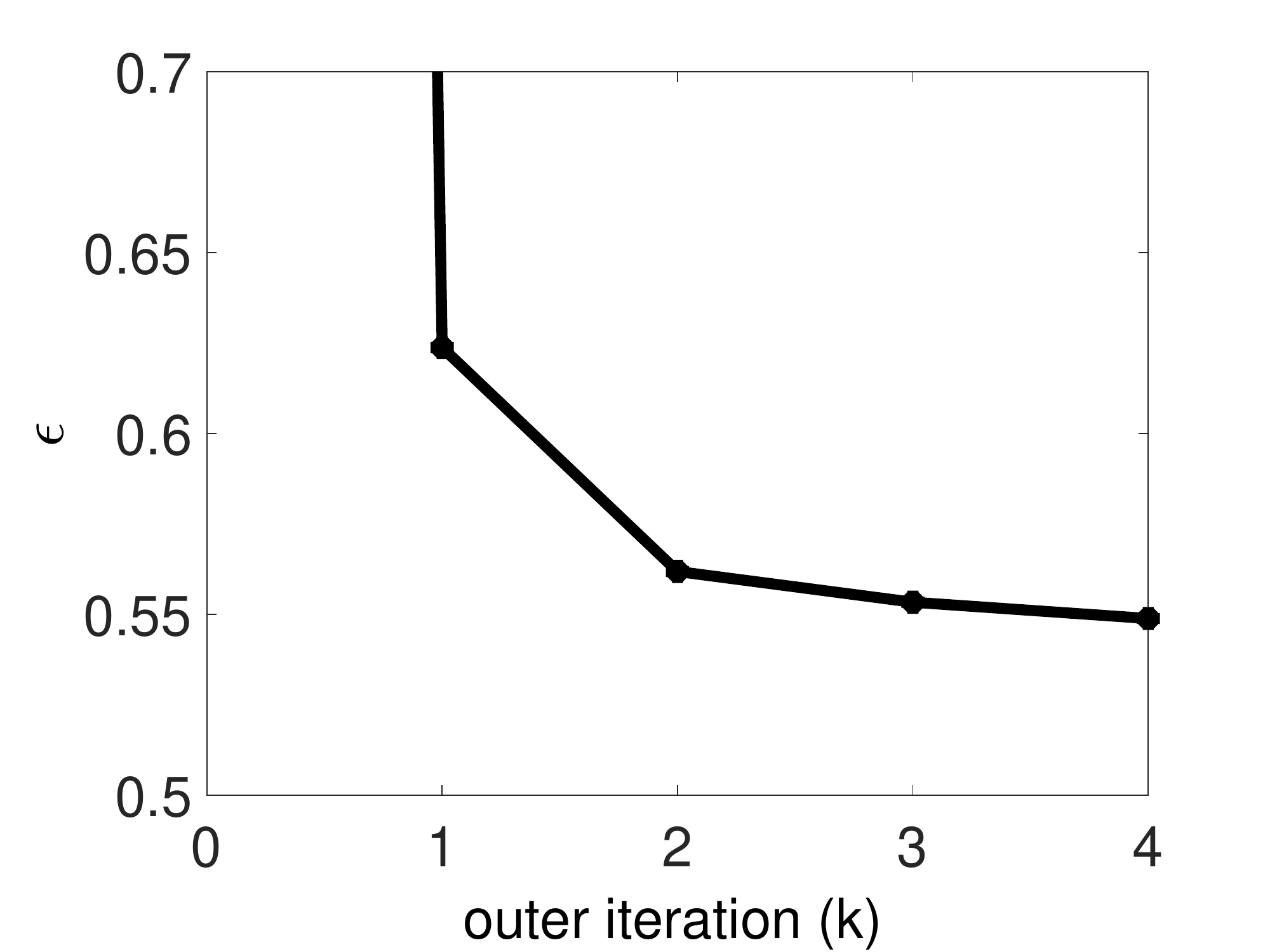}\label{F4d}}
		}
		\caption{$\epsilon$ versus outer iteration $k$ for 2D case. (a) LD  (b) LD from an enlarged view around the optimal point (c) PD-IPM (d) PD-IPM from an enlarged view around the optimal point.}
	\end{figure}
	
	Table \ref{tab2} shows the RE values for the final reconstructed images shown in figure~\ref{fig:42}.

	\begin{table}\centering
		\caption {RE(\%) of the final reconstructed images for 2D case.} \label{tab2} 
		\begin{tabular}{cclcl}
			Methods & \multicolumn{2}{c}{\textbf{$\mu$}} & \multicolumn{2}{c}{$\kappa$} \\ \hline
			ADMM    & \multicolumn{2}{c}{11.1882}        & \multicolumn{2}{c}{11.5316}  \\ 
			LD      & \multicolumn{2}{c}{9.8799}         & \multicolumn{2}{c}{8.6899}   \\ 
			PD-IPM & \multicolumn{2}{c}{9.8513}          & \multicolumn{2}{c}{8.6547}   \\ \hline
		\end{tabular}
	\end{table}
	
	\subsection{3D phantom}
	We performed our 3D simulation on a cubic domain $[-5,5] \times [-5,5]   \times  [-5, 5] \ \text{mm}^3 $.

	\subsubsection{Optical excitation}
	We used three optical excitation patterns, i.e., $N_q=3$. For each optical excitation $q$, we used a discretisation of an inward directed diffuse boundary current $I_{s,q} \ (J/\text{mm}^2)$ that obeys \eqref{isn} with $\iota_q \in \partial \Omega $ two confronting faces of the grid, i.e. the left-right, posterior-anterior, and bottom-top faces.  
	
	\subsubsection{Discretisation for data generation}
	For data generation, the cubic domain was discretised using a grid with $37 \times 37 \times 37$ nodes and an even separation distance of $2.78 \times 10^{-1}$mm along all Cartesian coordinates. For the optical portion of the problem, an FE mesh was simulated so that each set of six tetrahedral voxels forms a cubic pixel with a centre matching an associated node on the acoustic grid. For the acoustic portion of the problem, we added a perfectly matched layer (PML) with a thickness of $8$ grid points and a maximum attenuation coefficient of $2$ nepers per grid point. The acoustic wavefield was detected in $292$ time instants using $2145$ detectors that are equidistantly placed on two (the left and posterior) faces of the grid. A $30 \ \text{dB}$ Additive White Gaussian Noise (AWGN) was then added to the simulated data.
	
	\subsubsection{Discretisation for image reconstruction}
	For image reconstruction, we avoided an inverse crime for discretisation by using a grid made up of $33 \times 33 \times 33$ nodes with a homogeneous separation distance of $3.125 \times 10^{-1}$mm along all Cartesian coordinates.
	
	\subsubsection{Acoustic properties}
	As discussed in section \ref{acous2d}, to avoid an inverse crime for acoustic properties of the medium, we corrupted the sound speed and ambient density with $30 \ \text{dB}$ AWGN noise for simulation of data, whereas we used the clean acoustic maps for image reconstruction. Table \ref{tab3} shows the minimial and maximal values for these maps. Using this table, the grid for data generation (resp. image reconstruction) supports a maximal frequency of $2.0766$ MHz (resp. $2.156$ MHz). The distributions of the sound speed and ambient density for the 3D phantom for data generation (resp. image reconstruction) are shown from a top view in figures \ref{F5a} and \ref{F5b} (resp. \ref{F5c} and \ref{F5d}), respectively. Additionally, the acoustic attenuation coefficient and the associated exponent factor (cf. equation \eqref{Twelve}) was simulated the same as the 2D phantom.

	\begin{table} \centering
		\caption {The minimal and maximal values for acoustic properties of the 3D phantom.} \label{tab3} 
		\begin{tabular}{ccccc}
			\hline
			& \multicolumn{2}{c}{\text{$c_0 (\text{m} \text{s}^{-1}$) }}            & \multicolumn{2}{c}{$\rho_0 (\text{kg} \text{m}^{-3}$)}                 \\ \cline{2-5} 
			& min                   & max                    & min                   & max                   \\ \hline
			Data generation      & $1.186 \times 10^3$ & $2.077 \times 10^3$ & $0.853 \times 10^3 $ & $1.374 \times 10^3$ \\ 
			Image reconstruction  & $1.390 \times 10^3$ & $1.897 \times 10^3$ & $0.956 \times 10^3 $ & $1.252 \times 10^3$                                           \\ \hline
		\end{tabular}
	\end{table}
	
	\subsubsection{Optical phantom}
	The left columns in figures \ref{F6a} and \ref{F6b} show the distributions of $\mu$ and $\kappa$ for the 3D phantom, respectively. The images are obtained in horizontal planes (slices)
	$z=\left\{3,2,1,0,-1,-2,-3\right\}$ mm, and the colorbars are shown to the right of images.

	\begin{figure}    \centering
		{\subfigure[]{\includegraphics[width=.49\textwidth]{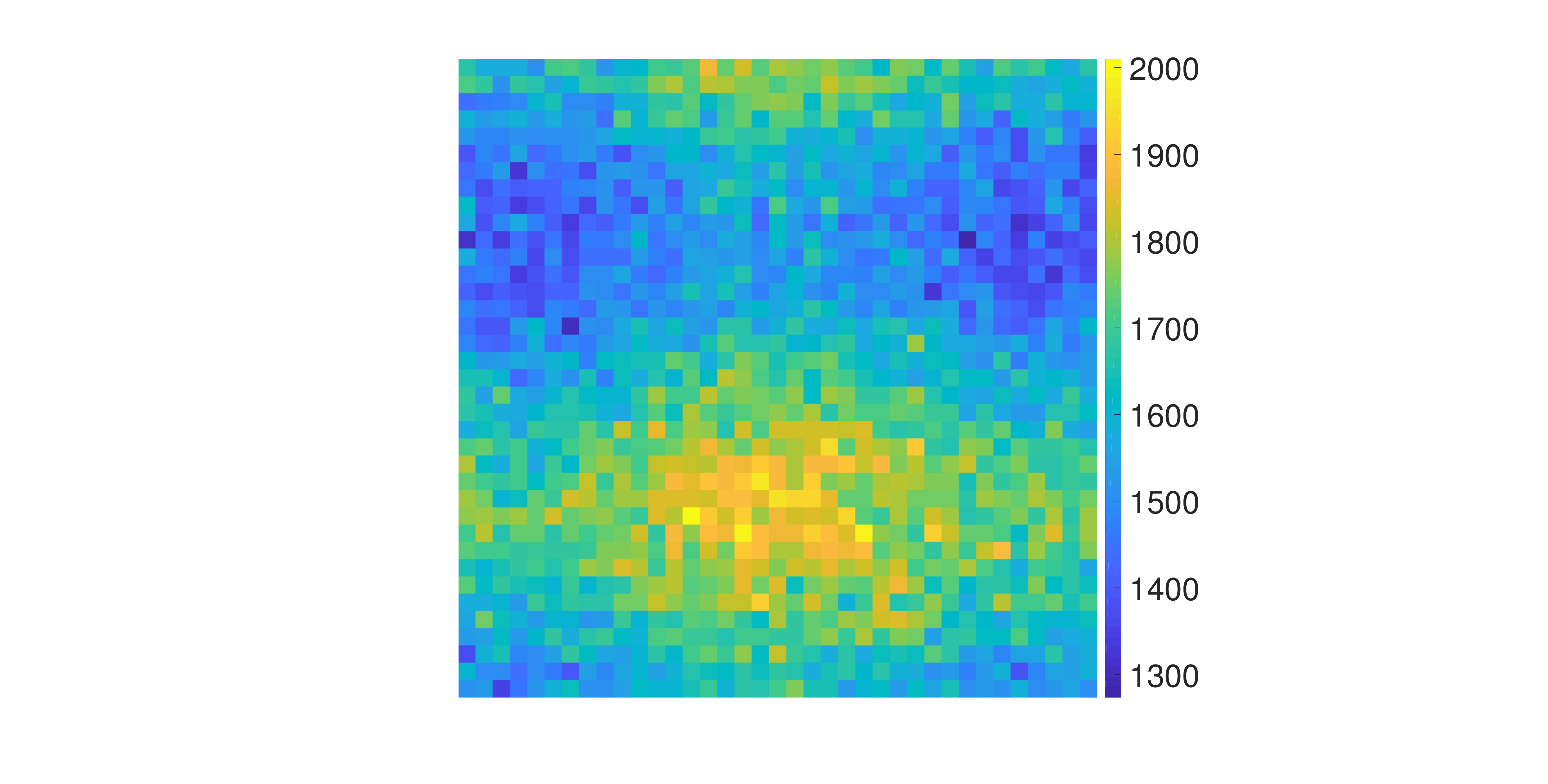}\label{F5a}}
			\subfigure[]{\includegraphics[width=.49\textwidth]{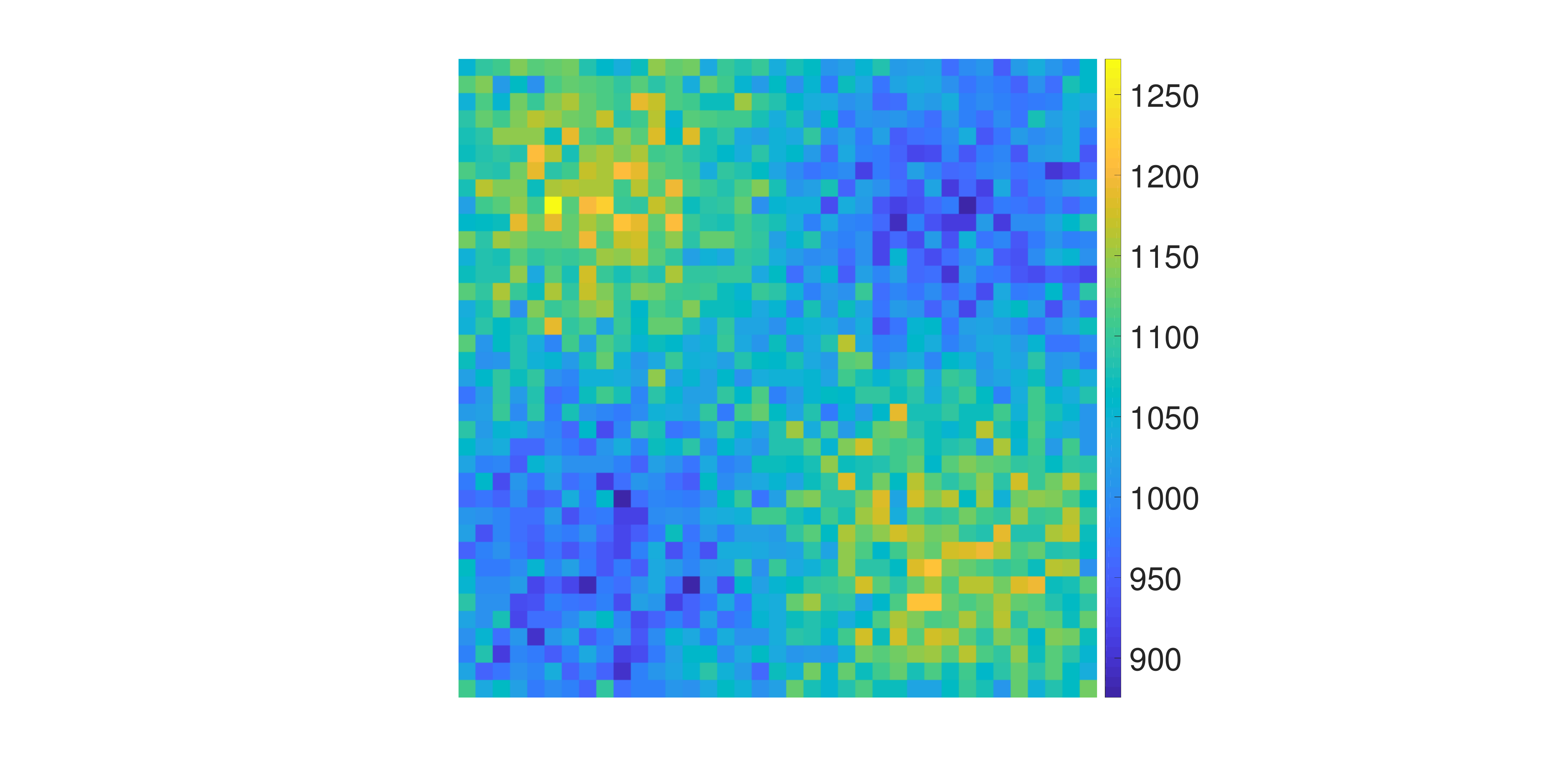}\label{F5b}}
			\subfigure[]{\includegraphics[width=.49\textwidth]{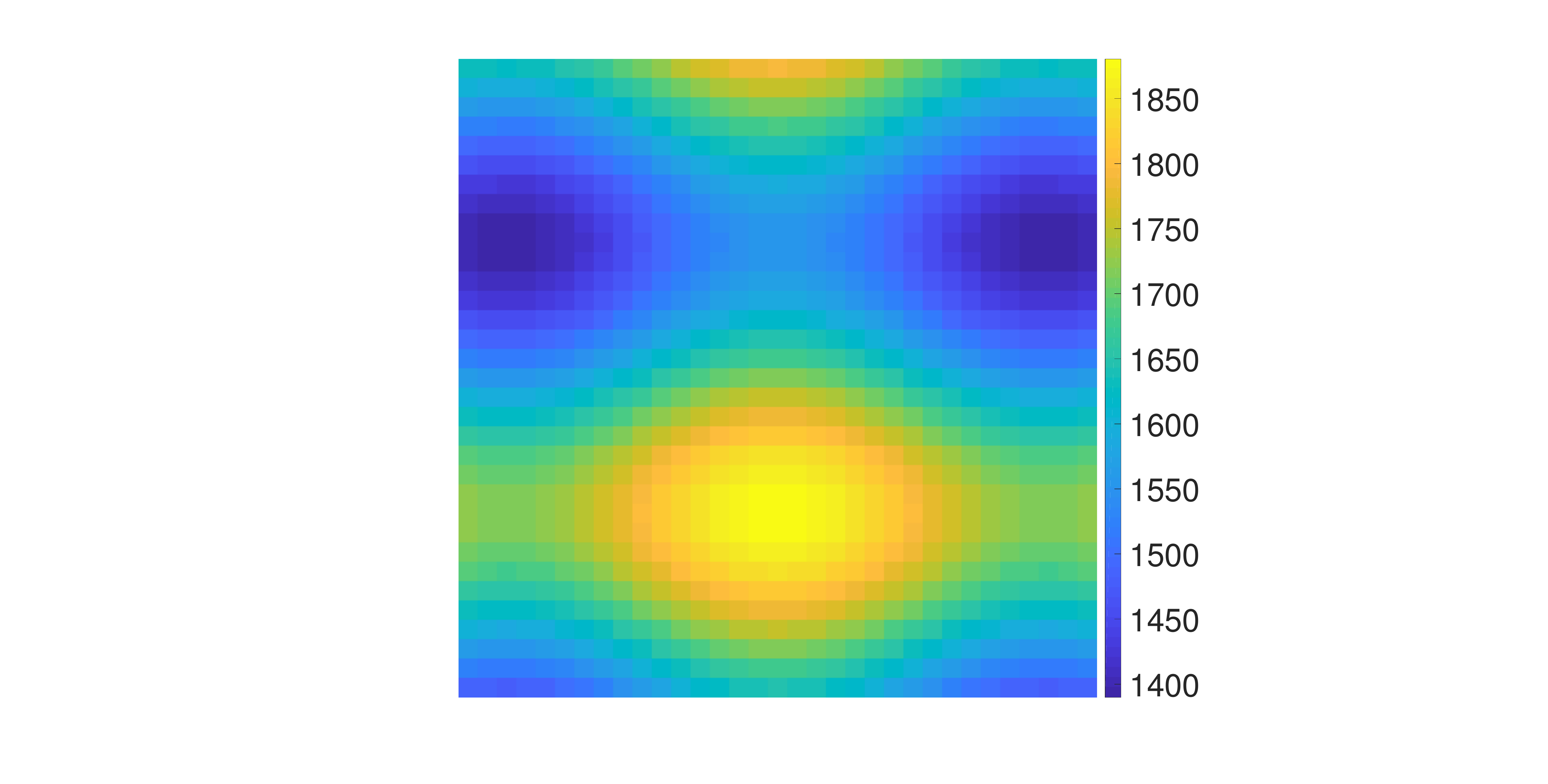}\label{F5c}}
			\subfigure[]{\includegraphics[width=.49\textwidth]{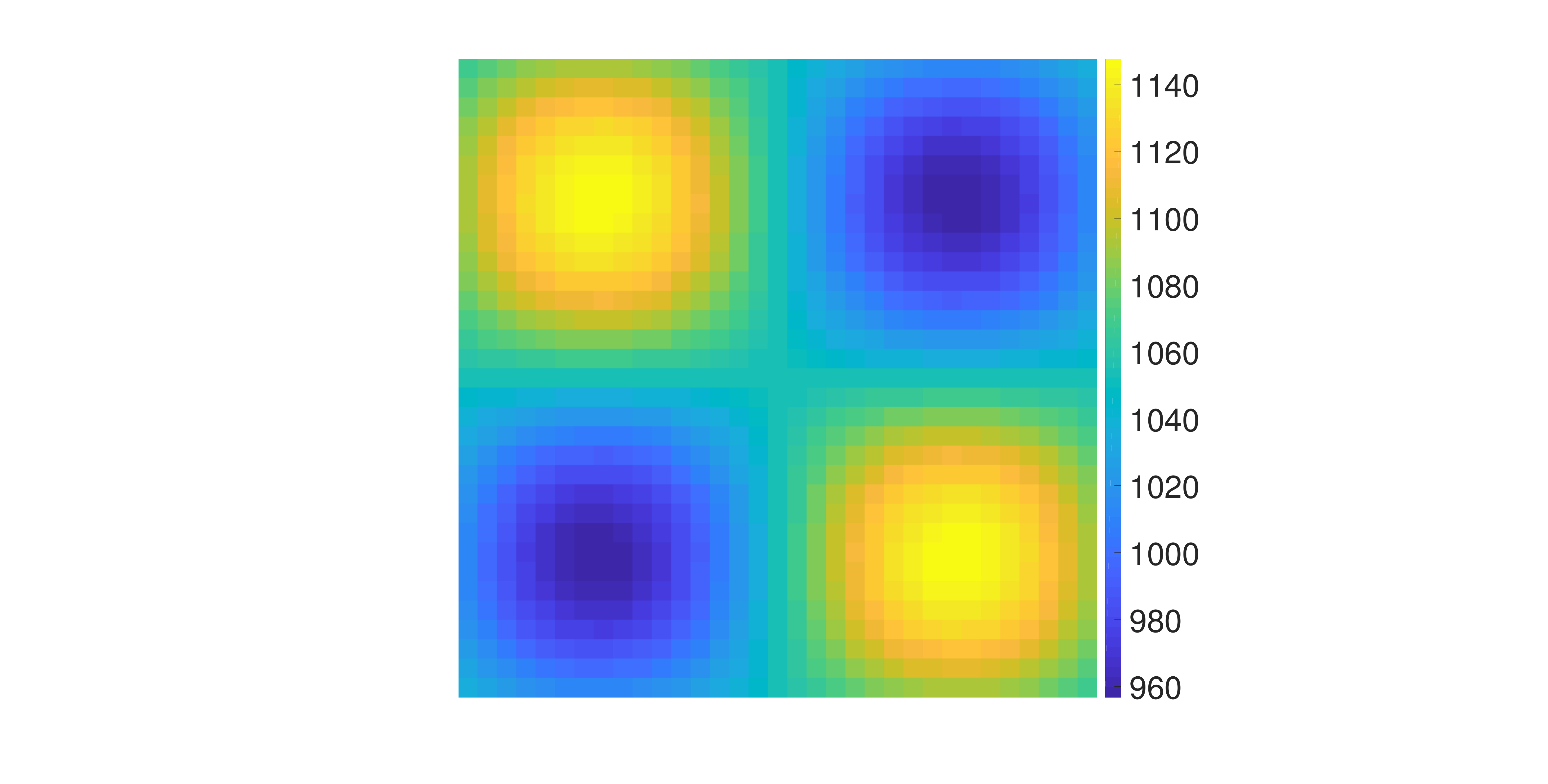}\label{F5d}}
		}
		\caption{Acoustic properties for 3D case from a top view. Data generation: (a) $c_0$ (b) $\rho_0$, and image reconstruction: (c) $c_0$ (d) $\rho_0$.}
		\label{fig:45}
	\end{figure}

	\subsubsection{Image reconstruction}  \label{imr3}
	All algorithms were initialised using values $1.2$ times more than the mean of optical coefficients for the 3D phantoms.
	
	\textit{ADMM.} The parameters for an implementation of the ADMM algorithm were chosen carefully in order to obtain almost the best possible image. All these parameres match our choices for the 2D phantom.
	
	\textit{LD.} For an implementation of the LD algorithm, the arising linearised subproblems were solved using a PCG algorithm with the same parameters as the 2D phantom, except that for stopping each PCG loop, we used $i_m=3$. The TV preconditioner $M$ was applied using $\gamma=1 \times 10^{-8}$ and $\beta= 1 \times 10^{-6}$. We stopped our LD algorithm using $Tol_\text{out}=1 \times 10^{-3}$.
	
	\textit{PD-IPM.} For an implementation of the PD-IPM algorithm, each outer linearised problem was solved using a sequence of inner linearised subproblems. The TV preconditioner was applied using $\gamma=1 \times 10^{-8}$ and $\beta= 1 \times 10^{-8}$. For termination of each PCG loop associated with each inner linearised problem, we used $i_m=2$ and $i_\text{max}=30$ (cf. Algorithm~\ref{alg:44}). We also terminated each outer linearised subproblem using $k'_\text{max}=25$ and $Tol_{\text{med}}= 1 \times 10^{-3}$. Also, our PD-IPM algorithm was terminated using a stopping threshold $Tol_{out}= 1 \times 10^{-3}$.

	\subsubsection{Observations}
	
	\textit{ADMM.} The stopping criterion for the ADMM algorithm was satisfied after outer iteration 7. In figures \ref{F6a} and \ref{F6b}, the second columns (from the left side) show the final reconstructed images for $\mu$ and $\kappa$, respectively.
	These images are shown using horizonal slices the same as the first columns for the phantom. Figure \ref{F7a} shows the RE of solutions computed by the ADMM algorithm for $\mu$ and $\kappa$ at outer iterations.

	\textit{LD.} The LD algorithm was terminated after 10 iterations. The final reconstructed images for $\mu$ and $\kappa$ are shown in the 3rd columns of figures \ref{F6a} and \ref{F6b}, respectively. The images are shown in the same way as the ADMM method. Figure \ref{F7b} shows the RE of the solutions computed by LD for outer iterations $k$. The computed values for $\epsilon(X^k)$ are shown in figure \ref{F8a}, and \ref{F8b} from an enlarged view around the optimal point. As shown in these figures, $\epsilon$ monotonically converges to a minimiser for all outer iterations.

	\textit{PD-IPM.} The stopping criterion for the PD-IPM algorithm was satisfied after 6 outer iterations. The final reconstructed images for $\mu$ and $\kappa$ are shown in the 4th columns of figures \ref{F6a} and \ref{F6b}, respectively. Figure \ref{F7c} shows the RE of solutions (optical coefficients) computed by the PD-IPM algorithm for outer iterations $k$. Figures \ref{F8c} and \ref{F8d} show the obtained values for $\epsilon(X^k)$. As shown in these figures, our choices for the step sizes associated with outer (resp. inner) subproblems, i.e., $\alpha^k=1$ (resp. $\alpha_{k'}=1$), provided a monotonic reduction for $\epsilon^k$ (resp. $\epsilon_{k'}$).

	Table \ref{tab4} shows the RE values for the final reconstructed images shown in figure~\ref{fig:46}.

	\begin{figure} \centering
		{\subfigure[]{\includegraphics[width=1.0\textwidth]{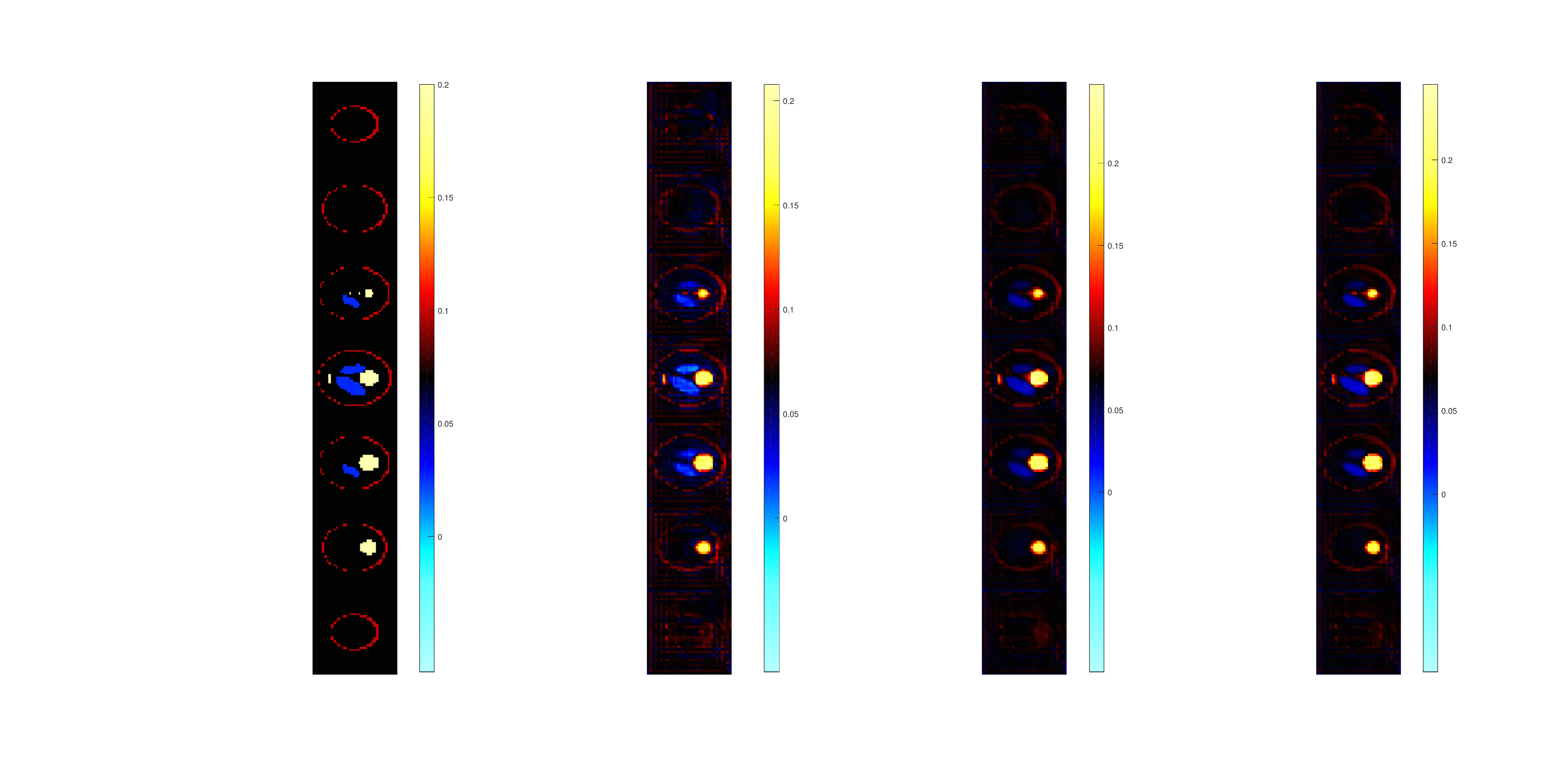}\label{F6a}}
			\subfigure[]{\includegraphics[width=1.0\textwidth]{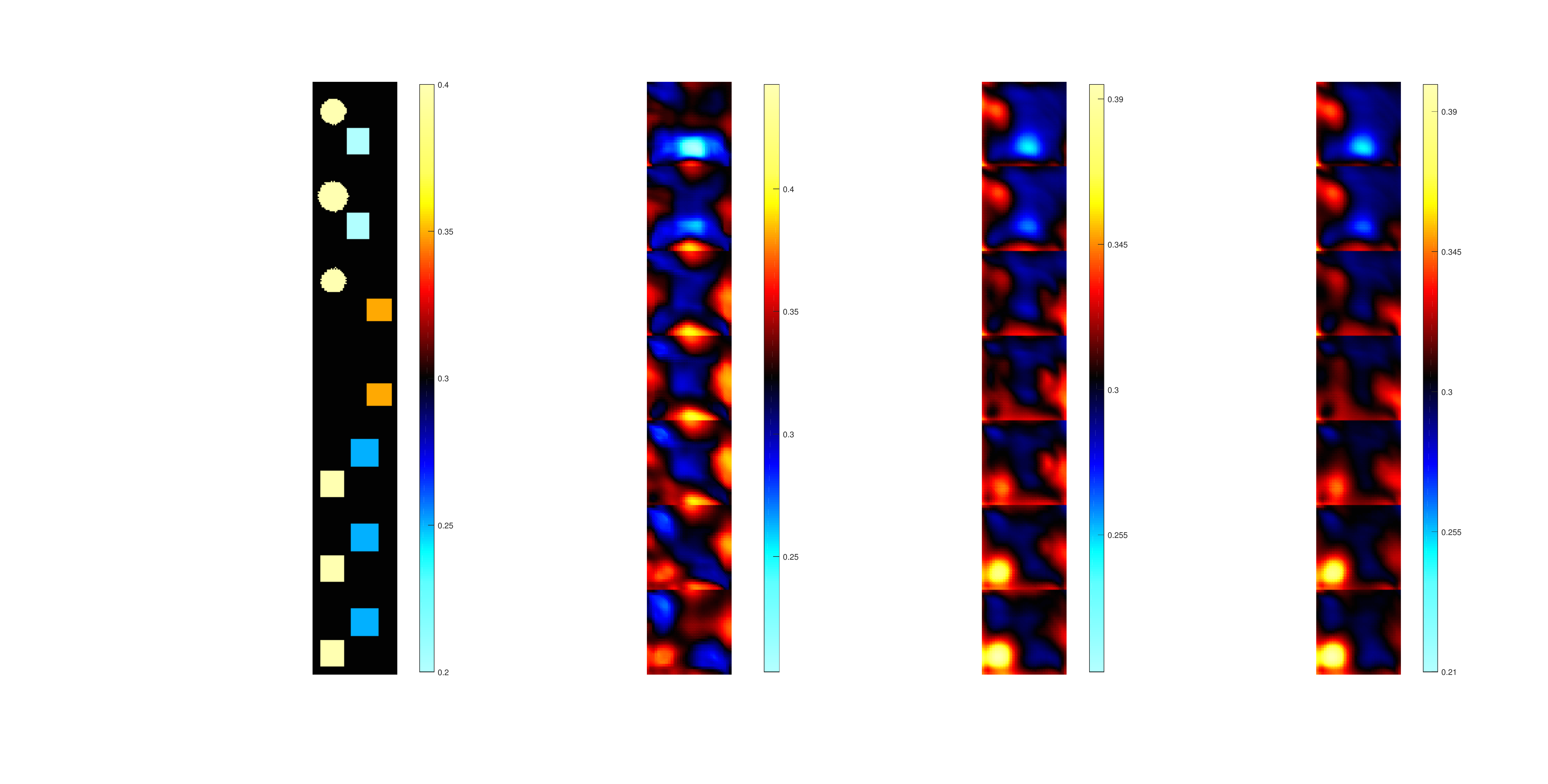}\label{F6b}}
		}
		\caption{Optical coefficients for 3D case. (a)  $\mu$, from the left to right: phantom, ADMM, LD and PD-IPM. (b)  $\kappa$, from the left to right: phantom, ADMM, LD and PD-IPM.
		}
		\label{fig:46}
	\end{figure}

	\begin{figure}  \centering
		{\subfigure[]{\includegraphics[width=.49\textwidth]{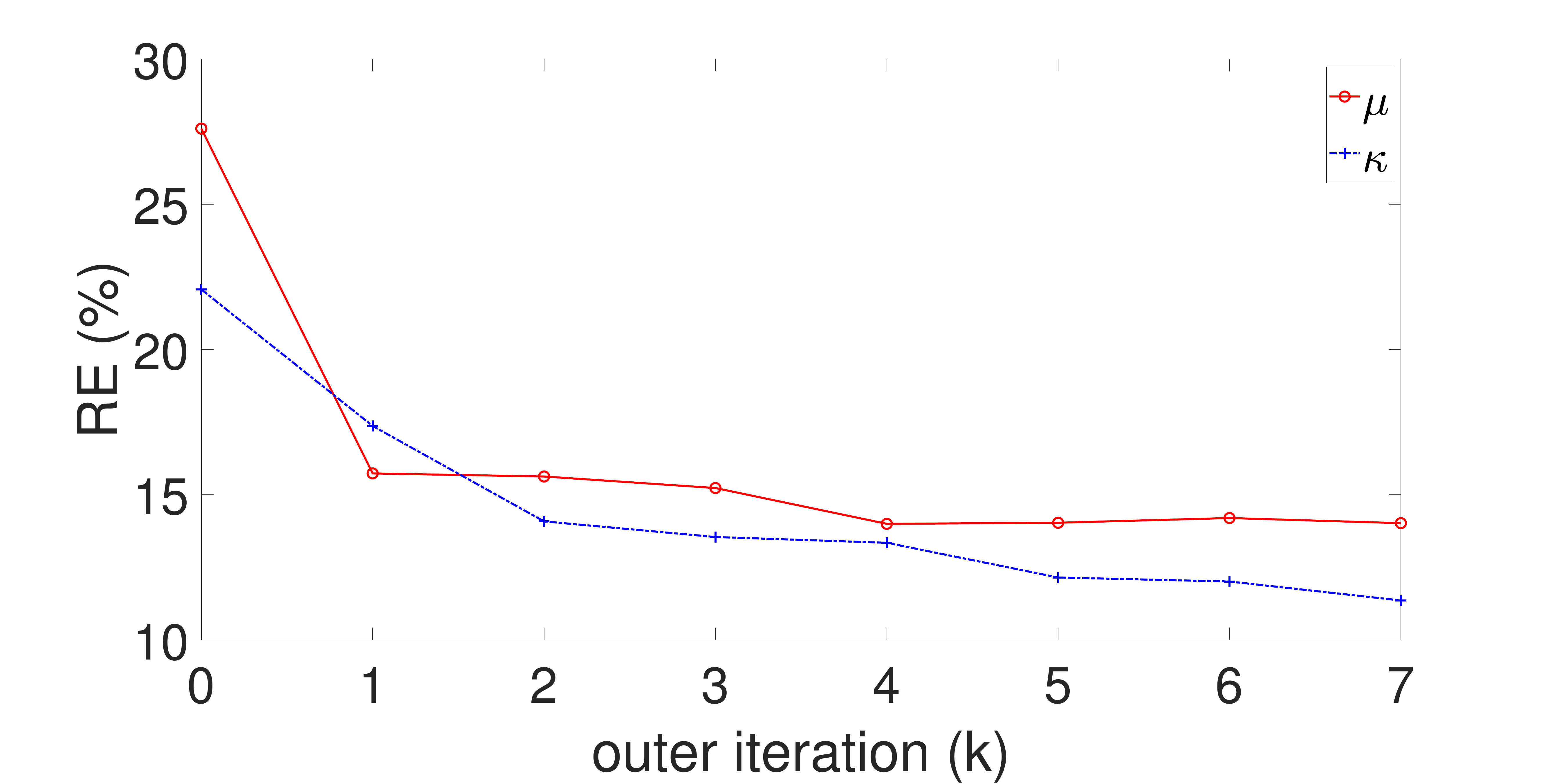}\label{F7a}}
			\subfigure[]{\includegraphics[width=.49\textwidth]{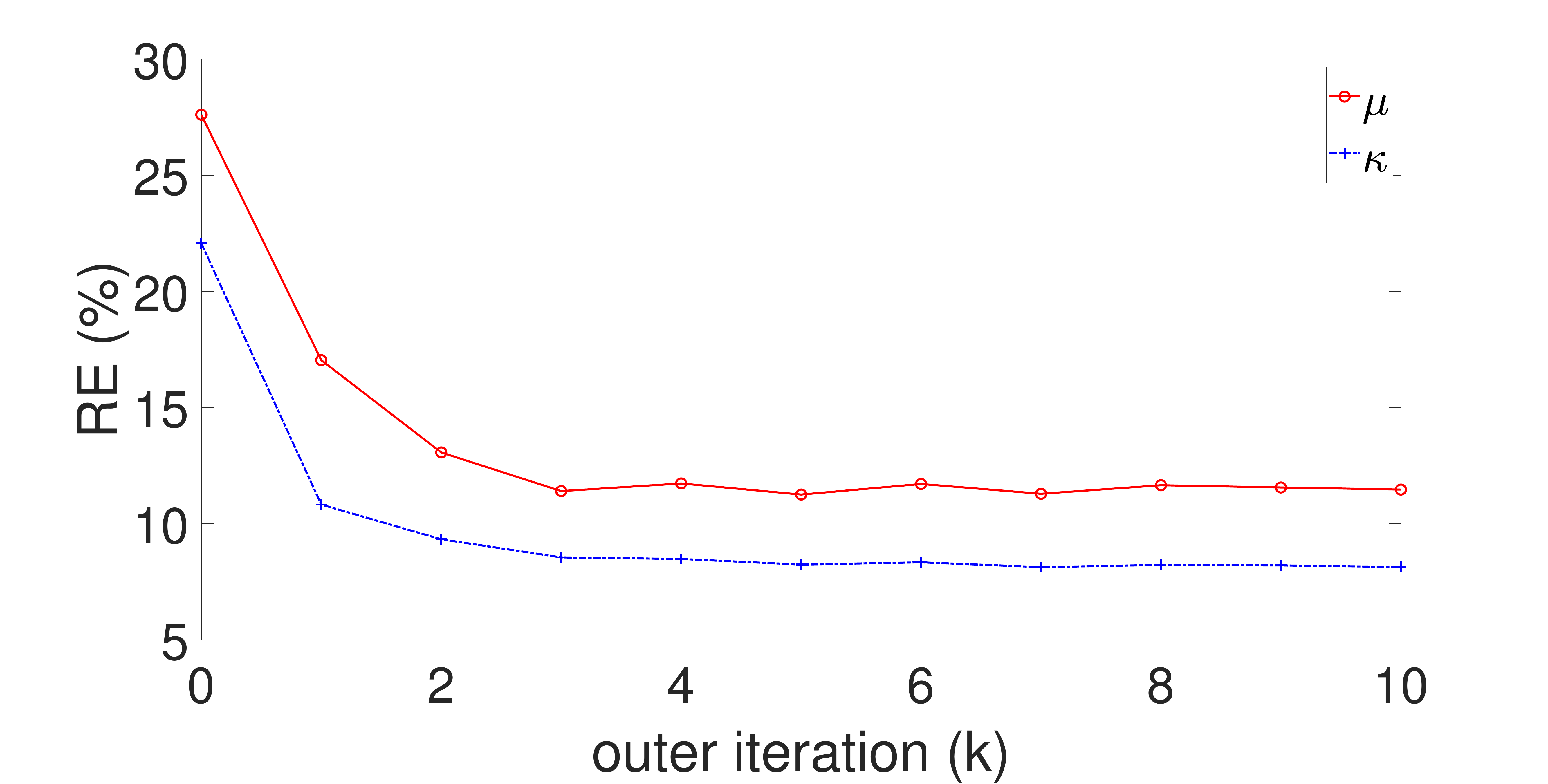}\label{F7b}}
			\subfigure[]{\includegraphics[width=.49\textwidth]{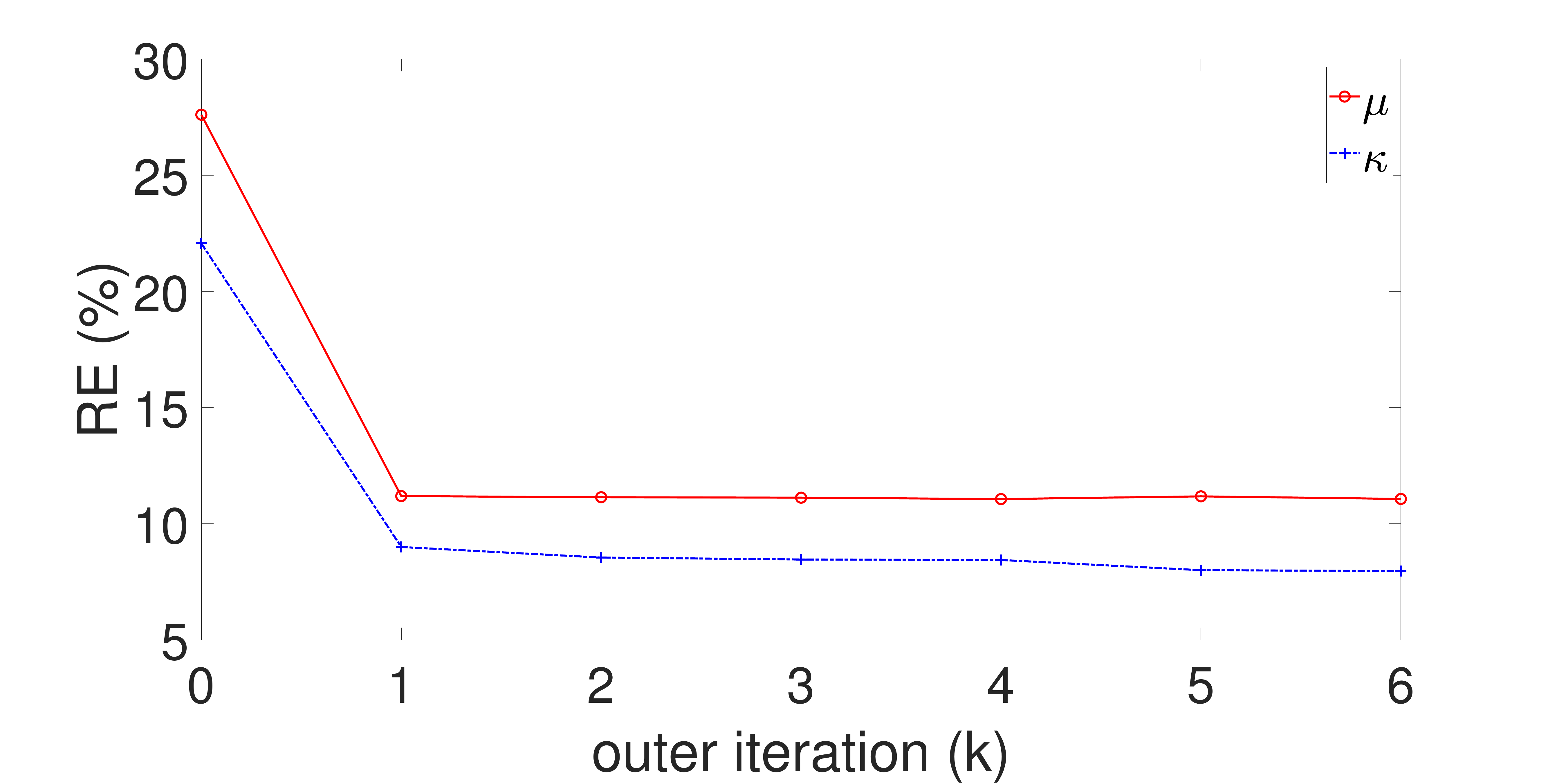}\label{F7c}}
		}
		\caption{RE versus outer iteration $k$ for 3D case. (a) ADMM (b) LD (c) PD-IPM.}
	\end{figure}

	\begin{table}  \centering
		\caption {RE(\%) of the final reconstructed images for 3D case.} \label{tab4} 
		\begin{tabular}{cclcl}
			Methods & \multicolumn{2}{c}{\textbf{$\mu$}} & \multicolumn{2}{c}{$\kappa$} \\ \hline
			ADMM    & \multicolumn{2}{c}{14.0210}        & \multicolumn{2}{c}{11.3598}  \\ 
			LD      & \multicolumn{2}{c}{11.4718}         & \multicolumn{2}{c}{8.1426}   \\ 
			PD-IPM & \multicolumn{2}{c}{11.0694}          & \multicolumn{2}{c}{7.9647}   \\ \hline
		\end{tabular}
	\end{table}

	\begin{figure}  \centering
		{\subfigure[]{\includegraphics[width=.35\textwidth]{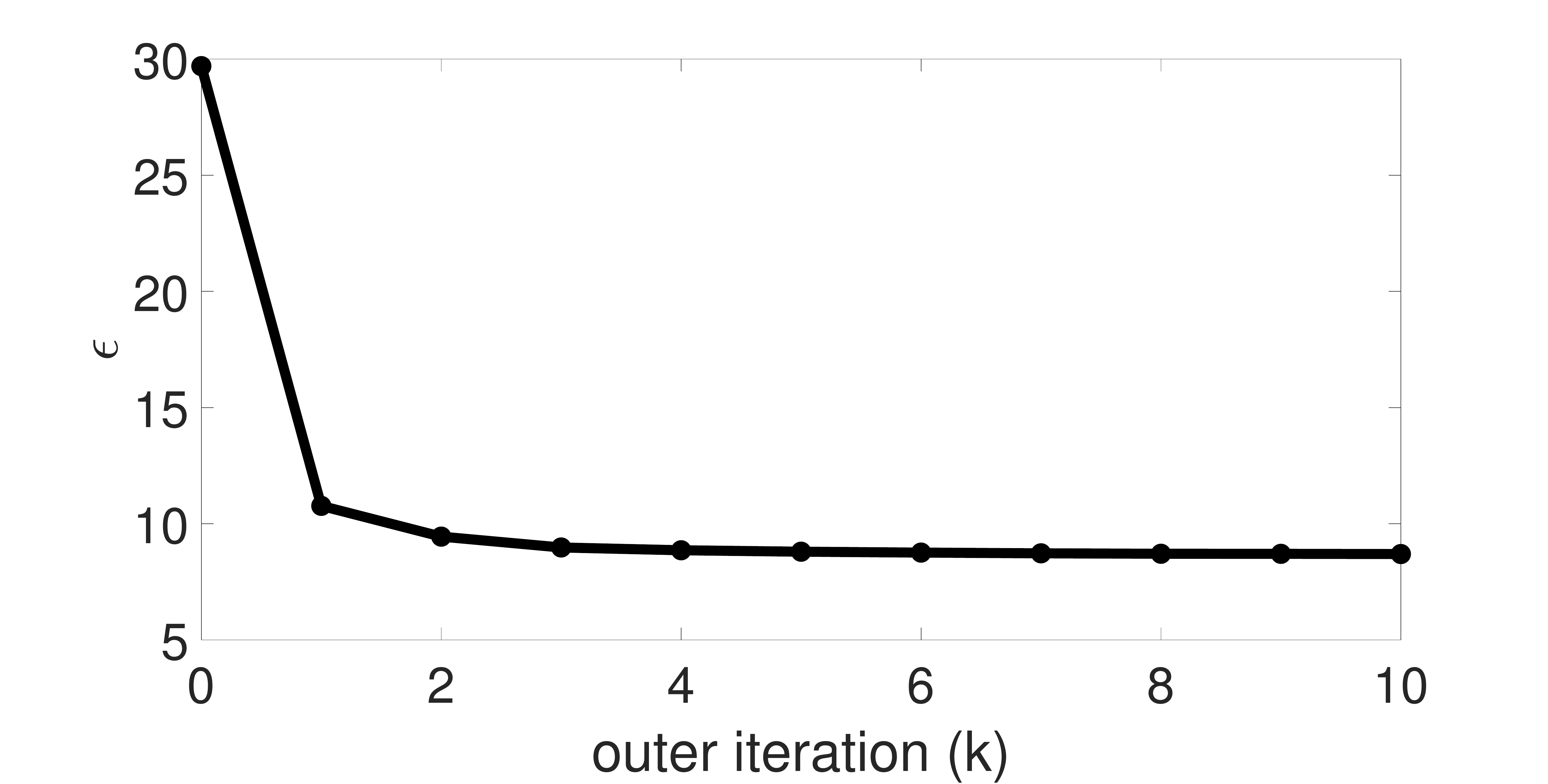}\label{F8a}}
			\subfigure[]{\includegraphics[width=.35\textwidth]{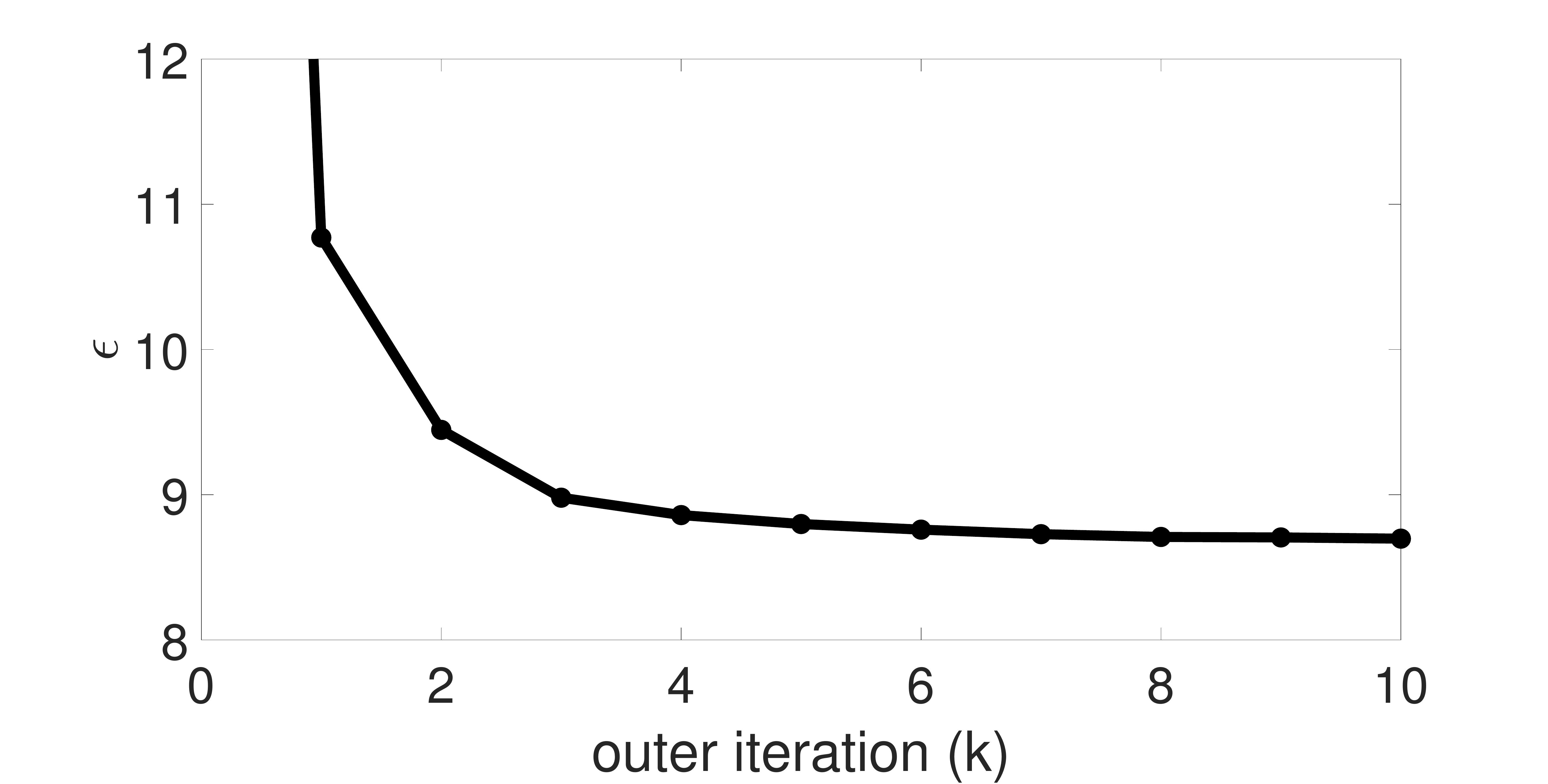}\label{F8b}}
			\subfigure[]{\includegraphics[width=.35\textwidth]{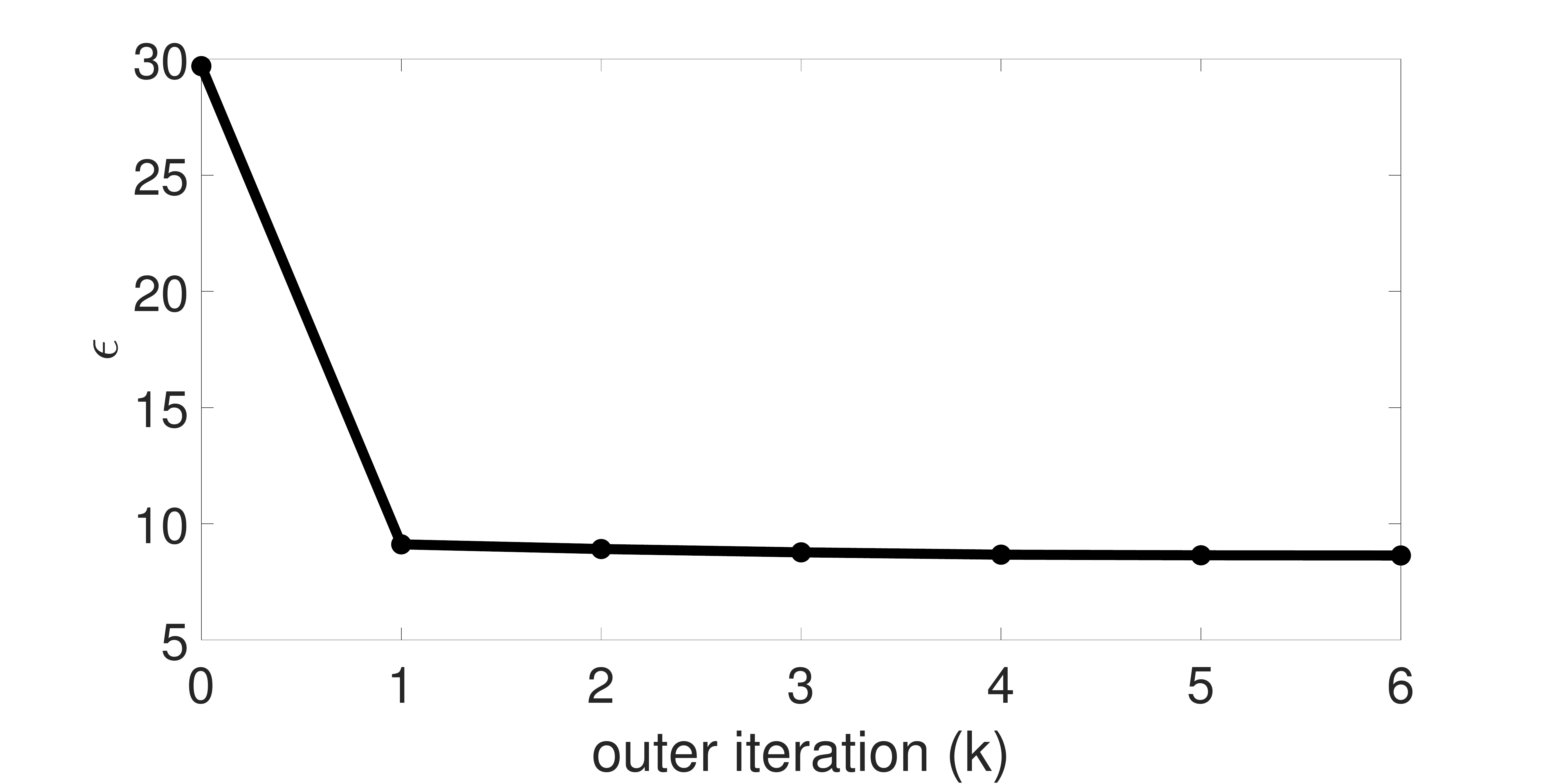}\label{F8c}}
			\subfigure[]{\includegraphics[width=.35\textwidth]{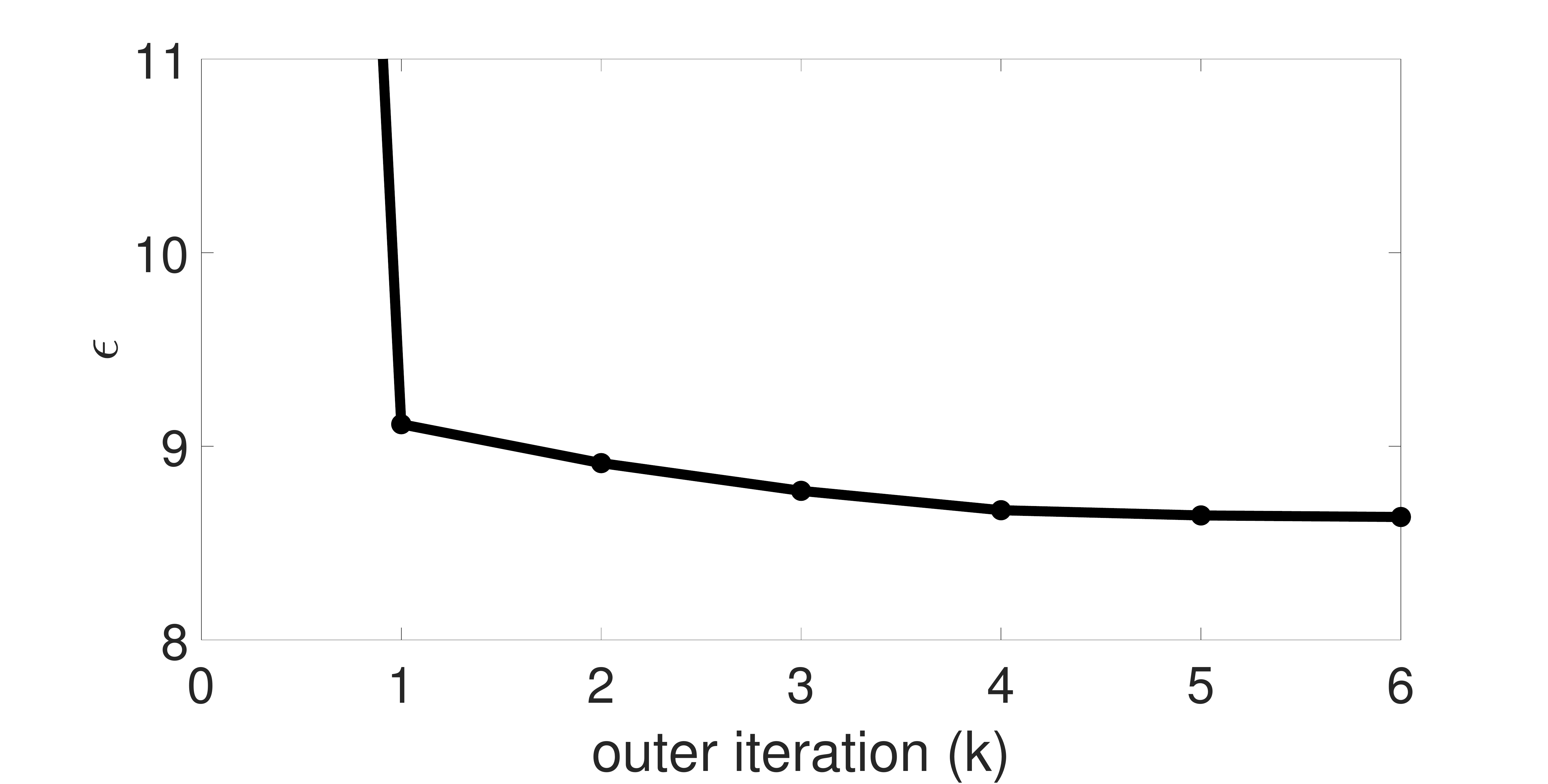}\label{F8d}}
		}
		\caption{$\epsilon$ versus outer iteration $k$ for 3D case. (a) LD  (b) LD from an enlarged view around the optimal point (c) PD-IPM (d) PD-IPM from an enlarged view around the optimal point.}
	\end{figure}

	\section{Discussion}
	
	In the previous section, we numerically evaluated the performance of the used iterative algorithms for a direct problem of QPAT for realistic acoustic media. In this section, we give further details on our numerical results. 
	
	The quality of images reconstructed by the ADMM algorithm was sensitive to choices for $\varrho$ and $Tol_{\text{out}}$, and thus these parameters were chosen very carefully. This leads to several repetition of the entire reconstruction. For exmple, by using a smaller $Tol_{\text{out}}$ (further proceeding of the iterations), the $RE$ values started to increase. One way for avoiding this may be increasing the amount of regularisation via an increase in $\varrho$, but our numerical experience shows that choosing greater values for $\varrho$ negatively affects the convergence of the algorithm. As shown in the second columns in figure~\ref{fig:46}, the reconstructed image for $\mu$ includes a high level of artifact, and also the ADMM algorithm failed to produce an accurate image for $\kappa$. (Note that better images were reconstructed using ADMM for 2D case.) The poor performance of the ADMM algorithm, especially in 3D case, may be because of assuming a high level of errors in the acoustic properties, as shown in figures~\ref{fig:41} and~\ref{fig:45}. 
	
	The LD algorithm was not sensitive to $\gamma$, $Tol_{\text{in}}$ and $Tol_{\text{out}}$. Both $\epsilon$ and $RE$ will proceed with a monotonic reduction, if we choose smaller values for $Tol_{\text{out}}$. This is indicated by figures \ref{F4b} and \ref{F8b}. However, the convergence of the LD algorithm will be deteriorated, if we use very small values for $\beta$, for example the values that we used for the PD-IPM algorithm.
	
	Our numerical experience showed that the PD-IPM algorithm was also not sensitive to choices for $\gamma$, $Tol_{\text{in}}$, $Tol_{\text{med}}$ and $Tol_{\text{out}}$.
	(For example, by choosing $Tol_{\text{med}}=0$, the performance of the algorithm is almost the same.) Additionally, this algorithm converged well using very small values for $\beta$.

	Regarding the computational cost, the ADMM algorithm reconstructed final images using 100-150 gradient-based iterations for both 2D and 3D cases. (Note that each iteration involves at least an implementation of the forward operator and the adjoint of the Jacobian matrix.) Note also that for ADMM, as discussed in the first paragraph of this section, using a smaller stopping threshold deteriorates the quality of reconstructed images. The LD and PD-IPM approaches produced the final images in 150-200 inner iterations. Note that the major cost of each inner iteration for the PCG loop involves an implicit implementation of the Jacobian matrix and its transpose, and has almost the same computational cost as the gradient.
	
	For a direct problem of QPAT for realistic acoustic media with an error in estimation of acoustic properties, our numerical results show that the developed matrix-free Jacobian-based inexact-Newton methods outperform gradient-based approaches that utilise a search direction using Quasi-Newton approaches like the L-BFGS method \cite{Gao,Gao-b}, at the same time does not impose a large computational cost due to an explicit construction of the Jacobian matrix \cite{Pulkkinen-b}.
	
	Our next goal is an extension of multi-source QPAT to multi-spectral QPAT \cite{Bal-d}, which is more practical for biomedical cases, especially when a limited view is accessible for optical excitations. A simultaneous reconstruction of the optical coefficients and the sound speed using adjunct information obtained from ultrasound computed tomography may be promising for improving the quality of reconstructed images \cite{Matthews}.
	
	\section*{Acknowledgement}
	This work was supported in part by a Dean's award from the Faculty of Science and Engineering at the University of Manchester, and the Engineering and Physical Sciences Research Council (EP/M016773/1). There were no data used or collected for this research.
	
	\bibliographystyle{abbrv}
	\bibliography{QPAT_inexact_Newton.bib}

\end{document}